\documentclass[12pt]{article}
\usepackage{amsthm,amsmath,amssymb,amscd,verbatim}
\usepackage{graphics,pifont}
\usepackage{amssymb}
\usepackage{graphicx}
\date{}
\newcommand{\al}{\alpha}
\newcommand{\be}{\beta}
\newcommand{\ov}{\overline}
\newcommand{\sq}{$\square$}
\newcommand{\va}{\varepsilon}
\begin{document}
\title{Continuous Transitive Maps on the Interval Revisited}                       
\author{Bau-Sen Du\footnote {Bau-Sen Du is a retired research fellow at the Institute of Mathematics, Academia Sinica, Taiwan} \\ 
dubs@gate.sinica.edu.tw \cr}
\maketitle

\begin{abstract}
In this note, continuous transitive maps $f$ on the interval $I$ are re-addressed, where $I$ denotes one of the intervals: $(-\infty, \infty)$, $(-\infty, a]$, $[b, \infty)$, $[a, b]$, where $a < b$ are real numbers.  Such maps must have a fixed point, say $z$, in the interior of $I$.  Some well-known properties of such maps are re-proved in a systematic way according to the following : 
\begin{itemize}
\item[(1)] $f$ moves some point $c \ne z$ away from $z$, i.e., fo some point $c \ne z$, we have $f(c) \le c < z$ or $z < c \le f(c)$; 

\item[(2)] $f$ moves some point $\hat c \ne z$ towards but not "over" $z$, i.e., for some point $\hat c \ne z$, we have $\hat c < f(\hat c) < z$ or $z < f(\hat c) < \hat c$; and 

\item[(3)] $f$ moves all points $x \ne z$ to the other side of $z$, i.e., for all points $x \ne z$, we have $x < z \le f(x)$ and $f(x) \le x < z$.  
\end{itemize}
Those maps satisfying (1) are bitransitive, turbulent, and has periodic points of all periods.  Those satisfying (2) are bitransitive with some periodic points of odd periods $\ge 3$, but need not be turbulent while their second iterates are doubly turbulent.  Those satisfying (3) are never bitransitive nor turbulent, have periodic points of all even periods while their second iterates are doubly turbulent.  The proofs are arranged in such ways that they yield the same results.  For example, Theorem 3 treats maps satisfying Condition (1) or Condition (2) while Theorem 4 treats separately maps satisfying Condition (1) and, Conditions (2) or (3).  Characterizations of continuous bitransitive maps on an interval are re-addressed and a new chaotic property of such maps is introduced (Theorem 8, p.21).  For generalizations of Theorem 8 to continuous weakly mixing maps on infinite separable locally compact metric spaces, we refer to {\bf\cite{du5}}.
\end{abstract}

In this note, by an interval, we always mean a non-degenerate interval in the real line.  Let $I$ denote an interval and let $f : I \rightarrow I$ be a continuous map.  For any interval $J$ in $I$, let $int(J)$ and $\ov J$ denote the interior and the closure of $J$ respectively with respct to the topology on the real line.  For each integer $n \ge 1$, let $f^n$ be defined by: $f^1 = f$ and $f^n = f \circ f^{n-1}$ when $n \ge 2$.  For $x_0$ in $I$, we call $x_0$ a periodic point of $f$ with least period $m$ or a period-$m$ point of $f$ if $f^m(x_0) = x_0$ and $f^i(x_0) \ne x_0$ when $0 < i < m$.  If $f(x_0) = x_0$, then we call $x_0$ a fixed point of $f$.  We let $O_f(c) = \{ f^i(c) : i = 0, 1, 2, \cdots \}$ denote the orbit of the point $c \in I$ with respect to $f$.  

We say that $f$ is (topologically) transitive on $I$ if, for any two nonempty open sets $U$ and $V$ in $I$, there exists a positive integer $k$ such that $f^k(U) \cap V \ne \emptyset$, or, equivalently, if there exist a positive integer $k$ and a nonempty open set $W$ in $U$ such that $f^k(W) \subset V$.  We say that $f$ is bitransitive on $I$ if $f^2$ is transitive on $I$.  

When $f$ is a continuous transitive map defined on the interval $I$, the transitivity of $f$ is also equivalent to the existence of a point whose orbit is dense in $I$.  Indeed, if $f$ has a dense orbit in $I$ then it is clear that $f$ is transitive.  Conversely, suppose $f$ is transitive and let $V_1, V_2, \cdots$ be an enumeration of all open intervals in $I$ with rational endpoints.  Let $a < b$ be any two points in the {\it interior} of $I$ and let $U = (a, b)$.  So, the closure $\overline{U} = [a, b]$ is a compact interval in the interior of $I$.  Since $f$ is transitive, there exist a positive integer $n_1$ and an open interval $W_1$ whose length $|W_1|$ is $< 1/2$ such that $\overline{W_1} \subset U$ and  $f^{n_1}(\overline{W_1}) \subset V_1$.  Similarly, there exist a positive integer $n_2 \, (> n_1)$ and an open interval $W_2$ with length $|W_2| < 1/2^2$ such that $\overline{W_2} \subset W_1$ and $f^{n_2}(\overline{W_2}) \subset V_2$.  Inductively, we obtain positive integers $n_1 < n_2 < n_3 < \cdots$ and open intervals $W_i$ with length $|W_i| < 1/2^i, i \ge 1$ such that $U \supset \overline{W_1} \supset W_1 \supset \overline{W_2} \supset W_2 \supset \overline{W_3} \supset W_3 \supset \cdots$ and $f^{n_i}(\overline{W_i}) \subset V_i, i \ge 1$.  Now since $< \overline{W_i} >_{i \ge 1}$ is a nested sequence of nonempty compact intervals whose lengths tend to zero, the set $\bigcap_{i \ge 1} \overline{W_i}$ is nonempty and consists of exactly one point which lies in $\bigcap_{i \ge 1} \overline{W_i}$ ($\subset U$).  It is obvious that this point is a {\it transitive} point of $f$, i.e., it has a dense orbit in $I$. 

Let $J$ be a compact interval in $I$.  If there are two compact intervals $J_0$ and $J_1$ in $J$ with at most one point in common such that $f(J_0) \cap f(J_1) \supset J_0 \cup J_1$, then we say that $f$ is turbulent on $J$ and on $I$ {\bf\cite{bc}}.  If there are two compact intervals $I_0$ and $I_1$ in $I$ with at most one point in common such that $f$ is turbulent on $I_0$ and on $I_1$, then we say that $f$ is doubly turbulent on $I$.  It is well-known that if $f$ is a continuous transitive map on a compact interval then $f$ has dense periodic points {\bf\cite{bc, vb}}, has periodic points of all {\it even} periods {\bf\cite{li}} and $f^2$ is turbulent \big(Corollary VI.4 of {\bf\cite{bc}}, see also {\bf{\cite{du2}}}\big).  On the other hand, continuous transitive maps and continuous maps with dense periodic points are characterized in {\bf\cite{barge1, barge2}} \big(see also {\bf\cite{bl}}\big) by using the notion of snakelike continua introduced in {\bf\cite{bing}}.  From this characterization of continuous maps with dense periodic points {\bf\cite{barge1}}, we can see that {\it nontrivial} continuous maps with dense periodic points are closely related to transitive maps.  In this note, we shall re-prove these results with standard arguments (without resorting to the celebrated Sharkovsky's theorem {\bf\cite{sh}} or the notion of snakelike continua\big) which are more accessible to the readers \big(cf. {\bf\cite{ru}}\big). 

We first note that each continuous transitive map $f$ on the interval $I$ must have a fixed point $z$ in the interior of $I$.  So, $f$ satisfies exactly one of the following 3 mutually exclusive conditions: 
\begin{itemize}
\item[{\rm(1)}]
$f$ moves some point $c \ne z$ away from $z$, i.e., fo some point $c \ne z$, we have $f(c) \le c < z$ or $z < c \le f(c)$; 

\item[{\rm(2)}] 
$f$ moves some point $\hat c \ne z$ towards but not "over" $z$, i.e., for some point $\hat c \ne z$, we have $\hat c < f(\hat c) < z$ or $z < f(\hat c) < \hat c$; and 

\item[{\rm (3)}] $f$ moves all points $x \ne z$ to the other side of $z$, i.e., for all points $x \ne z$, we have $x < z \le f(x)$ and $f(x) \le x < z$. \end{itemize}

Since different conditions sometimes produce same results, for the sake of clarity, those conditions which will yield the same results are treated simultaneously.  For example, Theorem 1 treats separately maps satisfying Condition (1) and, Conditions (2) or (3) while Theorem 4 treats separately maps satisfying Conditions (1) or (2) and, Condition (3).  

\noindent
{\bf Theorem 1.} 
{\it Let $I$ denote one of the following intervals: $(a, b)$, $(a, b]$, $[a, b)$, $[a, b]$, $(-\infty, b)$, $(-\infty, b]$, $(a, \infty)$, $[a, \infty)$ and $(-\infty, \infty)$, where $a < b$ are real numbers, and let $f$ be a continuous transitive map on $I$.  Then exactly one of the following statements hold:  
\begin{itemize}
\item[{\rm (1)}]
If there exist a fixed point $z$ of $f$ and a point $c \ne z$ in $I$ such that $\bigr(f(c)-z\bigr)/(c-z) \ge 1$, then $f$ is turbulent (and so $f^2$ is doubly turbulent) and has periodic points of all periods.

\item[{\rm (2)}]
If $f$ has a unique fixed point $z$ (which lies in the interior of $I$) and $\bigr(f(x) - x\bigl)/(z - x) > 0$ (i.e., both $f(x)$ and $z$ lie on the same side of $x$) for all $x \ne z$ in $I$, then $f^2$ is doubly turbulent and $f$ has periodic points of all even periods. 
\end{itemize}
\noindent
Consequently, if $f$ is a continuous transitive map on $I$, then $f$ has periodic points of all even periods and $f^2$ is doubly turbulent on $I$.}

\noindent
{\it Proof.}
For the proof of (1), assume that there exist a fixed point $z$ of $f$ and a point $c$ in $I$ such that $f(c) \le c < z$ \, \big(if $z < c \le f(c)$, the proof is similar\big).  If $f(c) = c$, then, since $f$ is not constant on $[c, z]$, there exists a point $c^*$ in $(c, z)$ such that $f(c^*) \ne c^*$ and so, either $f(c) = c < c^* < f(c^*)$ or $f(c^*) < c^* < z$.  Without loss of generality, we may assume that $f(c) < c < z$  and assume that there is no fixed points of $f$ in $[c, z)$ and so, $f(x) < x$ for all $c \le x < z$.  We have two cases to consider: 

Case 1. The point $d = \max \big\{ x \in I : x \le c$ and $f(x) = z \big\}$ exists (if $I$ is a compact interval, then this case must hold).  In this case, $$\text{if} \,\,\, f(x) > d \,\,\, \text{for all} \,\,\, d \le x \le z,$$ then $f\big([d, z]\big) \subset (d, z] \subsetneq [d, z]$ and so, $f\bigl(f([d,z])\bigr) \subset f\big([d, z]\big)$.  Thus, $f\bigl([d, z]\bigr)$ is a {\it proper} invariant compact interval of $f$ in $[d. z] \, (\subset I)$ which contradicts the transitivity of $f$.  This implies that there is a point $b$ in $(d, z)$ such that $$f(b) = \min \big\{ f(x) : d \le x \le z \big\} \le d.$$Consequently, $$f\bigl([d, b]\bigr) \cap f\bigl([b, z]\bigr) \supset [d, z] = [d, b] \cup [b, z].$$  Therefore, $f$ is turbulent on $[d, z]$ and clearly $f^2$ is doubly turbulent on $I$.

On the other hand, suppose, for some integer $k \ge 1$, there is a point $q_k$ in $[d, b]$ such that $f^k(q_k) = q_k$.  Then $f^k\bigl([d, q_k]\bigr) \supset [q_k, z] \supset [b,z]$.  So, $f^{k+1}\bigl([d, q_k]\bigr) \supset f\bigl([b, z]\bigr) \supset [d, z] \supset [d, q_k]$.  Consequently, there is a point $q_{k+1}$ such that $d < q_{k+1} < q_k$ and $f^{k+1}(q_{k+1}) = q_{k+1}$.  Since $f\bigl([d, b]\bigr) \supset [d, z] \supset [d, b]$, there is a fixed point of $f$ in $[d, b]$.  It follows from the above that, for each $m \ge 1$, the point $$p_m = \min \big\{ d \le x \le b : f^m(x) = x \big\}$$ is a period-$m$ point of $f$ in $[d, b] \,\,\, \big(\subset [d, z]\big)$ and $d < \cdots < p_3 < p_2 < p_1 < b < z$.  

Case 2.  $f(x) < z$ for all $x$ in $I \cap (-\infty, c)$.  In this case, since we already have $f(x) < x$ for all $c \le x < z$, we actually have $f(x) < z$ for all $x$ in $I \cap (-\infty, z)$.  Consequently, $f\bigr(I \cap (-\infty, z]\bigl) \subset I \cap (-\infty, z]$.  Since $f$ is transitive on $I$, we have $I \cap (-\infty, z] = I$.  

If $I \cap (-\infty, c]$ were compact, then $I$ is bounded below.  Let $w$ be a point in $(c, z)$ such that $w > f(x)$ for all $x$ in $[\min I, c]$.  Then, $f\bigl([\min I, w]\bigr) \subset [\min I, w) \subsetneq [\min I, w]$ and so, $f\bigl([\min I, w]\bigr)$ is a {\it proper} invariant compact interval of $f$ in $[\min I, w]$ $(\subset I)$ which contradicts the transitivity of $f$.  So, in this case, $I \cap (-\infty, c]$ cannot be compact.  Therefore, we have shown that 
$$
\text{if} \,\,\, f(x) < z \,\,\, \text{for all} \,\,\, x \in I \cap (-\infty, c], \, \text{then} \,\,\, I = (-\infty, z] \,\,\, \text{or} \,\,\, I = (a, z] \,\,\, \text{for some point} \,\,\, a < z.
$$
\indent The following proof works for $I = (-\infty, z]$ and for $I = (a, z]$ for some point $a < z$.  

Since $f(c) < c < z$, if $f(x) < x$ for all $x$ in $I \cap (-\infty, c]$, then $f\bigl(I \cap (-\infty, c]\bigr) \subset I \cap (-\infty, c]$ which contradicts the transitivity of $f$ on $I$.  So, $(z >)$ $f(c') \ge c'$ for some $c'$ in $I \cap (-\infty, c)$. 

If ($z >) \, f(x) > x$ on $I \cap (-\infty, c']$, let $\delta$ be a point in $I \cap (-\infty, c')$ such that $$\delta < \min \big\{ f(x) : c' \le x \le z \big\}.$$  Then, since $\delta < f(x) < z$ on $[c', z)$ and $f(x) > x$ on $[\delta, c']$, we have$$\delta < f(x) \le z \,\,\, \text{for all} \,\,\, \delta \le x \le z.$$ This implies that $f\bigl([\delta, z]\bigr) \subset (\delta, z] \subsetneq [\delta, z]$ and so, $f\bigl([\delta, z]\bigr)$ is a {\it proper} invariant compact interval of $f$ in $I$ which contradicts the transitivity of $f$ on $I$.  Consequently, there is a point $\hat c_1$ in $I$ such that $f(\hat c_1) < \hat c_1 < c' < z$.  Let $$\hat z_1 = \min \big\{ \hat c_1 \le x \le c' : f(x) = x \big\}.$$  If $f(x) \le \hat z_1$ for all $x \in I \cap (-\infty, \hat z_1]$, then $I \cap (-\infty, \hat z_1]$ is a {\it proper} invariant interval of $f$ in $I$ which contradicts the transitivity of $f$ on $I$.  So, there is a point $\hat t_1 < \hat c_1$ such that $f(\hat t_1) > \hat z_1$.  Consequently, the point $$\hat d_1 = \max \big\{ \hat t_1 \le x \le \hat c_1 : f(x) = \hat z_1 \big\}$$ exists.  If $f(x) > \hat d_1$ for all $\hat d_1 \le x \le \hat z_1$, then $f\bigl([\hat d_1, \hat z_1]\bigr) \subset (\hat d_1, \hat z_1] \subsetneq [\hat d_1, \hat z_1]$ which contradicts the transitivity of $f$ on $I$.  

Therefore, $\min \big\{ f(x): \hat d_1 \le x \le \hat z_1 \big\} \le d_1$ which implies that $f$ is turbulent on $[\hat d_1, \hat z_1]$ and, by arguing as those in Case 1, has periodic points of all periods in $[\hat d_1, \hat z_1]$.  

We note in passing that, in this case, by repeating the above arguments indefinitely, we actually have sequences of points $$\cdots < \hat t_k < \hat d_k < \hat c_k < \hat z_k < \cdots < \hat t_2 < \hat d_2 < \hat c_2 < \hat z_2 < \cdots < \hat t_1 < \hat d_1 < \hat c_1 < \hat z_1 < c < z$$such that, for each $k \ge 1$, $$\hat z_k < f(\hat t_k) < z, \,\,\, f(\hat d_k) = f(\hat z_k) = \hat z_k, \,\,\, f(\hat c_k) < \hat c_k < \hat z_k$$ and, $f$ is turbulent on $[\hat d_k, \hat z_k]$ and has periodic points of all periods in $[\hat d_k, \hat z_k]$. 

By combining the above Case 1 and Case 2, we confirm (1).   
$$--------------------------$$
\indent As for the proof of (2), assume that
$$
z \,\,\, \text{is an interior point of} \,\,\, I \,\,\, \text{such that} \,\,\, x < f(x) \,\,\, \text{for all} \,\,\, x < z \,\,\, \text{and} \,\,\, f(x) < x \,\,\, \text{for all} \,\,\, z < x.\qquad (*)
$$
\indent Let $\breve u$ be a point in $I \cap (-\infty, z)$ with {\it dense orbit} in $I$ and let $\breve \beta$ be a point in $I$ such that 
$$
\breve \beta = \max \big\{ f(x) : \breve u \le x \le z \big\}.
$$  
\big(On the interval $I \cap (z, \infty)$, we can do and argue similarly to obtain similar results\big). \\ If $\breve \beta \le z$, then $f\bigl([\breve u, z]\big) \subset (\breve u, z] \subsetneq [\breve u, z]$ and $f\big([\breve u, z]\big)$ is a {\it proper} compact invariant interval of $f$ in $[\breve u, z]$ which contradicts the transitivity of $f$ on $I$.  So, $\breve \beta > z$.  If $\min \big\{ f(x) : z \le x \le \breve \beta \big\} > \breve u$, then by (*), $f\big([\breve u, \breve \beta]\big) \subset (\breve u, \breve \beta] \subsetneq [\breve u, \breve \beta]$.  Consequently, $f\big([\breve u, \breve \beta]\big)$ is a {\it proper} invariant compact interval of $f$ in $[\breve u, \breve \beta]$ $(\subset I)$.  This, again, contradicts the transitivity of $f$ on $I$.  Therefore, we have $\min \big\{ f(x) : z \le x \le \breve \beta \big\} \le \breve u$.  

Let $\breve v$ be a point in $[\breve u, z)$ such that $f(\breve v) \in (z, \breve \beta]$ and $f^2(\breve v) = \breve u < \breve v$.  Let $$\breve z = \min \big\{ \breve v \le x \le z : f^2(x) = x \big\}.$$  Then, since $f(\breve v) > z$ and $f(z) = z$, we have $f(x) > z$ for all $\breve v \le x \le \breve z$ and, since 
$f^2(\breve v) < \breve v < \breve z \le z$, we have $f^2(x) < x < \breve z \le z$.  This implies that 
$$
f^2(x) < x < \breve z \le z < f(x) \,\,\, \text{for all} \,\,\, \breve v \le x < \breve z.
$$  
In particular, $\breve u = f^2(\breve v) < \breve v < \breve z \le z < f(\breve v)$.  We have two cases to consider:

Case I.  The point $\breve d = \max \big\{ x \in I \cap (-\infty, \breve v] : f^2(x) = \breve z \big\}$ exists which, combined with the fact that $f^2(x) < x < \breve z \le z < f(x)$ for all $\breve v \le x < \breve z$, implies that 
$$
f^2(x) < x < \breve z \le z < f(x) \,\,\, \text{for all} \,\,\, \breve d < x < \breve z.
$$  
If $\breve d < f^2(x) \le \breve z \,\,\, \text{for all} \,\,\, \breve d \le x \le \breve z$, 
then, since $f(x) > z \ge \breve z$ for all $\breve d \le x < \breve z$, we have
$$
f\Big([\breve d, \breve z] \cup f\big([\breve d, \breve z]\big)\Big) \,\,= f\big([\breve d, \breve z]\big) \cup f^2\bigl([\breve d, \breve z]\big) \subset f\big([\breve d, \breve z]\big) \cup (\breve d, \breve z] \subsetneq f\big([\breve d, \breve z]\big) \cup [\breve d, \breve z].
$$  
Thus, $f\big([\breve d, \breve z] \cup f([\breve d, \tilde z])\big)$ is a union of two compact intervals which is a {\it proper} invariant subset of $f$ in $I$.  This contradicts the transitivity of $f$ on $I$.  So, $$\text{there exists a point} \,\,\, \breve b \,\,\, \text{in} \,\,\, (\breve d, \breve z) \,\,\, \text{such that} \,\,\, f^2(\breve b) \le \breve d.$$  
Thus, $f^2\bigl([\breve d, \breve b]\bigr) \cap f^2\bigl([\breve b, \breve z]\bigr) \supset [\breve d, \breve z] = [\breve d, \breve b] \cup [\breve b, \breve z]$ which implies that $f^2$ is turbulent on $[\breve d, \breve z]$.  

We now show that, since $f(x) > z \ge \breve z$ for all $\breve d < x < \breve z$, for each $k \ge 1$, the point 
$$
\breve p_{2k} = \min \big\{ \breve d \le x \le \breve b : f^{2k}(x) = x \big\} = \min \big\{ \breve d \le x \le \breve z : f^{2k}(x) = x \big\}
$$ 
is a period-$(2k)$ point of $f$ in $[\breve d, \breve z] \subset I \cap (-\infty, z]$ by arguing as follows:

Note that $f^{2n}(\breve d) = \breve z$ for all $n \ge 1$.  Let $\breve u_2 = \min \big\{ \breve d \le x \le \breve b: f^2(x) = \breve d \big\}$.  Then 
$$
\breve d < f^2(x) < \breve z \,\,\, \text{for all} \,\,\, \breve d < x < \breve u_2.
$$  
\indent Let $\breve p_2 = \min \big\{ \breve d \le x \le \breve u_2: f^2(x) = x \big\}$ and let $\breve \nu_2$ be a point in $(\breve d, \breve p_2)$ such that $f^2(\breve \nu_2) = \breve u_2$.  Then since  $f^4\big([\breve d, \breve u_2]\big) \supset f^4\big([\breve d, \breve \nu_2]\big) \supset f^2\bigl([\breve u_2, \breve z]\bigr) \supset [\breve d, \breve z]$, the point $\breve u_4 = \min \big\{ \breve d \le x \le \breve u_2: f^4(x) = \breve d \big\}$ exists and  
$$
\breve d < f^4(x) < \breve z \,\,\, \text{for all} \,\,\, \breve d < x < \breve u_4.
$$ 
\indent Let $\breve p_4 = \min \big\{ \breve d \le x \le \breve u_4: f^4(x) = x \big\}$ and let $\breve \nu_4$ be a point in $(\breve d, \breve p_4)$ such that $f^4(\breve \nu_4) = \breve u_2$.  Then since  $f^6\bigl([\breve d, \breve u_4]\bigr) \supset f^6\bigl([\breve d, \breve \nu_4]\bigr) \supset f^2\big([\breve u_2, \breve z]\big) \supset [\breve d, \breve z]$, the point $\breve u_6 = \min \big\{ \breve d \le x \le \breve u_4: f^6(x) = \breve d \big\}$ exists and 
$$
\breve d < f^6(x) < \breve z \,\,\, \text{for all} \,\,\, \breve d < x < \breve u_6.
$$
\indent By proceeding in this manner indefinitely, we obtain 2 strictly decreasing sequences 
$$
\breve d < \cdots < \breve p_6 < \breve u_6 < \breve p_4 < \breve u_4 < \breve p_2 < \breve u_2 < \breve v < \breve z
$$ 
such that, for each $k \ge 1$, (let $\breve u_0 = \breve b$) 
$$
\breve u_{2k} = \min \big\{ \breve d \le x \le \breve u_{2k-2}: f^{2k}(x) = \breve d \big\},
$$
$$\breve p_{2k} = \min \big\{ \breve d \le x \le \breve u_{2k}: f^{2k}(x) = x \big\} \,\,\, \text{and}
$$ 
$$
\breve d < f^{2n}(x) < \breve z \,\,\, \text{for all} \,\,\, \breve d < x < \breve u_{2k} \,\,\, \text{and all} \,\,\, 1 \le n \le k.
$$
In particular, since $f(x) > z \ge \breve z$ for all $\breve d < z < \breve z$, we obtain that 
$$
\breve d < f^i(p_{2k}) < \breve z \le z < f^j(p_{2k}) \,\,\, \text{for all even} \,\,\, i \,\,\, \text{and all odd} \,\,\, j \,\,\, \text{in} \,\,\, [0, 2k].
$$
Therefore, the least period of $p_{2k}$ with respect to $f$ is {\it even} and divides $2k$.  Since 
\begin{multline*}
$$
\cdots < \breve p_{2k} < \breve p_{2k-2} < \cdots < \breve p_4 < \breve p_2 < \breve z \,\,\, \text{and} \\ 
\text{for each} \,\,\, k > i \ge 1, \, \breve p_{2i} \,\,\, \text{is the smallest point in} \,\,\, [\breve d, z] \,\,\, \text{such that} \,\,\, f^{2i}(x) = x, 
$$ 
\end{multline*}
we obtain that, for each $k \ge 1$, the point $\breve p_{2k}$ is a period-$(2k)$-point of $f$ in $[\breve d, \breve z]$. \vspace{.1in}

Case II.  $f^2(x) < \breve z$ for all $x$ in $I \cap (-\infty, \breve v]$.  In this case, since we already have $$f^2(x) < x < \breve z \le z < f(x) \,\,\, \text{for all} \,\,\, \breve v < x < \breve z,$$we actually have $$\qquad\qquad\qquad\qquad\qquad\quad\,\, f^2(x) < \breve z \le z < f(x) \,\,\, \text{for all} \,\,\, x \,\,\, \text{in} \,\,\, I \cap (-\infty, \breve z]. \qquad\qquad\qquad\qquad\,\,\, (**)$$
\noindent
In particular, $f^2\big(I \cap (-\infty, \breve z]\big) \subset I \cap (-\infty, \breve z]$.  Since 
$$
f\Big(I \cap (-\infty, \breve z] \cup f\big(I \cap (-\infty, \breve z]\big)\Bigr) = f\big(I \cap (-\infty, \breve z]\big) \cup f^2\big(I \cap (-\infty, \breve z]\big) \subset f\big(I \cap (-\infty, \breve z]\big) \cup I \cap (-\infty, \breve z],
$$ 
the set $I \cap (-\infty, \breve z] \cup f\big(I \cap (-\infty, \breve z]\big)$ is a union of two intervals which is invariant under $f$.  Since $f$ is transitive on $I$, we have $I \cap (-\infty, \breve z] \cup f\big(I \cap (-\infty, \breve z]\big) \,= I$ and so, by (**), we must have $$
\breve z = z \,\,\, \text{and} \,\,\, f^2(x) < z < f(x) \,\,\, \text{for all} \,\,\, x \,\,\, \text{in} \,\,\, I \cap (-\infty, z).
$$
\noindent
\indent If $I \cap (-\infty, \breve v]$ were compact, then $I$ is bounded below.  Let $\breve w$ be a point in $(\breve v, z)$ such that $\breve w > f(x)$ for all $x$ in $[\min I, \breve v]$.  Then, the set $[\min I, \breve w] \cup f\big([\min I, \breve w]\big)$ is a union of two compact intervals such that 
$$
f\Big([\min I, w] \cup f\big([\min I, w]\big)\Big) = f\big([\min I, w]\big) \cup f^2\big([\min I, w]\big) \subsetneq f\big([\min I, w]\big) \cup [\min I, w)
$$
which, by (**), contradicts the transitivity of $f$ on $I$.  So, in this case, $I \cap (-\infty, \breve v]$ cannot be compact.  Therefore, we have shown that 
\begin{multline*}
$$
\text{if} \,\,\, f^2(x) < \breve z = z \,\,\, \text{for all} \,\,\, x \in I \cap (-\infty, \breve v], \, \text{then} \,\, I = (-\infty, \infty) \,\,\, \text{or}\, \\ \, I = (a, b), \, \text{where} \,\,\, a < b \,\,\, \text{are real numbers.\qquad\,\,\,\,}
$$
\end{multline*}
Furthermore, by (**), we have 
$$
f\big(I \cap (-\infty, z]\big) \subset I \cap [z, \infty), \, f\big(I \cap [z, \infty)\big) \subset (-\infty, z] \,\,\, \text{and} \,\,\, f^2\big(I \cap (-\infty, z)\big) \subset I \cap (-\infty, z).
$$  
By arguing as those in Case 2 of the proof of (1) above, we obtain sequences of points 
$$
\cdots < \breve t_k < \breve d_k < \breve v_k < \breve z_k < \cdots < \breve t_2 < \breve d_2 < \breve v_2 < \breve z_2 < \cdots < \breve t_1 < \breve d_1 < \tilde v_1 < \breve z_1 < \breve z = z
$$
such that, for each $k \ge 1$, 
$$
z > f^2(\breve t_k) > \breve z_k, \, f^2(\breve d_k) = f^2(\breve z_k) = \breve z_k, \, f^2(\breve v_k) < \breve v_k < \breve z_k,
$$ 
$f^2$ is turbulent on $[\breve d_k, \breve z_k]$ and, by arguing as those in Case I above, $f$ has periodic points of all even periods in $[\breve d_k, \breve z_k] \cup f\bigl([\breve d_k, \breve z_k]\bigr)$.  

Finally, similar arguments as those on the interval $I \cap (-\infty, z]$ show that $f^2$ is turbulent on $I \cap [z, \infty)$ and so $f^2$ is {\it doubly} turbulent on $I$.  

By combining the above Case I and Case II, we establish (2). 
\hfill\sq

We need two preliminary results.

\noindent
{\bf Lemma 2.}
{\it Let $I$ denote one of the following intervals: $(a, b)$, $(a, b]$, $[a, b)$, $[a, b]$, $(-\infty, b)$, $(-\infty, b]$, $(a, \infty)$, $[a, \infty)$ and $(-\infty, \infty)$, where $a < b$ are real numbers.  Let $f$ be a continuous transitive map on $I$ and let $z$ be a fixed point of $f$ in the interior of $I$.  If there exists a point $c \ne z$ in $I$ such that $\big(f(c) - z\big)/(c - z) > 0$, then, for every open interval $V$ in $I$, we have $\bigcup_{i \ge 0} f^i(V) \supset int(I)$.  In particular, there is an integer $j \ge 0$ such that $z \in f^j(V)$.}

\noindent
{\it Proof.}
Let $V$ be an open interval in $I$.  Let $G$ be the {\it maximal} connected set in $\bigcup_{i \ge 0} f^i(V)$ which contains $V$.  Then $G$ is an {\it interval}.  Note that if $f^k(G) \cap G \ne \emptyset$ for some integer $k \ge 1$, then, since $f^k(G) \cup G$ is an interval in $\big(\bigcup_{i \ge 0} f^{k+i}(V)\big) \, \bigcup \, \big(\bigcup_{i \ge 0} \, f^i(V)\big) = \bigcup_{i \ge 0} f^i(V)$ which contains $V$, by the maximality of $G$, we have $f^k(G) \subset f^k(G) \cup G \subset G$.

Now since $f$ is transitive on $I$, $\emptyset \ne  f^m\big(int(G)\big) \cap int(G) \subset f^m(G) \cap G$ for some $m \ge 1$.  So, $f^m({G}) \subset {G}$.  Since $f$ is transitive on $I$, the finite union $\bigcup_{i=0}^{m-1} f^i(G)$ $\bigl(= \bigcup_{i\ge0} f^i(G)\bigr)$ is dense in $I$.  Since $z$ is a fixed point of $f$ which lies in the interior of $I$, we have, for some integer $0 \le j < m$,$$\text{either} \,\,\, z \in f^j(G) \,\,\, \text{or} \,\,\, z \,\,\, \text{is an endpoint of} \,\,\, f^j(G).$$In either case, there exists a sequence $<x_n>_{n \ge 1}$ of points in $G$ such that $\lim_{n \to \infty} f^j(x_n) = z$.  Consequently, $\lim_{n \to \infty} f^m(x_n) = z$.  Since $f^m(G) \subset G$, we obtain that 
$$
\text{either} \,\,\, z \in G \,\,\, \text{or} \,\, z \,\,\, \text{is an endpoint of} \,\,\, G.
$$
\indent Since $c$ is a point such that $\big(f(c) - z\big)/(c - z) > 0$, without loss of generality, we may assume that $c$ lies in the {\it interior} of $I$ and both $c$ and $f(c)$ lie on the right side of $z$, i.e., we may assume that $z < \min \big\{ c, f(c) \big\}$.  Let $u$ be a point in $G$ with dense orbit in $I$ and let $\ell$ be a positive integer such that $z < c < f^\ell(u)$.  Then $f^\ell(G) \supset \big(z, f^\ell(u)\big] \supset (z, c]$.  Since $z$ is a fixed point of $f$ and
\begin{multline*}
$$
\text{both} \,\,\, c \,\,\, \text{and} \,\,\, f(c) \,\,\, \text{lie on the right side of} \,\,\, z, \, \text{we obtain that} \\
\emptyset \ne \big(z, \min \{ c, f(c) \}\big) \subset f^{\ell+1}(G) \cap f^\ell(G). 
$$
\end{multline*}
Let $n$ be a positive integer such that $nm - \ell > 0$.  Then, since $f^m(G) \subset G$, we have 
$$
\emptyset \ne f^{nm-\ell}\bigg(\big(z, \min \{ c, f(c) \}\big)\bigg) \subset f^{nm-\ell}\bigg(f^{\ell+1}(G) \cap f^\ell(G)\bigg) \subset f^{nm+1}(G) \cap f^{nm}(G) \subset f(G) \cap G. 
$$
Consequently, $f(G) \subset G$.  Since $f$ is transitive and $f(G) \subset G$, $\bigcup_{i \ge 0} f^i(G) \, (= G)$ is dense in $I$.  Since $G$ is an interval, this implies that $G = \bigcup_{i \ge 0} f^i(G) \supset int(I)$ and so, $\bigcup_{i \ge 0} f^i(V) \supset G \supset int(I)$.    
\hfill\sq

\noindent
{\bf Lemma 3.} 
{\it Let $I$ denote one of the following intervals: $(a, b)$, $(a, b]$, $[a, b)$, $[a, b]$, $(-\infty, b)$, $(-\infty, b]$, $(a, \infty)$, $[a, \infty)$ and $(-\infty, \infty)$, where $a < b$ are real numbers, and let $f$ be a continuous map on $I$.  Let $J$ be an interval in $I$ which contains a fixed point $z$ of $f$ and a periodic point $q$ of $f$ with least period $m \ge 3$.  Then at least one of the following holds:
\begin{itemize}
\item[{\rm (1)}]
$f^n(J) \supset O_f(q)$ for all $n \ge 2m-2$;

\item[{\rm (2)}]
$m$ is even, the sets $O_{f^2}(q)$ and $f\big(O_{f^2}(q)\big)$ lie on opposites sides of $z$ and $f^{2n}(J) \supset O_{f^2}(q)$ for all $n \ge m/2 - 1$.
\end{itemize}
\noindent
Consequently, if $J$ contains a fixed point of $f$ and $f^s(J) \subset J$ for some positive integer $s$, then $f^2(P \cap J) \subset J$, where $P$ is the set of all periodic points of $f$ in $I$.}
\noindent

\noindent
{\bf Remark 1.} 
In the above lemma, if we only require $J$ to contain at least two points of a period-$m$ orbit of $f$, then the conclusion need not hold.  For example, let $f : [1, 9] \to [1, 9]$ be the continuous map defined by putting 
$$
f(x) = \begin{cases}
       x + 3 & \text{for $1 \le x \le 6$}, \\
       -7x + 51 & \text{for $6 \le x \le 7$}, \\
       x - 5 & \text{for $7 \le x \le 8$, and} \\ 
       -2x + 19 & \text{for $8 \le x \le 9$} \\
       \end{cases}
$$and let $J = [1, 3]$.  Then it is easy to see that the set $\{ 1, 2, 3, 4, 5, 6, 7, 8, 9 \}$ forms a period-9 orbit of $f$ and $f\big([1, 3]\big) = [4, 6], f\big([4, 6]\big) = [7, 9]$ and $f\big([7, 9]\big) = [1, 3]$.  Therefore, the conclusion of Lemma 3 does not hold.  

\noindent
{\it Proof.}
Let $a$ and $b$ be two {\it adjacent} points in the orbit $O_f(q)$ such that $f(b) \le a < b \le f(a)$ and let $i$ and $j$ be integers in $[0, m-1]$ such that $f^i(q) = a$ and $f^j(q) = b$.  If, for some integer $0 \le k \le m-1$, both $f^k(q)$ and $f^{k+1}(q)$ lie {\it on the same side} of $z$, then we have four cases to consider: $$\text{(i)} \,\, f^i(q) = a < b < z; \,\,\, \text{(ii)} \,\, f^k(q) < a < z < b; \,\,\, \text{(iii)} \,\, a < z < b < f^k(q); \, \text{and} \,\,\, \text{(iv)} \,\, z < a < b = f^j(q).$$In Case (i), we have $f^i(J) \supset [a, b]$.  In Case (iv), we have $f^j(J) \supset [a, b]$.  In Cases (ii) $\&$ (iii), since $f^k(J)$ contains the compact interval with $z$ and $f^k(q)$ as endpoints, we have $f^{k+1}(J) = f\big(f^k(J)\big) \supset [a, b]$.  In either case, since $f\big([a, b]\big) \supset [a, b]$, we have $f^m(J) \supset [a, b]$ and, since $f^{m-2}\big([a, b]\big) \supset O_f(q)$, we have $f^n(J) \supset O_f(q)$ for all $n \ge 2m-2$.  

On the other hand, if, for all $0 \le i \le m-1$, $f^i(q)$ and $f^{i+1}(q)$ lie on opposite sides of $z$, then $m$ is even.  Let $0 \le r \le m-1$ be the {\it even} integer such that $f^r(q)$ is the farthest point from $z$ among those points of $O_f(q)$ which lie on the same side of $z$ as $q$.  Then $f^{r+2j}(J) \supset \big\{ f^2(q), f^4(q), \cdots, f^{m-2}(q), f^{m}(q) = q  \big\}$ for all integers $j \ge 0$.  Therefore, $f^{2n}(J) \supset O_{f^2}(q)$ for all $n \ge (m/2) - 1$.  
\hfill\sq

Now we use Lemmas 2 $\&$ 3 to prove the following result of Barge and Martin {\bf\cite{barge1, barge3}} which implies that if $f$ is transitive then $f$ has dense periodic points {\bf\cite{bc, vb}}.

\noindent
{\bf Theorem 4.}
{\it Let $I$ denote one of the following intervals: $(a, b)$, $(a, b]$, $[a, b)$, $[a, b]$, $(-\infty, b)$, $(-\infty, b]$, $(a, \infty)$, $[a, \infty)$ and $(-\infty, \infty)$, where $a < b$ are real numbers, and let $f$ be a continuous transitive map on $I$, or, equivalently, assume that there is a point $u$ in $I$ whose orbit $O_f(u)$ is dense in $I$.  Then exactly one of the following holds:
\begin{itemize}
\item[{\rm (1)}]
If, for some fixed point $z$ of $f$ and some point $c \ne z$, we have $\big(f(c) - z\big)/(c - z) > 0$, then, for any compact interval $K$ in $I$ and any compact interval $L$ in the interior of $I$, there exists a positive integer $N$ such that $f^n(K) \supset L$ for all $n \ge N$.  

\item[{\rm (2)}]
There exists a unique fixed point $z$ of $f$ in the interior of $I$ such that $\big(f(x) - z\big)/(x - z) \le 0$ for all $x \ne z$ in $I$ and, on each of the intervals $I \cap (-\infty, z]$ and $I \cap [z, \infty)$, $f^2$ is transitive and has at least two fixed points.  Consequently, in this case \big(i.e., $f$ is transitive but not bitransitive on $I$\big), $I$ must be a compact interval or an open interval in the real line.
\end{itemize}
Consequently, if $f$ is transitive on $I$, then $f$ has dense periodic points in $I$.}

\noindent
{\bf Remark 2.}
For every integer $n \ge 2$, let $f_n(x)$ be the classical continuous transitive map {\bf{\cite{barge2, bc}}} from $[1, 2n+1]$ onto itself defined by putting $f_n(1) = n+1$; $f_n(i) = 2n+3-i$ for $2 \le i \le n+1$; $f_n(j) = 2n+2-j$ for $n+2 \le j \le 2n+1$; and by linearity on each interval $[k, k+1]$, $1 \le k \le 2n$ which has exactly one period-$(2n+1)$ orbit but no period-$(2n-1)$ orbits. Note that each $f_n$ is also bitransitive, but none is turbulent.  

\noindent
{\it Proof.}
Suppose there is no point $c \ne z$ such that $z < \min \big\{ c, f(c) \big\}$ or $\max \big\{ c, f(c) \big\} < z$.  Then $$f\big(I \cap (-\infty, z]\big) \subset I \cap [z, \infty) \,\,\, \text{and} \,\,\, f\big(I \cap [z, \infty)\big) \subset I \cap (-\infty, z].$$ So, $$f^2\big(I \cap (-\infty, z]\big) \subset I \cap (-\infty, z] \,\,\, \text{and} \,\,\, f^2\big(I \cap [z, \infty)\big) \subset I \cap [z, \infty).$$  Consequently, since $f$ has a dense orbit in $I$, $f^2$ has a dense orbit in $I \cap (-\infty, z]$ and in $I \cap [z, \infty)$.  Since $f^2$ has a dense orbit in $I \cap (-\infty, z]$ \big($I \cap [z, \infty)$ respectively\big), it has a fixed point $\hat z$ in the interior of $I \cap (-\infty, z)$ \big($I \cap [z, \infty)$ respectively\big).  So, we are back to (1) with $f^2$ replacing $f$ and $\hat z$ replacing $c$.  Therefore, it suffices to prove just (1).  

Now we prove (1).  Let $z$ be a fixed point of $f$ and let $c$ be a point such that $z < \min \big\{ c, f(c) \big\}$ or $\max \{ c, f(c) \} < z$.  If $z$ is not an interior point of $I$, then since $f$ is transitive, $f$ must have a fixed point $z^*$ in the interior of $I$.  In this case, we consider $z^* < z$ or $z < z^*$ instead of $z < \min\{ c, f(c) \}$ or $\max\big\{ c, f(c) \big\} < z$.  So, we may assume that $z$ is a fixed point of $f$ in the {\it interior} of $I$.  We may also assume that $c$ is a point in the {\it interior} of $I$.  

Let $K$ be a compact subinterval of $I$ and let $L$ be a compact subinterval of $int(I)$.  By considering a subinterval of $int(K)$ if necessary, we may assume that$$K \,\,\, \text{lies in the interior of} \,\,\, I \,\,\, \text{and does not contain a fixed point of} \,\,\, f.$$By Lemma 2, there is a point $v$ in the interior of $K$ and a positive integer $n_1$ such that $f^{n_1}(v) = z$.  Since $int(K)$ contains a point $u$ with dense orbit, ther exist positive integers $(n_1 <)$ $n_2 < n_3 < n_4$ such that the point $f^{n_2}(u)$ is so close to the point $c$ that $\big(f^{n_2+1}(u) - z\big)/\big(f^{n_2}(u) - z\big) > 0$, i.e, both the points $f^{n_2}(u)$ and $f^{n_2+1}(u)$ lie on the same side of $z$, and$$f^{n_3}(u) < \min \big\{\min K, \min L, z \big\} < \max \big\{\max K, \max L, z \big\} < f^{n_4}(u).$$If $\max K < z$, then $f^{n_3}(K) \supset \big\{ f^{n_3}(u), z \big\}$.  If $z < \min K$, then $f^{n_4}(K) \supset \big\{ f^{n_4}(u), z \big\}$.  In either case, $f^i(K) \supset K \cup \big\{ z \big\}$ for $i = n_3$ or $n_4$ depending on the relative locations of $K$ with respect to $z$.  In particular, since $f^i(K) \supset K$, $K$ contains a periodic point of $f$.  Since $K$ is arbitrary in $I$, we have actually shown that the periodic points of $f$ are dense in $I$.  

Therefore, we can choose a period-$j$ point $q$ of $f$ in $K$ which is so close to the point $u$ that
$$
\frac {f^{n_2+1}(q) - z}{f^{n_2}(q) - z} > 0, \, \text{i.e., both points} \,\,\, f^{n_2+1}(q) \,\,\, \text{and} \,\,\, f^{n_2}(q) \,\,\, \text{lie on the same side of} \,\,\, z, \, \text{and}
$$
$$
f^{n_3}(u) \approx f^{n_3}(q) < \min \big\{\min K, \min L, z \big\} < \max \big\{\max K, \max L, z \big\} < f^{n_4}(q) \approx f^{n_4}(u).\vspace{.1in}
$$ 
\indent Consequently, since the compact interval $f^{n_2}(K)$ contains the period-$j$ point $f^{n_2}(q)$ and the fixed point $z$ of $f$ (since $v \in K$, $n_1 < n_2$ and $f^{n_1}(v) = z$), and $f^{n_2}(q)$ and $f^{n_2+1}(q)$ lie on the same side of $z$, it follows from Lemma 3 that   
$$
f^n(K) \supset O_f(q) \,\,\, \text{for all} \,\,\, n \ge n_2+2j-2.
$$
\indent Since $f\big(O_f(q)\big) \supset O_f(q)$ and $$\min O_f(q) \le f^{n_3}(q) < \min L < \max L < f^{n_4}(q) \le \max O_f(q),$$we obtain that$$f^n(K) \supset L \,\,\, \text{for all} \,\,\, n \ge n_2+2j-2.$$
\hfill\sq

To characterize continuous transitive maps in terms of dense periodic points, we shall need the following result.  For simplicity, in the following result, we shall assume that $I$ denotes one of the intervals: $(-\infty, \infty)$, $(-\infty, b]$, $[a, \infty)$, $[a, b]$, where $a < b$ are real numbers.  So, we do not need to use relative topology.

\noindent
{\bf Lemma 5.} 
{\it Let $I$ denote one of the intervals: $(-\infty, \infty)$, $(-\infty, b]$, $[a, \infty)$, $[a, b]$, where $a < b$ are real numbers.  Assume that $f$ is a continuous map on $I$ whose periodic points are dense in $I$.  Let $K$ be a proper closed subinterval of $I$ such that $f^s(K) \subset K$ for some positive integer $s$.  Then $f^2(K) = K$ and $f^2(\overline{L}) = \overline{L}$ for each (of the at most two) component $L$ of $I \setminus K$.}  

\noindent
{\it Proof.}  
Let $P$ denote the set of all periodic points of $f$ in $I$.  Let $p$ be a period-$i$ point of $f$ in $K$.  Then, since $f^s(K) \subset K$, $p = (f^i)^s(p) = (f^s)^{i-1}\big(f^s(p)\big) \in f^s(P \cap K)$.  This shows that $P \cap K \subset f^s(P \cap K) \subset P \cap f^s(K) \subset P \cap K$.  Therefore, $$P \cap f^s(K) = f^s(P \cap K) = P \cap K.$$ 
\indent On the other hand, let $q$ be a period-$j$ point of $f$ in $I$.  

If $f^s(q) \in K$, then $$q = (f^j)^s(q) = (f^s)^{j-1}\big(f^s(q)\big) \in (f^s)^{j-1}(K) \subset K \,\,\, \text{since} \,\,\, f^s(K) \subset K.$$  That is, if $q$ is a periodic point of $f$ in $I$ such that $f^s(q) \in K$, then $q$ must be a periodic point of $f$ in $K$.  So, $f^s\big(P \cap (I \setminus K)\big) \subset P \cap (I \setminus K)$.  Consequently, if $q \in P \cap (I \setminus K)$, then 
$$
q = (f^j)^s(q) = (f^s)^{j-1}\big(f^s(q)\big) \in (f^s)^{j-1}\big(P \cap (I \setminus K)\big) \subset f^s\big(P \cap (I \setminus K)\big).
$$ 
This, combined with the above, implies that 
$$
\quad\qquad\qquad\qquad\,\,\, f^s(P \cap K) = P \cap K \,\,\, \text{and} \,\,\, f^s\big(P \cap (I \setminus K)\big) = P \cap (I \setminus K). \quad\qquad\qquad\qquad\qquad (\star)
$$
In particular, the map $f^s(x)$ takes those periodic points of $f$ in $K$ into $K$ and those in $I \setminus K$ into $I \setminus K$.  Since the periodic points of $f$ are dense in $I$, this implies that
$$
\text{if} \,\,\, c \,\,\, \text{is a {\it common} boundary point of} \,\,\, K \,\,\, \text{and} \,\,\, I \setminus K, \, \text{then so is the point} \,\,\, f^s(c).
$$  
\indent Now we have two cases to consider:

Case 1. $I \setminus K$ is an interval.  In this case, $K$ and $I \setminus K$ has a unique common boundary point $d$.  Consequently, $f^s(d) = d$.  In particular, $d$ is also a perodic point of $f$.  Since $I$ denotes one of the intervals: $(-\infty, \infty)$, $(-\infty, b]$, $[a, \infty)$, $[a, b]$, where $a < b$ are real numbers, the closure $\overline{I \setminus K}$ is also a closed interval in $I$.

If there is a point $v$ in $\overline{I \setminus K}$ such that $f^s(v)$ lies in the interior of $K$, then, since the periodic points of $f$ are dense in $I$, there exists a periodic point $q$ of $f$ in $I \setminus K$ such that $f^s(q) \in int(K)$.  This contradicts, by ($\star$), the fact that $f^s\big(P \cap (I \setminus K)\big) = P \cap (I \setminus K)$.  Therefore, we obtain that 
$$
f^s(K) \subset K \,\, \text{(by assumption)} \,\,\, \text{and} \,\,\,  f^s\big(\overline{I \setminus K}\big) \subset \overline{I \setminus K}.
$$  
Since both $K$ and $\overline{I \setminus K}$ are closed {\it subintervals} of $I$ which are invariant under $f^s(x)$, without loss of generality, we may assume that $K$ contains a fixed point of $f$.  By Lemma 3, we obtain that 
$$
f^2(P \cap K) \subset K.
$$
\indent Since $f^2(x)$ is continuous and the periodic points of $f$ are dense in $K$, we obtain that 
$$
f^2(K) \subset K.
$$  
By arguing as those in the beginning of the proof with $s$ replaced by 2, we obtain that 
$$
f^2(K) \subset K \,\,\, \text{and} \,\,\, f^2\big(\overline{I \setminus K}\big) \subset \overline{I \setminus K} \,\,\, \text{and} \,\,\, f^2(d) = f^2\big(K \cap (\overline{I \setminus K})\big) \subset K \cap (\overline{I \setminus K}) = \{ d \}.
$$
This implies that 
$$
f^2(K) = K \,\,\, \text{and} \,\,\, f^2(\overline{I \setminus K}) = \overline{I \setminus K}
$$
as long as $I \setminus K$ is an interval. 

Case 2. $I \setminus K$ consists of exactly two intervals which lie on the opposite sides of $K$.  In this case, $K = [a, b]$ is a compact interval in the {\it interior of} $I$.  Then 
$$
I \setminus K = \big(I \cap (-\infty, a)\big) \cup \big(I \cap (b, \infty)\big)
$$ 
and so, both $a$ and $b$ are common boundary points of $K$ and $I \setminus K$.  Thus, $\big\{ f^s(a), f^s(b) \big\} \subset \{ a, b \}$ and two cases arise:

(A).  Suppose $f^s(a) = a$.  Then by ($\star$), we have
$$
f^s\big(P \cap (I \cap (-\infty, a])\big) = f^s\big(P \cap (I \cap (-\infty, a))\big) \,\, \cup \,\, \{ f^s(a) \} \qquad\qquad\qquad
$$ 
$$
\qquad\qquad\qquad\qquad \subset f^s\big(P \cap (I \setminus K)\big) \cup \big\{ f^s(a) \big\} = \big(P \cap (I \setminus K)\big) \cup \big\{ a \big\} \subset (I \setminus K) \cup \big\{ a \big\},
$$
Since $f^s$ is continuous and the periodic points of $f$ are dense in $I$, we obtain that $f^s\big(I \cap (-\infty, a]\big)$ $\big(\subset (I \setminus K) \cup \big\{ f^s(a) \big\}$\big) is a {\it connected} interval which contains the point $a \, \big(= f^s(a)$\big) and is disjoint from the interior of $K \, \big(= [a, b]\big)$.  Consequently, 
$$
f^s\big(I \cap (-\infty, a]\big) \subset I \cap (-\infty, a].
$$  
Since $I \setminus \big(I \cap (-\infty, a]\big) \, \big(= I \cap (a, \infty)\big)$ is an interval, it follows from Case 1 with $K$ replaced by the closed subinterval $I \cap (-\infty, a]$ of $I$ that 
$$
f^2(a) = a, \, f^2\big(I \cap (-\infty, a]\big) = I \cap (-\infty, a] \,\,\, \text{and} \,\,\, f^2\big(I \cap [a, \infty)\big) = I \cap [a, \infty).
$$  
\indent On the other hand, $\big(I \cap (-\infty, a]\big) \cup K = I \cap (-\infty, b]$ is a closed subinterval of $I$ such that 
\begin{multline*}
$$
f^s\big(I \cap (-\infty, b]\big) = f^s\big(I \cap (-\infty, a] \cup K\big) \subset f^s\big(I \cap (-\infty, a]\big) \cup f^s(K) \subset \big(I \cap (-\infty, a]\big) \cup K = I \cap (-\infty, b] \\ \text{and} \,\,\, I \setminus \big(I \cap (-\infty, b]\big) = I \cap (b, \infty) \,\,\, \text{is an interval}.
$$
\end{multline*}
It follows from Case 1 with $K$ replaced by $I \cap (-\infty, b]$ that 
$$
f^2(b) = b, \, f^2\big(I \cap (-\infty, b]\big) = I \cap (-\infty, b] \,\,\, \text{and} \,\,\, f^2\big(I \cap [b, \infty)\big) = I \cap [b, \infty).
$$  
Consequently, $f^2(K) = f^2\big([a, b]\big) \subset f^2\big(I \cap (-\infty, b]\big) \cap f^2\big(I \cap [a, \infty)\big) = \big(I \cap (-\infty, b]\big) \cap \big(I \cap [a, \infty)\big) = K$.  Since $f^2(a) = a$ and $f^2(b) = b$, we obtain that $f^2(K) = f^2\big([a, b]\big) = [a, b] = K$.  This, together with the above, shows that 
$$
f^2(K) = K, \, f^2\big(I \cap (-\infty, a]\big) = I \cap (-\infty, a], \, \text{and} \,\,\, f^2\big(I \cap [b, \infty)\big) = I \cap [b, \infty).
$$  

(B). Suppose $f^s(a) = b$.  Then, by ($\star$), 
$$
f^s\big(P \cap (I \cap (-\infty, a))\big) \,\, \big(\subset f^s(P \cap (I \setminus K))\big) = P \cap (I \setminus K) \subset I \setminus K.
$$
\indent Similar arguments as those in (A) imply that $f^s\big(I \cap (-\infty, a]\big)$ $\subset \big(\, \overline{I \setminus K} = \big(I \cap (-\infty, a]\big) \cup \big(I \cap [b, \infty)$\big) is a connected interval which contains the point $b$ and is disjoint from the interior of $K \, \big(= [a, b]\big)$.  Therefore, we obtain that 
$$
f^s\big(I \cap (-\infty, a]\big) \subset I \cap [b, \infty).
$$
If $f^s(b) = b$, then similar arguments imply that 
$$
f^s\big(I \cap [b, \infty)\big) \subset I \cap [b, \infty)
$$
which, combined with the above, implies that
$$
f^s(I \setminus K) = f^s\big(P \cap (I \cap (-\infty, a))\big) \cup f^s\big(I \cap (b, \infty)\big) \subset I \cap [b, \infty).
$$
Consequently, it follows from ($\star$) that 
$$
P \cap \big(I \cap (-\infty, a)\big) \subset P \cap (I \setminus K) = f^s\big(P \cap (I \setminus K)\big) \subset f^s(I \setminus K) \subset I \cap [b, \infty)
$$
which is a contradiction.  

So, if $f^s(a) = b$, then $f^s(b) = a$ and similar arguments as those in the beginning of the proof of $(B)$ show that $$
f^s\big(I \cap (-\infty, a]\big) \subset I \cap [b, \infty) \,\,\, \text{and} \,\,\, f^s\big(I \cap [b, \infty)\big) \subset I \cap (-\infty, a].
$$ 
Therefore, we have $f^{2s}(a) = a$ and $f^{2s}(b) = b$ and $$f^{2s}\big(I \cap (-\infty, a]\big) \subset I \cap (-\infty, a], \, f^{2s}\big(I \cap [b, \infty)\big) \subset I \cap [b, \infty) \,\,\, \text{and} \,\,\, f^{2s}(K) \subset K.$$  By similar arguments as those in Part (A) above, we obtain that $f^2(K) = K$, $f^2\big(I \cap (-\infty, a]\big) = I \cap (-\infty, a]$ and $f^2\big(I \cap [b, \infty)\big) = I \cap [b, \infty)$.
\hfill\sq

In the following characterizations of bitransitive maps (see also {\bf\cite{du1}}), statements (2), (3), (8) are taken from {\bf\cite{barge1}}.  Statements (9) and (10) are from {\bf\cite{barge3}}.  Statements (4) and (5) in terms of dense periodic points are derived from {\bf\cite{barge2}}.  Statements (6) and (7) are well-known.  For simplicity, in the following statements (3), (4) $\&$ (5), we assume that $I$ denotes one of the intervals: $[a, b]$, $[a, \infty)$, $(-\infty, b]$, $(-\infty, \infty)$, where $a < b$ are real numbers.  So, we do not need to use relative topology.

\noindent
{\bf Theorem 6.}
{\it Let $I$ denote one of the intervals: $(a, b)$, $(a, b]$, $[a, b)$, $[a, b]$, $(-\infty, b)$, $(-\infty, b]$, $(a, \infty)$, $[a, \infty)$ and $(-\infty, \infty)$, where $a < b$, are real numbers, except in statements (3), (4) $\&$ (5) where $I$ denotes one of $[a, b]$, $[a, \infty)$, $(-\infty, b]$, $(-\infty, \infty)$, and let $f : I \rightarrow I$ be a continuous map.  Then the following statements are equivalent:
\begin{itemize}
\item[{\rm (1)}]
$f$ is transitive on $I$ and there exist a fixed point $z$ of $f$ and a point $c \ne z$ such that either $z < \min \big\{ c, f(c) \big\}$ or $\max \big\{ c, f(c) \big\} < z$, i.e., both points $c$ and $f(c)$ lie on the same side of $z$.  

\item[{\rm (2)}]
For any compact subinterval $K$ of $I$ and any compact subinterval $L$ of int$(I)$, there exists a positive integer $M$ such that $f^n(K) \supset L$ for all $n \ge M$.  Consequently, for any positive integer $m$, any compact intervals $K_i$, $1 \le i \le m$, and any intervals $V_i$, $1 \le i \le m$, in $I$, by taking, for each $1 \le i \le m$, a compact interval $L_i$ in the interior of $V_i$, we obtain that there exists a positive integer $N_m$ such that $f^n(K_i) \supset L_i$ for all $1 \le i \le m$ and all $n \ge N_m$, and so, there exist compact intervals $K_{n,i}^*$ in $K_i$, $1 \le i \le m$, which satisfy that 
$$
f^n(K_{n,i}^*) \subset L_i \subset V_i \,\,\, \text{for all} \,\,\, 1 \le i \le m \,\,\, \text{and all} \,\,\, n \ge N_m.
$$
\big(We phrase it this way for later use in the proof of Theorem 8 below.  Note that since we only require the set $f^n(K_{n,i}^*)$ to be contained in $V_i$, we can choose the compact interval $K_{n,i}^*$ as small as we wish\big).

\item[{\rm (3)}]
For any open interval $U$ in $I$, there exists a positive integer $m(U)$ such that $f$ has a period-$n$ point in $U$ for each integer $n \ge m(U)$ and $f$ has no proper invariant closed intervals in $I$.

\item[{\rm (4)}]
For some infinite set $\mathcal N$ of positive integers which contains an odd integer $\ge 3$, the periodic points of $f$ with periods $n \in \mathcal N$ are dense in $I$ and $f$ has no proper invariant closed intervals in $I$. 

\item[{\rm (5)}]
The set of periodic points of $f$ is dense and $f^2$ has no proper invariant closed interval in $I$. 

\item[{\rm (6)}]
$f$ is (topologically) mixing, i.e., for any two non-empty open subsets $U$ and $V$ of $I$, there exists a positive integer $m(U, V)$ such that $f^n(U) \cap V \ne \emptyset$ for all $n \ge m(U, V)$.  
  
\item[{\rm (7)}]
$f$ is (topologically) weakly mixing, i.e., the map $f \times f$ defined by $\big(f \times f\big)\big((x, y)\big) = \big(f(x), f(y)\big)$ is transitive on $I \times I$.   

\item[{\rm (8)}]
For any infinite subset $\mathcal M$ of the positive integers, there is a residual set $S$ \big(i.e., a set containing the intersection of countably infinitely many open dense sets\big) in $I$ such that $S$ is dense in $I$ and if $x \in S$ then the set $\big\{ f^n(x) : n \in \mathcal M \big\}$ is dense in $I$.  

\item[{\rm (9)}]   
$f^{2k}$ is transitive on $I$ for some positive integer $k$.   

\item[{\rm (10)}]
$f$ is totally transitive, i.e., $f^k$ is transitive on $I$ for each integer $k \ge 1$. 
\end{itemize}}

\noindent
{\it Proof.}
The implication ${(1)} \Rightarrow {(2)}$ follows from Theorem 4(1) and the implication ${(3)} \Rightarrow {(4)}$ is clear.  

${(2)} \Rightarrow {(3)}$.  Suppose $U$ is an open interval in $I$ whose closure is contained in the interior of $I$.  It is clear that there exist two disjoint compact intervals $K$ and $L$ in $U$ such that 
$$
f(K) \cap (K \cup L) = \emptyset.
$$
By {(2)}, there is a positive integer $N_0$ such that $f^n(K) \cap f^n(L) \supset K \cup L$ for all $n \ge N_0$.  Let 
$$
N = m_1^{n_1}m_2^{n_2} \cdots m_k^{n_k} \, (> N_0),
$$
where $\{ m_1, m_2, \cdots, m_k \}$ is the set of all prime numbers $\le N_0$ and each $n_i, 1 \le i \le k$, is a fixed positive integer such that $m_i^{n_i} > N_0$. Let $m$ be any integer $\ge N > N_0$.  

If $m$ has a prime factor $s > N_0$, then since $f^s(K) \cap f^s(L) \supset K \cup L$, by considering, with respect to the map $f^s$, the cycle $K^{(m/s)-1}LK$ of length $m/s$, we obtain a period-$(m/s)$ point $w$ of $f^s$ in $K$ whose {\it whole orbit} $O_{f^s}(w)$ is contained in $K \cup L$.  Since $f(K) \cap (K \cup L) = \emptyset$, the orbit $O_f(w) \,\, \big(\supset \big\{ f(w) \big\} \cup O_{f^s}(w)\big)$ contains at least $1 + (m/s)$ points.  So, $w$ is a period-$m$ point of $f$ in $K \subset U$.  

If all prime factors of $m$ are $\le N_0$, write $m = m_1^{n_1'}m_2^{n_2'} \cdots m_k^{n_k'}$, with each $n_i' \ge 0$ and $n_j' > n_j$ for some $1 \le j \le k$.  Let $\ell = m_j^{n_j}$.  Then $\ell > N_0$ and so $f^\ell(K) \cap f^\ell(L) \supset K \cup L$.  Consequently, $f^\ell$ has a periodic point of least period $m/\ell$ in $K \subset U$ which turns out to be a period-$m$ point of $f$.  On the other hand, (2) implies trivially that $f$ has no proper invariant closed interval in $I$.  

${(4)} \Rightarrow {(5)}$.  Suppose $f^2$ had a proper invariant closed subinterval $L$ of $I$.  By Lemma 5, we may assume that $L$ contains a fixed point $z$ of $f$.  Thus, $L \cap f(L)$ is an interval or consists of one point.  Since $f\big(L \cap f(L)\big) \subset L \cap f(L)$, we have $f\big(\overline{L \cap f(L)}\big) \subset \overline{L \cap f(L)}$.  It follows from assumption that $L \cap f(L) = \{ z \}$.  Since we also have $f\big(f(L)\big) \subset L$, $f$ can only have periodic points of even periods (except fixed points) in $L$.  This contradicts the assumption.  Therefore, $f^2$ has no proper invariant closed intervals in $I$.  

${(5)} \Rightarrow {(2)}, {(6)}$.  Let $U$ and $V$ be nonempty open intervals in $I$ and let $K$ and $L$ be compact intervals in $U$ and $V$ respectively.  Let $p$ be a periodic point of $f$ in $K$.  Let $q$ be a periodic point of $f$ in $I$ with $q_1 = \min O_f(q), q_2 = \max O_f(q)$ and $[q_1, q_2] \subset int(I)$.  Let $m$ be a fixed positive integer such that $f^m(p) = p$ and $f^m(q) = q$.  Let 
$$
W = \bigcup_{k=0}^\infty \, (f^m)^k(K).
$$
Then since $f^m(p) = p \in K$, $W$ is an {\it interval} and $f^m(W) \subset W$.  So, $f^m(\ov W) \subset \ov W$.  By Lemma 5, $f^2(\ov W) = \ov W$.  By assumption, $\ov W = I$.  Since $W$ is an interval, $W \supset int(W) = int(I) \supset [q_1, q_2]$.  Therefore, there exist two positive integers $i$ and $j$ such that $q_1 \in (f^m)^i(K)$ and $q_2 \in (f^m)^j(K)$ and, since $f^m(q_1) = q_1$ and $f^m(q_2) = q_2$, we have $q_1 \in f^{ijm}(K)$ and $q_2 \in f^{ijm}(K)$ and so $O_f(q) \subset f^{ijm}(K)$.  Since $f\big(O_f(q)\big) = O_f(q)$, we have
$$
O_f(q) \subset f^n(K), \, \text{and so}, \, [q_1, q_2] \subset f^n(K) \,\,\, \text{for all} \,\,\, n \ge ijm.
$$
Consequently,
$$
\text{if we choose} \,\,\, q \,\,\, \text{to be a periodic point of} \,\,\, f \,\,\, \text{such that} \,\,\, O_f(q) \subset int(L),
$$
then we easily obtain that (6) holds, i.e., $f$ is topologically mixing on $I$.  On the other hand,\begin{multline*}$$\text{if we choose} \,\,\, q \,\,\, \text{and} \,\,\, \hat q \,\,\, \text{to be periodic points of} \,\,\,f \,\,\, \text{with} \\ q < \min L < \max L < \hat q \,\,\, \text{and} \,\,\, O_f(q) \cup O_f(\hat q) \subset f^k(K) \,\,\, \text{for some} \,\,\, k \ge ijm.$$\end{multline*}then we obtain that $f^n(K) \supset L$ for all sufficiently large integers $n$.  Therefore, (2) holds. 

${(6)} \Rightarrow {(7)}$.  Let $U_1, U_2, V_1, V_2$ be nonempty open sets in $I$.  By {(6)}, there exist two positive integers $n_1, n_2$ such that $f^n(U_i) \cap V_i \ne \emptyset$, for all $n \ge n_i$ and all $i = 1, 2$.  In particular, $f^{n_1+n_2}(U_1 \times U_2) \cap (V_1 \times V_2) \ne \emptyset$.  So, $f \times f$ is transitive on $I \times I$.  That is, $f$ is weakly topologically mixing on $I$.  

${(7)} \Rightarrow {(2), (6), (10)}$.  Firstly we show that\begin{multline*}$$\text{if} \,\,\, a < b \,\,\, \text{are any two points in} \,\,\,I \,\,\, \text{and} \,\,\, c < d \,\,\, \text{are any two points in} \,\,\, int(I), \, \text{then} \\ f^i\big([a, b]\big) \supset [c, d] \,\,\, \text{for some positive integer} \,\,\, i.$$\end{multline*}Indeed, let $U$ and $V$ be two nonempty open intervals in $I$ such that $U \subset (\inf I, c)$ and $V \subset (d, \sup I)$.  Since $f$ is weakly topologically mixing, there is a positive integer $i$ such that $\big(f \times f\big)^i\big((a, b) \times (a, b)\big) \cap (U \times V) \ne \emptyset$.  Therefore, $f^i\big((a, b)\big) \cap U \ne \emptyset$ and $f^i\big((a, b)\big) \cap V \ne \emptyset$.  Consequently,
$$
f^i\big([a, b]\big) \supset f^i\big((a, b)\big) \supset [c, d] \,\,\, \text{and so}, \, f \,\,\, \text{is transitive on} \,\,\, I.
$$
If we choose $c = a$ and $d = b$, then this shows that $f$ has a periodic point in $[a, b]$.  Since $a < b$ are arbitrary in $I$, this implies that $f$ has dense periodic points in $I$.  Furthermore, since $f$ is transitive on $I$, there is a positive integer $j$ such that $f^j(U) \cap V \ne \emptyset$.  So, there is a nonempty open interval $W$ in $U$ such that $f^j(W) \subset V$.  Since the periodic points of $f$ are dense in $I$, there is a periodic point $p$ in $W \subset U$ such that $f^j(p) \in f^j(W) \subset V$ and the whole orbit $O_f(p)$ is contained in the interior of $I$.  Let $p_1 = \min O_f(p)$ and $p_2 = \max O_f(p)$.  Then it follows from what we just proved that $f^k\big([a, b]\big) \supset [p_1, p_2]$ for some $k \ge 1$.  Since
$$
p_1 \le p \in W \subset U \subset (\inf I, c] \,\,\, \text{and} \,\,\, f^j(p) \in V \subset [d, \sup I] \,\,\, \text{and} \,\,\, f^j(p) \le p_2,
$$
we obtain that $f^n\big([a, b]\big) \supset [p_1, p_2] \supset [c, d]$ for all $n \ge k$.  This establishes trivially {(2), (6)} and {(10)}.  

${(2)} \Rightarrow (8)$.  Let $\big\{ V_k : k \ge 1 \big\}$ be an enumeration of all open intervals in the interior of $I$ with rational endpoints.  For each $V_k$ and any open interval $U$ in $I$, (2) implies that $f^n(U) \cap V_k \ne \emptyset$ for all sufficiently large $n$ (in $\mathcal M$, in particular).  So, the set $\bigcup_{n \in \mathcal M} \, \big(f^{-n}(V_k)\big)$ is open and dense in $I$.  

To show that the set
$$
S = \bigcap_{k \ge 1} \, \big(\bigcup_{n \in \mathcal M} f^{-n}(V_k)\big)
$$
is dense in $I$ {\it without} resorting to Baire Category Theorem, we argue as follows using an idea of Babilonov\'a {\bf\cite{bab}} (cf. the proof of Theorem 8 below): 
$$
\text{Let} \,\,\, K \,\,\, \text{be any compact interval in} \,\,\, I \,\,\, \text{and, for each} \,\,\, k \ge 1, \, \text{let} \,\,\, L_k \,\,\, \text{be a compact interval in} \,\,\, V_k.
$$
By (2), there is an integer $n_1$ in $\mathcal M$ such that $f^{n_1}(K) \supset L_1$.  Let $K_1$ be a compact interval in $K$ with length $< 1/2$ such that $f^{n_1}(K_1) \subset L_1$.  Again, by (2), there is an integer $n_2 > n_1$ in $\mathcal M$ such that $f^{n_2}(K_1) \supset L_2$.  Let $K_2$ be a compact interval in $K_1$ with length $< 1/4$ such that $f^{n_2}(K_2) \subset L_2$.  Proceeding in this manner indefinitely, we obtain a sequence $n_1 < n_2 < \cdots$ of positive integers in $\mathcal M$ and a sequence $< K_i >$ of compact intervals such that, for each $i \ge 1$, $K_i \subset K_{i-1}$ the length of $K_i$ is $< 1/2^i$ and $f^{n_i}(K_i) \subset L_i \subset V_i$, where $K_0 = K$.  It is clear that the set $\bigcap_{i \ge 0} K_i$ consists of exactly one point which is contained in $\bigcap_{k \ge 1} K_k \subset \bigcap_{k \ge 1} \big(\bigcup_{n \in \mathcal M} f^{-n}(V_k)\big) = S$.  Since $K$ is arbitrary in $I$, $S$ is dense in $I$.  This confirms (8).  

The implication ${(8)} \Rightarrow {(9)}$ is clear and, by Theorem 4, we have ${(9)} \Rightarrow {(10)} \Rightarrow {(1)}$. 

To prove the implication $(9) \Rightarrow (2)$, we may proceed as follows: Let $z$ be a fixed point of $f$ in the interior of $I$ and let $K$ be any compact interval in the interior of $I$.  Without loss of generality, we may assume that $z < x$ for all $x \in K$.  By Lemma 2, there exist a point $v$ in $int(K)$ and a positive integer $s$ such that $f^s(v) = z$.  Since, for some positive integer $k$, $f^{2k}$ is transitive on $I$, there exist a point $u$ in the {\it interior} of $K$ whose orbit $O_{f^{2k}}(u)$ with respect to $f^{2k}$ is dense in $I$ and a positive integer $t > s$ such that $\max K < f^t(u)$.  Thus, $f^t(K) \supset K$.  Therefore, $K$ contains a periodic point of $f$.  Since $K$ is arbitrary in the interior of $I$, this shows that $f$ has dense periodic points in $I$.  

Let $L$ be any compact interval in $int(I)$.  Since the orbit $O_{f^{2k}}(u)$ of $u (\in int(K))$ with respect to $f^{2k}$ is dense in $I$, for some positive integers $s < n_1 < n_2$, we have$$f^{2kn_1}(u) < \min \big\{ z, \min L \big\} < \max \big\{ z, \max L \big\} < f^{2kn_2}(u).$$Since $u \in int(K)$ and since we have just shown that the periodic points of $f$ are dense in $I$,
$$
\text{there is a period-}m \,\,\, \text{point} \,\,\, p \,\,\, \text{of} \,\,\, f \,\,\, \text{in} \,\,\, K \,\,\, \text{which is so close to} \,\,\, u
$$
that the orbit $O_f(p)$ is contained in $int(I)$ and
$$
f^{2kn_1}(u) \approx f^{2kn_1}(p) < \min \big\{ z, \min L \big\} < \max \big\{ z, \max L \big\} < f^{2kn_2}(p) \approx f^{2kn_2}(u).
$$
Consequently, the interval $f^s(K)$ contains the fixed point $z$ $(= f^s(v))$ and the period-$m$ point $f^s(p)$ of $f$ whose even iterates are distributed on both sides of $z$ and so, the orbit $O_f(p)$ is not {\it separated} by the fixed point $z$. By Lemma 3, we must have $f^i\big(f^s(K)\big) \supset O_f(p)$ for all $i \ge 2m-2$.  Therefore, $f^n(K) \supset \big[\min O_f(p), \max O_f(p)\big] \supset L$ for all $n \ge s + 2m - 2$.  This confirms $(9) \Rightarrow (2)$.

$(10) \Rightarrow (2)$.  Suppose $f$ is totally transitive on $I$.  By following arguments as those in the above proof of the implication $(9) \Rightarrow (2)$, we obtain that the periodic points of $f$ are dense in $I$.  

Let $K$ and $L$ be any compact intervals in $I$ and $int(I)$ respectively.  Since $f$ is (totally) transitive on $I$, there is a point $w$ in the {\it interior} of $K$ such that $O_f(w)$ is dense in $I$.  Let $0 < m_1 < m_2 < m_3$ be integers such that $$f^{m_1}(w) < \min L < \max L < f^{m_2}(w) \,\,\, \text{and} \,\,\, f^{m_3}(w) \in int(K).$$Since the periodic points of $f$ are dense in $I$, there is a period-$k$ point $q$ of $f$ in the {\it interior} of $K$ which is so close to $w$ that $$f^{m_1}(w) \approx f^{m_1}(q) < \min L < \max L < f^{m_2}(q) \approx f^{m_2}(w) \,\,\, \text{and} \,\,\, f^{m_3}(w) \approx f^{m_3}(q) \in K$$and the orbit $O_f(q)$ is contained in the {\it interior} of $I$.    

Let $\hat q = f^{m_3}(q)$.  Then $f^k(q) = q$, $f^k(\hat q) = \hat q$, and $f^k\big([q, \hat q]\big) \supset [q, \hat q]$.  Consequently, the sequence$$[q, \hat q] \subset f^k\big([q, \hat q]\big) \subset (f^k)^2\big([q, \hat q]\big) \subset (f^k)^3\big([q, \hat q]\big) \subset \cdots$$forms an increasing sequence of compact intervals in $I$.  Since $f$ is totally transitve and so, $f^k$ is transitive on $I$, the union $\cup_{i \ge 0} (f^k)^i\big((q, \hat q)\big) \, \big(\subset \cup_{i \ge 0} (f^k)^i([q, \hat q])\big)$ is dense in $I$.  So, $\cup_{i \ge 0} (f^k)^i\big([q, \hat q]\big)$ is dense in $I$.  Therefore, for some positive integer $j$,$$L \subset \big[\min O_f(q), \max O_f(q)\big] \subset (f^k)^j\big([q, \hat q]\big) \subset (f^k)^j(K).$$Since $f\big(O_f(q)\big) = O_f(q)$, we see that $L \subset \big[\min O_f(q), \max O_f(q)\big] \subset f^n(K)$ for all $n \ge jk$.  This shows that $(10) \Rightarrow (2)$.  
\hfill\sq

\vspace{.2in}
In Theorem 6 (9) $\&$ (10), we see that $f^2$ is transitive on $I$ if and only if $f$ is totally transitive on $I$.  In fact, we can say more. 

\noindent
{\bf Proposition 7.}
{\it Let $I$ denote one of the following intervals: $(a, b)$, $(a, b]$, $[a, b)$, $[a, b]$, $(-\infty, b)$, $(-\infty, b]$, $(a, \infty)$, $[a, \infty)$ and $(-\infty, \infty)$, where $a < b$ are real numbers.  Suppose $u$ is a point in $I$.  Then the following hold:
\begin{itemize}
\item[{\rm (1)}]
If the orbit $O_{f^2}(u)$ is dense in $I$, then, for each $k \ge 1$, the orbit $O_{f^k}(u) = \big\{ f^{ik}(u) : i = 0, 1, 2, \cdots \big\}$ is dense in $I$.

\item[{\rm (2)}]
If the orbit $O_f(u)$ is dense in $I$, then, for each odd $\ell \ge 1$, the orbit $O_{f^\ell}(u) = \big\{ f^{i\ell}(u) : i = 0, 1, 2, \cdots \big\}$ is dense in $I$.
\end{itemize}} 

\noindent
{\it Ptoof.}
Suppose $O_{f^2}(u)$ is dense in $I$.  Let $K$ be any compact interval in the interior of $I$ which contains no fixed points of $f$ and let $L = K$.  It follows from Theorem 4(1) that $f^n(K) \supset L=K$ for all sufficiently large integer $n$.  Consequently, $K$ contains periodic points of $f$ whose periods are sufficiently large {\it prime} numbers.  Therefore, since $K$ is arbitrary, this shows that the set of all periodic points of $f$ whose periods are large odd {\it prime} numbers is dense in $I$.  Let $k \ge 2$ be a fixed integer.  Let $V$ be any open interval in $I$ and let $q \in V$ be a period-$\ell$ point of $f$ such that $\ell$ is a prime number and $\ell > k$.  Write $\ell = ik + r$, where $i$ and $r$ are integers with $i \ge 1$ and $0 < r < k$.  Then $k$ and $r$ are coprime.  Let $t_1$ and $t_2$ be positive integers such that $rt_1 - kt_2 = 1$.  Since $O_f(u)$ is dense in $I$, let $f^n(u)$ be a point so close to the period-$\ell$ point $q$ that 
$$
\{ f^n(u), f^{n+\ell}(u), f^{n+2\ell}(u), \cdots, f^{n+(kt_1-1)\ell}(u), f^{n+(kt_1)\ell}(u) \} \subset V.
$$ 
Write $n = jk + s$, where $j$ and $s$ are integers such that $j \ge 0$ and $0 \le s < n$.  Then $n + \big((k-s)t_1\big)\ell = jk + s + (k-s)t_1(ik+r) = jk + s + (k-s)t_1ik + (k-s)rt_1 = jk + s + (k-s)t_1ik + (k-s)kt_2 + (k-s) \equiv 0 \mod k$.  So, $f^{n + ((k-s)t_1)\ell}(u) \in V$.  This implies that the orbit $O_{f^k}(u)$ is dense in $I$.  So, (1) holds.  

As for the proof of (2), without loss of generality, we may assume that $O_f(u)$ is dense in $I$ while $O_{f^2}(u)$ is not.  It follows from Theorem 4 that there exists a fixed point $z$ of $f$ such that 
$$
f\big(I \cap (-\infty, z]\big) = I \cap [z, \infty) \,\,\, \text{and} \,\,\, f\big(I \cap [z, \infty)\big) = I \cap (-\infty, z].
$$
Without loss of generality, we may assume that $u \in I \cap (-\infty, z]$.  Let $g = f^2$.  Then $O_g(u)$ is dense in $I \cap (-\infty, z]$.  Suppose $O_{g^2}(u)$ were not dense in $I \cap (-\infty, z]$.  Then it follows from Theorem 4 that there exists a fixed point $\hat z$ of $g$ in the interior of $I \cap (-\infty, z]$ such that 
$$
g\big(I \cap (-\infty, \hat z]\big) = I \cap [\hat z, z] \,\,\, \text{and} \,\,\, g\big(I \cap [\hat z, z]\big) = I \cap (-\infty, \hat z]
$$
which is impossible because $z$ is a fixed point of $f$ and of $g$ $(= f^2)$.  Therefore, $O_{g^2}(u)$ is dense in $I \cap (-\infty, z]$.  It follows from (1) above that, 
$$
\text{for each} \,\,\, k \ge 1, \, \text{the orbit} \,\,\, O_{g^k}(u) = \{ f^{2ik}(u): i = 0, 1, 2, \cdots \} \,\,\, \text{is dense in} \,\,\, I \cap (-\infty, z].
$$
Since, for each odd $\ell \ge 1$, $f^\ell\big(I \cap (-\infty, z]\big) = I \cap [z, \infty)$ and the orbit $O_{f^{2\ell}}(u)$ is dense in $I \cap (-\infty, z]$, we obtain that the orbit $O_{f^\ell}(u) = \{ f^{i\ell}(u): i = 0, 1, 2, \cdots \}$ is dense in $I$.  This confirms (2).
\hfill{\sq}

\noindent
{\bf Remark 3.}
Let $I$ denote one of the intervals: $(-\infty, \infty)$, $(-\infty, b]$, $[a, \infty)$, $[a, b]$, where $a < b$ are real numbers and let $f : I \rightarrow I$ be a continuous map.  Assume that $f$ has dense periodic points in $K$ and, for some point $w$ in $I$, $f^2(w) \ne w$.  Let 
$$
\mathcal C = \big\{ L : L \,\,\, \text{is a closed interval in} \,\,\, I \,\,\, \text{such that} \,\,\, w \in L \,\,\, \text{and} \,\,\, f^2(L) \subset L \big\}.
$$ 
Then $I \in \mathcal C$.  Let $J = \bigcap_{L \in \mathcal C} L$ and let $g = f^2$.  Then $w \in J$, $\bar J = \overline {\bigcap_{L \in \mathcal C} L} \subset \overline L = L$ and $f^2\bigl(\overline {\bigcap_{L \in \mathcal C} L}\bigr) \subset f^2(\overline L) = f^2(L) \subset L$ for each $L \in \mathcal C$.  Therefore, 
$$
\bar J = J, \, w \in J \,\,\, \text{and} \,\,\, f^2(\bar J) \subset \bigcap_{L \in \mathcal C} L = J = \bar J.
$$
Suppose $J$ were not a {\it minimal} invariant closed subinterval of $f^2$ in $I$ and let $J'$ be a {\it proper} invariant closed interval of $f^2$ in $J$.  Then $J'$ cannot contain the point $w$.  Otherwise, $J' \in \mathcal C$ and so, $J \subset J'$ which is a contradiction.  Therefore, $w \notin J'$.  By Lemma 5, let $J^*$ be the component of the at most two components of $J \setminus J'$ which contains the point $w$ and is also invariant with respect to $f^2$.  This is again a contradiction.  So, $J$ is a {\it minimal} invariant closed interval of $f^2 = g$ in $I$.  

Since $g$ has dense periodic points in $J$, again by Lemma 5, $J$ is also a {\it minimal} invariant closed interval of $g^2$ in $I$.  By Theorem 6(5)(9), $g^2 \, (= f^4)$ is transitive on $J$.  It is clear that, by Lemma 5, $f^2$ has at most countably infinitely many {\it minimal} invariant closed intervals $\hat J_i$ in $I$ such that $f^2(\hat w) \ne \hat w$ for some point $\hat w \in \hat J$.  Let $J_1, J_2, \cdots, J_n, \cdots$ be such an (maybe finite) enumeration.  Then $f^2(J_i) = J_i$ and $f^4$ is transitive on $J_i$ for each $i \ge 1$ and $f^2(x) = x$ for all $x \notin \bigcup_{i \ge 1} J_i$.  In particular, this proves the main result of {\bf{\cite{barge1, bl}}} for compact $I$ (see also {\bf\cite{ru}}) without using the notion of snakelike continua introduced in {\bf\cite{bing}}.  

$\qquad\qquad\qquad\qquad \aleph \qquad\qquad\qquad\qquad\qquad \aleph \qquad\qquad\qquad\qquad\qquad \aleph$

\indent Recall that, when we say an interval, we mean a non-degenerate interval in the real line.  Let $I$ denote any interval, i.e., let $I$ denote one of the intervals: $(a, b)$, $(a, b]$, $[a, b)$, $[a, b]$, $(-\infty, a)$, $(-\infty, a]$, $(b, \infty)$ $[b, \infty)$, $(-\infty, \infty)$, where $a < b$ are real numbers.  When we say an open (compact respectively) interval, we mean an open (compact respectively) interval with respect to the topology on the real line.  Same when we say the interior (the closure respectively) of an interval.  When we say a relative open interval, we mean an interval which is relatively open in the relative topology on $I$.  Let $c$ be a point in $I$.  If there exist a (finite) point $y$ in the {\it closure} of $I$ and a strictly increasing sequence $< n_i >$ of positive integers such that $\lim_{i \to \infty} f^{n_i}(c) = y$, then we call $y$ an $\omega$-limit point of $c$ (with respect to $f$).  Note that here we do not require $y$ to be a point of $I$.  If $y \in \omega_f(c)$ and $y \notin I$, then $y$ is an endpoint of $I$.  

On the other hand, let $S$ be a subset of $I$ with at least two distinct points and let $\beta$ be a positive number or $\beta = \infty$.  When $\beta$ is finite, we say that $S$ is a $\beta$-scrambled set of $f$ if, for any two distinct points $x$ and $y$ in $S$, we have 
$$\liminf_{n \to \infty} \big|f^n(x) - f^n(y)\big| = 0 \quad \text{and} \quad \limsup_{n \to \infty} \big|f^n(x) - f^n(y)\big| \ge \beta.
$$  
We say that $S$ is an $\infty$-scrambled set of $f$ if, for any two distinct points $x$ and $y$ in $S$, we have 
$$\liminf_{n \to \infty} \big|f^n(x) - f^n(y)\big| = 0 \quad \text{and} \quad \limsup_{n \to \infty} \big|f^n(x) - f^n(y)\big| = \infty.
$$
\indent Following Babilonov\'a {\bf\cite{bab}}, we call such a $\beta$-scrambled set $S$ an {\it extremely scrambled set} of $f$ if, for any two distinct points $c$ and $d$ in $S$ and any compact interval $K$ in the interior of $I$, there exists a positive integer $n$, depending on $x$, $y$ and $K$, such that the interval $[f^{n}(c) : f^{n}(d)]$ with $f^n(c)$ and $f^n(d)$ as endpoints contains $K$.  

In {\bf\cite{bab}}, Babilonov\'a uses Theorem 6(2) to show the existence of {\it extremely scrambled} sets for any continuous {\it bitransitive} map on a compact interval while in {\bf\cite{xiong}}, Xiong and Yang show a general result which, in particular, implies that if $f : I \longrightarrow I$ is a continuous {\it bitransitive} map, then for any given strictly increasing sequence $< p_n >$ of positive integers, there exists an extremely scrambled set $B$ such that, for any non-empty subset $A$ of $B$ and any continuous map $g : A \longrightarrow I$, there exists a subsequence $< q_n >$ of $< p_n >$ such that 
$$
\lim_{n \to \infty} f^{q_n}(x) = g(x) \quad \text{for every} \,\,\, x \in A.
$$
\indent Let $\mathcal M$ be an infinite set of positive integers and let $\beta$ be a positive number or $\beta = \infty$.  We need a terminology which is inspired by Xiong and Yang {\bf\cite{xiong}}.  We call a subset $S$ of $I$ a $\beta$-scrambled set (with respect to $\mathcal M$) if, for any two distinct points $x$ and $y$ in $S$, we have
$$
\liminf_{\substack{{n \to \infty} \\ n \in \mathcal M}} \big|f^n(x) - f^n(y)\big| = 0 \quad \text{and} \quad \limsup_{\substack{{n \to \infty} \\ n \in \mathcal M}} \big|f^n(x) - f^n(y)\big| \ge \beta.
$$
\indent Following Mai {\bf\cite{mai}}, we say that the set $C$ is synchronously proximal ($\mathcal M$-synchronously proximal respectively) to the point $v$ if there exists a strictly increasing sequence $< m_i >_{i \ge 1}$ of positive integers (in $\mathcal M$ respectively) such that the diameter of the set $f^{m_i}(C) \cup \{ v \}$ tends to zero as $i$ tends to $\infty$.  On the other hand, we say that $C$ is {\it dynamically} synchronously proximal to the point $v$ under $f$ if there exists a strictly increasing sequence $< n_i >_{i \ge 1}$ of positive integers such that the diameter of the set $f^{n_i}\big(C \cup \{ v \}\big)$ tends to zero as $i$ tends to $\infty$.  It is possible that $C$ is ({\it dynamically} respectively) synchronously proximal to many different points.  We say that $C$ is a Cantor set if it is a compact, perfect and totally disconnected set in the real line.

In the following, we extend ideas of Babilonov\'a {\bf\cite{bab}} and Xiong and Yang {\bf\cite{xiong}} to show that every continuous {\it bitransitive} map on $I$ has extremely scrambled sets with a special chaotic property which is different from that of Xiong and Yang {\bf\cite{xiong}}.  The construction of the scrambled sets is similar to that of the classical construction of the Cantor ternary set except that now, at each atage, we take ever more steps than just one step in the classical construction.

\noindent
{\bf Theorem 8.}
{\it Let $I$ denote one of the following intervals: $(a, b)$, $(a, b]$, $[a, b)$, $[a, b]$, $(-\infty, b)$, $(-\infty, b]$, $(a, \infty)$, $[a, \infty)$ and $(-\infty, \infty)$, where $a < b$ are real numbers, and let $f : I \rightarrow I$ be a continuous bitransitive map.  Let $\mathcal M$ be an infinite set of positive integers, let
$$
\beta =  \begin{cases}
               \infty, & \text{if $I$ is unbounded}, \\
               |I| = \text{the length of $I$}, & \text{if $I$ is bounded}, \\
       \end{cases}
$$
and let 
$$\delta = \inf_{n \ge 1} \big\{ \sup_{x \in I} \big|f^n(x) - x\big|  \big\}.$$ 
Let 
\begin{multline*}
$$
\{ v_1, \, v_2, \, v_3, \, \cdots \} \,\,\, \text{be a countably infinite subset of the closure of} \,\,\, I \\ \text{(later we shall choose these points} \,\,\, v_i\text{'s for some special purposes)}
$$
\end{multline*} 
and let
$$
\{ x_1, \, x_2, \, x_3, \, \cdots \} \,\,\, \text{be a countably infinite set of points in} \,\,\, I \,\,\, \text{with some finite} \,\,\, \omega\text{-limit points}
.$$  
Then $\delta > 0$ and there exists countably infinitely many pairwise disjoint Cantor sets $\mathcal S^{(1)}, \, \mathcal S^{(2)}, \, \mathcal S^{(3)}, \, \cdots$ of totally transitive points of $f$ in $I$ such that 
\begin{itemize}
\item[{\rm (1)}]
each nonempty open set in $I$ contains countably infinitely many these Cantor sets $\mathcal S^{(j)}$'s;

\item[{\rm (2)}]
for each integer $\ell \ge 1$, the set $\bigcup_{j=1}^\ell \mathcal S^{(j)}$ is $\mathcal M$-synchronously proximal to the point $v_m$ for each integer $m \ge 1$ and dynamically synchronously proximal to the point $f^i(x_m)$ for each integer $m \ge 1$ and each integer $i \ge 0$;

\item[{\rm (3)}]
the set $\mathbb S = \bigcup_{j=1}^\infty \, \mathcal S^{(j)}$ is a dense $\beta$-scrambled (i.e., dense extremely scrambled) set (with respect to $\mathcal M$) of $f$ in $I$;

\item[{\rm (4)}]
for any point $x$ in $\bigcup_{i=0}^\infty f^i\big(\{x_1, x_2, x_3, \cdots\}\big)$ and any point $c$ in $\widehat {\mathbb S} = \bigcup_{i=0}^\infty \, f^i(\mathbb S)$, 
we have $$\liminf_{n \to \infty} \big|f^n(x) - f^n(c)\big| = 0 \,\,\, \text{and} \,\,\, \limsup_{n \to \infty} \big|f^n(x) - f^n(c)\big| \ge \frac \beta2,$$i.e., the set $\{ x, c \}$ is a ($\frac {\beta}2$)-scrambled set (not necessarily with respect to $\mathcal M$) of $f$ in $I$.

\item[{\rm (5)}]
for any $c \ne d$ in $\mathbb S$ and any integer $i \ge 0$, we have 
$$
\limsup_{\substack{n \to \infty \\ n \in \mathcal M}} \big|f^n\big(f^i(c)\big) - f^n\big(f^i(d)\big)\big| = \beta.
$$

\item[{\rm (6)}]
the set $$\widehat {\mathbb S} = \bigcup_{i=0}^{\infty} \, f^i(\mathbb S)$$ is a dense {\rm invariant} $\delta$-scrambled set (with respect to $\mathcal M$) of $f$ in $I$;    
\end{itemize}
Consequently, by taking 
$$
\mathcal M = \{ n!: n = 1, 2, 3, \cdots \} \,\,\, \text{and appropriate} \,\,\, v_1, v_2, v_3, \cdots,
$$ 
the above Cantor sets $\mathcal S^{(1)}$, $\mathcal S^{(2)}$, $\mathcal S^{(3)}$, $\cdots$, of totally transitive points of $f$ can be chosen so that Parts (1) $\&$ (2) hold and the sets $\mathbb S = \bigcup_{j=1}^\infty \, \mathcal S^{(j)}$ and $\widehat {\mathbb S} = \bigcup_{i=0}^\infty \, f^i(\mathbb S)$ satisfy Parts (4) $\&$ (5) and, for each integer $n \ge 1$, the set $\mathbb S = \bigcup_{j=1}^\infty \, \mathcal S^{(j)}$ is a dense $\beta$-scrambled set of $f^n$ and the set $\widehat {\mathbb S} = \bigcup_{i=0}^\infty \, f^i(\mathbb S)$ is a dense {\rm invariant} $\delta$-scrambled set of $f^n$ in $I$.}

\noindent
{\it Proof.}
Since, for each integer $m \ge 1$, all points in the orbit $O_f(x_m)$ have the same $\omega$-limit set.  We can consider the invariant set $\{ f^i(x_m): m \ge 1, i \ge 0 \}$ and rename it $\{ x_1', x_2', x_3', \cdots \}$.  However, here we still consider the set $\{ x_1, x_2, x_3, \cdots \}$ instead of the invariant set $\{ f^i(x_m): m \ge 1, i \ge 0 \}$.  Furthermore, for any infinite set $\mathcal M$ of positive integers, by choosing the pertaining positive integers from $\mathcal M$ at odd stages and from the set $\{ n!: n = 1, 2, 3, \cdots \}$ at even stages, we can replace Part (3) with Part ($3'$): the set $\mathbb S = \bigcup_{j=1}^\infty \, \mathcal S^{(j)}$ is a dense $\beta$-scrambled set (with respect to $\mathcal M$) of $f$ in $I$ and, for each $n \ge 1$, the set $\mathbb S$ is a dense $\beta$-scrambled set of $f^n$.  However, here we only choose the pertaining positive integers from the set $\mathcal M$ for simplicity.
 
Firstly, we introduce some notations.  

Let $U_1$, $U_2$, $U_3$, $\cdots$ be an enumeration of all open intervals in $I$ with rational endpoints.  

Let $a_1$, $a_2$, $a_3$, $\cdots$ and $b_1$, $b_2$, $b_3$, $\cdots$ be two sequences of points in $I$ such that
$$
\lim_{n \to \infty} a_n = \inf I \,\,\, \text{and} \,\,\, \lim_{n \to \infty} b_n = \sup I.
$$
Then we have
$$
\lim_{n \to \infty} |a_n - b_n| = \beta = \begin{cases}
         \infty, & \text{if $I$ is unbounded}, \\
         |I| = \text{the length of $I$}, & \text{if $I$ is bounded}. \\
         \end{cases}
$$ 
\indent For any point $x_0$ in the {\it closure} of $I$ ($x_0$ need not be in $I$) and any positive integer $n$, let $V_n(x_0)$ denote a {\it relatively} open interval in the {\it closure} of $I$ of length $< 1/n$ which contains the point $x_0$.
 
For each $i \ge 1$, let $y_i$ be an $\omega$-limit point of $x_i$ in the {\it closure} of $I$ ($y_i$ need not be in $I$).  For each $n \ge 1$, let $\widehat V_n(y_i)$ be an open interval in $I$ with length $< 1/n$ such that 
$$
\lim_{n \to \infty} dist\big(\widehat V_n(y_i), V_n(y_i)\big) \ge \frac \beta2 = \, \begin{cases}
                    \infty, & \text{if $I$ is unbounded}, \\
                    \frac 12|I|, & \text{if $I$ is bounded}. \\
                    \end{cases}
$$
$$--------------------$$
\indent In the following, we describe the construction of the scrambled sets.

Roughly speaking, 

at the first stage, we start with any 2 disjoint compact intervals $K_0^{(1,0)}$ and $K_1^{(1,0)}$ in $U_1$.  Then we perform $(2 \cdot 1 +1)^2 + 1^2 = 10$ steps to obtain 10 positive integers and $1 \cdot 2^{1+1}$ pairwise disjoint compact intervals $K_{\al_0\al_1}^{(1,10)}$, $\al_i = 0, 1$, $0 \le i \le 1$, in $K_0^{(1,0)} \cup K_1^{(1,0)}$ with some properties such that $U_2 \setminus \bigcup \, \{ K_{\al_0\al_1}^{(1,10)}: \al_i = 0, 1$, $0 \le i \le 1 \} \ne \emptyset$,  

at the second stage, we start with any $2^2$ pairwise disjoint compact intervals $K_{\al_0\al_1}^{(2,10)}$, $\al_i = 0, 1$, $0 \le i \le 1$, in $U_2 \setminus \bigcup \, \{ K_{\al_0\al_1}^{(1,10)}: \al_i = 0, 1$, $0 \le i \le 1 \}$, together with the previously obtained $2^2$ pairwise disjoint compact intervals $K_{\al_0\al_1}^{(1,10)}$, $\al_i = 0, 1$, $0 \le i \le 1$.  Then we perform $(2 \cdot 2 +1)^2 + 2^2 = 29$ steps to obtain 29 positive integers and $2 \cdot 2^{2+1}$ pairwise disjoint compact intervals $K_{\al_0\al_1\al_2}^{(j,39)}$, $\al_i = 0, 1$, $0 \le i \le 2$, $1 \le j \le 2$, in $\bigcup \, \{ K_{\al_0\al_1}^{(j,10)}: \al_i = 0,1$, $0 \le i \le 1$, $1 \le j \le 2 \}$ with some properties such that $U_3 \setminus \bigcup \, \{ K_{\al_0\al_1\al_2}^{(j,39)}: \al_i = 0, 1$, $0 \le i \le 2$, $1 \le j \le 2 \} \ne \emptyset$,

we then proceed in this manner inductively.  

Now, we start the proof. 

At the first stage, 

let $n_1$ be any positive integer and let $K_0^{(1,0)}$ and $K_1^{(1,0)}$ be any two disjoint compact intervals in $U_1$ \big(the 1 in the superscript $(1,0)$ indicates that these sets lie in $U_1$\big).

By applying Theorem 6(2) to the following $2 \cdot 2^2$ intervals \big(note the subscripts of $K^{(1,0)}$\big)
$$K_0^{(1,0)}, K_0^{(1,0)}, K_1^{(1,0)}, K_1^{(1,0)}: \,\, V_{n_1}(a_1), V_{n_1}(b_1), V_{n_1}(a_1), V_{n_1}(b_1),$$ 
we obtain a positive integer $k_1 \, (> n_1)$ in $\mathcal M$ and $2^2$ pairwise disjoint compact intervals $K_{00}^{(1,1)}$ and $K_{01}^{(1,1)}$ in $K_{0}^{(1,0)}$, \, $K_{10}^{(1,1)}$ and $K_{11}^{(1,1)}$ in $K_{1}^{(1,0)}$ such that  
{\large 
\begin{multline*}
$$
\text{each length is so smaller than} \,\,\, \frac 1{2^3} \,\,\, \text{that the open set} \\ 
\qquad\quad U_2 \setminus \big(K_{00}^{(1,1)} \cup K_{01}^{(1,1)} \cup K_{10}^{(1,1)} \cup K_{11}^{(1,1)}\big) \,\,\, \text{is nonempty and},
$$
\end{multline*}}
for all $\al_0 = 0, 1$, $\al_1 = 0, 1$, 
{\large 
$$
f^{k_1}\big(K_{\al_0\al_1}^{(1,1)}\big) \subset W_{n_1}(\al_1), \, \text{where} \,\,\, W_{n_1}(0) = V_{n_1}(a_1) \,\,\, \text{and} \,\,\, W_{n_1}(1) = V_{n_1}(b_1).
$$}
\indent By applying Theorem 6(2) to the following $2 \cdot 2^2$ intervals
$$
K_{00}^{(1,1)}, K_{01}^{(1,1)}, K_{10}^{(1,1)}, K_{11}^{(1,1)}: \,\, V_{n_1}(a_1), V_{n_1}(b_1), V_{n_1}(a_1), V_{n_1}(b_1),
$$ 
we obtain a positive integer $N$ such that,
for all $\al_0 = 0, 1$, $\al_1 = 0, 1$, and all $n \ge N$, 
$$
f^n\big(K_{\al_0\al_1}^{(1,1)}\big) \cap W_{n_1}(\al_1) \ne \emptyset, \, \text{where} \,\,\, W_{n_1}(0) = V_{n_1}(a_1) \,\,\, \text{and} \,\,\, W_{n_1}(1) = V_{n_1}(b_1).
$$
\indent Let $k_2 \, (> k_1)$ be a positive integer in $\mathcal M$ such that $k_2 +1 \ge N$ and let $K_{\al_0\al_1}^{(1,2)} \subset K_{\al_0\al_1}^{(1,1)}$, $\al_0 = 0, 1$, $\al_1 = 0, 1$, be $2^2$ pairwise disjoint compact intervals such that,  
for all $\al_0 = 0, 1$, $\al_1 = 0, 1$, 
{\large 
\begin{multline*}
$$
f^{k_1}\big(f(K_{\al_0\al_1}^{(1,2)})\big) = f^{k_1+1}\big(K_{\al_0\al_1}^{(1,2)}\big) \subset W_{n_1}(\al_1), \\ \text{where} \,\,\, W_{n_1}(0) = V_{n_1}(a_1) \,\,\, \text{and} \,\,\, W_{n_1}(1) = V_{n_1}(b_1).
$$
\end{multline*}}
\indent By applying Theorem 6(2) to the following $2 \cdot 2^2$ intervals 
$$
K_{00}^{(1,2)}, K_{01}^{(1,2)}, K_{10}^{(1,2)}, K_{11}^{(1,2)}: \,\, V_{n_1}(v_1), V_{n_1}(v_1), V_{n_1}(v_1), V_{n_1}(v_1),
$$
we obtain a positive integer $k_3 \, (> k_2)$ in $\mathcal M$ and $2^2$ pairwise disjoint compact intervals $K_{\al_0\al_1}^{(1,3)} \subset K_{\al_0\al_1}^{(1,2)}$, $\al_0 = 0, 1$, $\al_1 = 0, 1$, such that 
$$
f^{k_3}\big(K_{00}^{(1,3)} \cup K_{01}^{(1,3)} \cup K_{10}^{(1,3)} \cup K_{11}^{(1,3)}\big) \subset V_{n_1}(v_1).%\vspace{.1in}
$$
\indent By applying Theorem 6(2) to the following $2 \cdot 2^2$ intervals 
$$
K_{00}^{(1,3)}, K_{01}^{(1,3)}, K_{10}^{(1,3)}, K_{11}^{(1,3)}: \,\, U_1, U_1, U_1, U_1,
$$
we obtain a positive integer $k_4 \, (> k_3)$ and $2^2$ pairwise disjoint compact intervals $K_{\al_0\al_1}^{(1,4)} \subset K_{\al_0\al_1}^{(1,3)}$, $\al_0 = 0, 1$, $\al_1 = 0, 1$, such that 
$$
f^{k_3}\big(K_{00}^{(1,4)} \cup K_{01}^{(1,4)} \cup K_{10}^{(1,4)} \cup K_{11}^{(1,4)}\big) \subset U_1 \quad \text{and} \quad 1! \quad \text{divides} \quad k_4.\vspace{.15in}
$$
\indent By applying Theorem 6(2) to the following $2 \cdot 2^2$ intervals 
$$
K_{00}^{(1,4)}, K_{01}^{(1,4)}, K_{10}^{(1,4)}, K_{11}^{(1,4)}: \,\, V_{n_1}(y_1), V_{n_1}(y_1), V_{n_1}(y_1), V_{n_1}(y_1),
$$ 
we obtain a positive integer $N'$ such that 
$$
f^n\big(K_{\al_0\al_1}^{(1,4)}\big) \cap V_{n_1}(y_1) \ne \emptyset \,\,\, \text{for all} \,\,\, \al_0 = 0, 1, \, \al_1 = 0, 1, \, \text{and all} \,\,\, n \ge N'.
$$
\indent Since $y_1$ is an $\omega$-limit point of $x_1$, the open neighborhood $V_{n_1}(y_1)$ of $y_1$ contains $f^i(x_1)$ for infinitely many $i$'s.  Let $k_5$ be one such $i$ which is $> N'+k_4$.  Then $f^{k_5}(x_1) \in V_{n_1}(y_1)$ and $k_5 > N'$.  Let $K_{\al_0\al_1}^{(1,5)} \subset K_{\al_0\al_1}^{(1,4)}$, $\al_0 = 0, 1$, $\al_1 = 0, 1$, be $2^2$ pairwise disjoint compact intervals such that 
$$
f^{k_5}\big(K_{00}^{(1,5)} \cup K_{01}^{(1,5)} \cup K_{10}^{(1,5)} \cup K_{11}^{(1,5)}\big) \subset V_{n_1}(y_1) \,\,\, \big(\text{and} \,\,\, f^{k_5}(x_1) \in V_{n_1}(y_1)\big).
$$
\indent By applying Theorem 6(2) to the following $2 \cdot 2^2$ intervals 
$$
K_{00}^{(1,5)}, K_{01}^{(1,5)}, K_{10}^{(1,5)}, K_{11}^{(1,5)}: \,\, V_{n_1}(y_1), V_{n_1}(y_1), V_{n_1}(y_1), V_{n_1}(y_1),
$$
we obtain a positive integer $N''$ such that 
$$
f^n\big(K_{\al_0\al_1}^{(1,5)}\big) \cap V_{n_1}(y_1) \ne \emptyset \,\,\, \text{for all} \,\,\, \al_0 = 0, 1, \, \al_1 = 0, 1, \, \text{and all} \,\,\, n \ge N''.
$$\indent Since $y_1$ is an $\omega$-limit point of $x_1$, the open neighborhood $V_{n_1}(y_1)$ of $y_1$ contains $f^i(x_1)$ for infinitely many $i$'s.  Let $k_6$ be one such $i$ which is $> N''+k_5$.  Then $f^{k_6}(x_1) \in V_{n_1}(y_1)$ and $k_6+1 \ge N''$.  Let $K_{\al_0\al_1}^{(1,6)} \subset K_{\al_0\al_1}^{(1,5)}$, $\al_0 = 0, 1$, $\al_1 = 0, 1$, be $2^2$ pairwise disjoint compact intervals such that 
\begin{multline*}
$$
f^{k_6}\big(f\big(K_{00}^{(1,6)}\big) \cup f\big(K_{01}^{(1,6)}\big) \cup f\big(K_{10}^{(1,6)}\big) \cup f\big(K_{11}^{(1,6)}\big)\big) \\ = f^{k_6+1}\big(K_{00}^{(1,6)}\big) \cup f^{k_6+1}\big(K_{01}^{(1,6)}\big) \cup f^{k_6+1}\big(K_{10}^{(1,6)}\big) \cup f^{k_6+1}\big(K_{11}^{(1,6)}\big)\subset V_{n_1}(y_1) \,\,\, \big(\text{and} \,\,\, f^{k_6}(x_1) \in V_{n_1}(y_1)\big).
$$
\end{multline*}
\indent By arguments similar to those in the previous two steps with $V_{n_1}(y_1)$ replaced by $\widehat V_{n_1}(y_1)$, we obtain successively two integers $k_7$ $(> k_6)$ and $k_8 > k_7$ and two $2^2$ pairwise disjoint compact intervals $K_{\al_0\al_1}^{(1,7)} \subset K_{\al_0\al_1}^{(1,6)}$, $\al_0 = 0, 1$, $\al_1 = 0, 1$, and $K_{\al_0\al_1}^{(1,8)} \subset K_{\al_0\al_1}^{(1,7)}$, $\al_0 = 0, 1$, $\al_1 = 0, 1$, such that 
$$
f^{k_7}\big(K_{00}^{(1,7)} \cup K_{01}^{(1,7)} \cup K_{10}^{(1,7)} \cup K_{11}^{(1,7)}\big) \subset \widehat V_{n_1}(y_1) \,\,\, \text{and} \,\,\, f^{k_7}(x_1) \in V_{n_1}(y_1),
$$and
$$
f^{k_8}\big(f\big(K_{00}^{(1,8)}\big) \cup f\big(K_{01}^{(1,8)}\big) \cup f\big(K_{10}^{(1,8)}\big) \cup f\big(K_{11}^{(1,8)}\big)\big) \subset \widehat V_{n_1}(y_1) \,\,\, \text{and} \,\,\, f^{k_8}(x_1) \in V_{n_1}(y_1).%\vspace{.1in}
$$
\indent On the other hand, by applying Theorem 6(2) to the following $2 \cdot 2^2$ intervals
$$
K_{00}^{(1,8)}, K_{01}^{(1,8)}, K_{10}^{(1,8)}, K_{11}^{(1,8)}: \,\, V_{n_1}(y_1), V_{n_1}(y_1), V_{n_1}(y_1), V_{n_1}(y_1),
$$ 
we obtain a positive integer $N'''$ such that 
$$
f^n\big(K_{\al_0\al_1}^{(1,8)}\big) \cap V_{n_1}(y_1) \ne \emptyset \,\,\, \text{for all} \,\,\, \al_0 = 0, 1, \, \al_1 = 0, 1, \, \text{and all} \,\,\, n \ge N'''.
$$
\indent Since $y_1$ is an $\omega$-limit point of $x_1$, the open neighborhood $V_{n_1}(y_1)$ of $y_1$ contains $f^i(x_1)$ for infinitely many $i$'s.  Let $k_9+1$ be one such $i$ which is $> N'''+k_8+1$.  Then $f^{k_9+1}(x_1) \in V_{n_1}(y_1)$ and $k_9 \ge N'''$.  Let $K_{\al_0\al_1}^{(1,9)} \subset K_{\al_0\al_1}^{(1,8)}$, $\al_0 = 0, 1$, $\al_1 = 0, 1$, be $2^2$ pairwise disjoint compact intervals such that 
$$
f^{k_8}\big(K_{00}^{(1,9)} \cup K_{01}^{(1,9)} \cup K_{10}^{(1,9)} \cup K_{11}^{(1,9)}\big) \subset V_{n_1}(y_1) \,\,\, \text{and} \,\,\, f^{k_9}\big(f(x_1)\big) = f^{k_9+1}(x_1) \in V_{n_1}(y_1).
$$
\indent By arguments similar to those in the previous step with {\it appropriate} $V_{n_1}(y_1)$ replaced by $\widehat V_{n_1}(y_1)$, we obtain a positive integer $k_{10} \, (> k_9)$ and $2^2$ pairwise disjoint compact intervals $K_{\al_0\al_1}^{(1,10)} \subset K_{\al_0\al_1}^{(1,9)}$, $\al_0 = 0, 1$, $\al_1 = 0, 1$, such that
$$
f^{k_9}\big(K_{00}^{(1,10)} \cup K_{01}^{(1,10)} \cup K_{10}^{(1,10)} \cup K_{11}^{(1,10)}\big) \subset \widehat V_{n_1}(y_1) \,\,\, \text{and} \,\,\, f^{k_{10}}\big(f(x_1)\big) = f^{k_{10}+1}(x_1) \in V_{n_1}(y_1).
$$

Since, for each $\al_0, \al_1 = 0, 1$, we have $K_{\al_0}^{(1,0)} \supset K_{\al_0\al_1}^{(1,1)} \supset K_{\al_0\al_1}^{(1,2)} \supset \cdots \supset K_{\al_0\al_1}^{(1,8)} \supset K_{\al_0\al_1}^{(1,9)} \supset K_{\al_0\al_1}^{(1,10)}$,
we can summarize the arguments in the (above) first stage as follows:

{\footnotesize \[
\left(\begin{array}{l}
\textrm{Let $n_1$ be any positive integer and let $K_0^{(1,0)}$ and $K_1^{(1,0)}$ any two disjoint compact intervals in $U_1$.} \\
\textrm{\qquad By applying Theorem 6(2) to the following $2 \cdot 2^2$ intervals}\\
\textrm{\qquad\qquad\qquad\quad $K_0^{(1,0)}, \, K_0^{(1,0)}, \, K_1^{(1,0)}, \, K_1^{(1,0)}; \, V_{n_1}(a_1), \, V_{n_1}(b_1), \, V_{n_1}(a_1), \, V_{n_1}(b_1)$,} \\ 
\textrm{we obtain a positive integer $k_1$ $(> n_1)$ in $\mathcal M$ and $2^2$ pairwise disjoint compact intervals}\\
\textrm{\qquad\qquad\qquad\quad $K_{00}^{(1,1)}$ and $K_{01}^{(1,1)}$ in $K_0^{(1,0)}$ and $K_{10}^{(1,1)}$ and $K_{11}^{(1,1)}$ in $K_1^{(1,0)}$ such that}\\
\textrm{\qquad each length is so smaller than $1/2^3$ that}\\
\textrm{\qquad\qquad\qquad\quad\qquad\qquad\quad the open set $U_2 \setminus \big(K_{00}^{(1,1)} \cup K_{01}^{(1,1)} \cup K_{10}^{(1,1)} \cup K_{11}^{(1,1)}\big)$ is nonempty and,}\\
\textrm{for all $\al_0 = 0, 1$, $\al_1 = 0, 1$, we have} \\
\textrm{\qquad\qquad\qquad\quad $f^{k_1}\big(K_{\al_0\al_1}^{(1,1)}\big) \subset W_{n_1}(\al_1), \, \text{where} \,\,\, W_{n_1}(0) = V_{n_1}(a_1) \,\,\, \text{and} \,\,\, W_{n_1}(1) = V_{n_1}(b_1)$.}\\
\textrm{Start with $K_{\al_0\al_1}^{(1,1)}$, $\al_0=0, 1$, $\al_1=0, 1$, by applying Theorem 6(2) successively and appropriately,}\\
\textrm{we obtain $(2 \cdot 1+1)^2 +1^2 = 10$ positive integers $(n_1 <)$ $k_1 < k_2 < k_3 < \cdots < k_8 < k_9 < k_{10}$ and}\\
\textrm{\qquad\quad\quad $2^2$ pairwise disjoint compact intervals $K_{\al_0\al_1}^{(1,10)}$, $\al_0= 0, 1$, $\al_1=0, 1$, such that}\\
\textrm{\qquad\quad\qquad\quad\qquad $\{ k_1, k_2, k_3 \} \subset \mathcal M$, $1!$ divides $k_4$ and} \\
\textrm{\qquad\quad\qquad\quad\qquad $f^{k_1}\big(K_{\al_0\al_1}^{(1,10)}\big) \subset W_{n_1}(\al_1)$, where $W_{n_1}(0) = V_{n_1}(a_1)$ and $W_{n_1}(1) = V_{n_1}(b_1)$,}\\
\textrm{\qquad\quad\qquad\quad\qquad $f^{k_2}\big(f(K_{\al_0\al_1}^{(1,10)})\big) \subset W_{n_1}(\al_1)$, where $W_{n_1}(0) = V_{n_1}(a_1)$ and $W_{n_1}(1) = V_{n_1}(b_1)$,}\\
\textrm{\qquad\quad\qquad\quad\qquad $f^{k_3}\big(K_{\al_0\al_1}^{(1,10)}\big) \subset V_{n_1}(v_1)$,}\\
\textrm{\qquad\quad\qquad\quad\qquad $f^{k_4}\big(K_{\al_0\al_1}^{(1,10)}\big) \subset U_1$,}\\
\textrm{\qquad\quad\qquad\quad\qquad $f^{k_5}\big(K_{\al_0\al_1}^{(1,10)}\big) \subset V_{n_1}(y_1)$ \,\,\,\,\, and \,\,\,\,\, $f^{k_5}(x_1) \in V_{n_1}(y_1)$,}\\
\textrm{\qquad\quad\qquad\quad\qquad $f^{k_6}\big(f(K_{\al_0\al_1}^{(1,10)})\big) \subset V_{n_1}(y_1)$ \, and \, $f^{k_6}(x_1) \in V_{n_1}(y_1)$,}\\
\textrm{\qquad\quad\qquad\quad\qquad $f^{k_7}\big(K_{\al_0\al_1}^{(1,10)}\big) \subset \widehat V_{n_1}(y_1)$ \,\,\,\,\,\, and \,\,\,\, $f^{k_7}(x_1) \in V_{n_1}(y_1)$, and}\\ 
\textrm{\qquad\quad\qquad\quad\qquad $f^{k_8}\big(f(K_{\al_0\al_1}^{(1,10)})\big) \subset \widehat V_{n_1}(y_1)$ \, and \, $f^{k_8}(x_1) \in V_{n_1}(y_1)$,}\\
\textrm{\qquad\quad\qquad\quad\qquad $f^{k_9}\big(K_{\al_0\al_1}^{(1,10)}\big) \subset V_{n_1}(y_1)$ \,\,\,\,\, and \,\,\,\,\, $f^{k_9}\big(f(x_1)\big) \in V_{n_1}(y_1)$,}\\
\textrm{\qquad\quad\qquad\quad\qquad $f^{k_{10}}\big(K_{\al_0\al_1}^{(1,10)}\big) \subset \widehat V_{n_1}(y_1)$ \,\,\,\, and \,\,\,\,\, $f^{k_{10}}\big(f(x_1)\big) \in V_{n_1}(y_1)$.}\\
\end{array}\right)  
\]}

At the second stage, 

let $n_2$ be any positive integer such that $n_2 > k_{10}$.  Since the open set, as indicated above,
$$
U_2 \setminus \big(K_{00}^{(1,10)} \cup K_{01}^{(1,10)} \cup K_{10}^{(1,10)} \cup K_{11}^{(1,10)}\big) \, \supset \, U_2 \setminus \big(K_{00}^{(1,1)} \cup K_{01}^{(1,1)} \cup K_{10}^{(1,1)} \cup K_{11}^{(1,1)}\big) \,\,\, \text{is nonempty},
$$
let $K_{00}^{(2,10)}, K_{01}^{(2,10)}, K_{10}^{(2,10)}, K_{11}^{(2,10)}$ \big(the 2 in the superscript $(2, 10)$ indicates that these sets lie in $U_2$ and the 10 in $(2, 10)$ is to keep up with the 10 in $(1, 10)$ above\big) be any $2^2$ pairwise disjoint compact intervals in 
$$
U_2 \setminus \big( K_{00}^{(1,10)} \cup K_{01}^{(1,10)} \cup K_{10}^{(1,10)} \cup K_{11}^{(1,10)}\big).
$$
\indent In the following, let $\big(K_{00}^{(1,10)}\big)^2$ denote $K_{00}^{(1,10)}, K_{00}^{(1,10)}$ and let $\big(V_{n_2}(a_2), V_{n_2}(b_2)\big)^8$ denote 8 copies of $\big(V_{n_2}(a_2), V_{n_2}(b_2)\big)$ and so on.  

By applying Theorem 6(2) to the following $2 \cdot 2^4$ intervals \big(note the subscripts of $K^{(j,10)}, \, j= 1, 2$\big) 
$$
\big(K_{\al_0\al_1}^{(1,10)}\big)^2, \, \al_i = 0, 1, \, 0 \le i \le 1; \, \big(K_{\al_0\al_1}^{(2,10)}\big)^2, \, \al_i = 0, 1, \, 0 \le i \le 1 \,\, : \,\, \big(V_{n_2}(a_2), V_{n_2}(b_2)\big)^8,
$$
\begin{footnotesize}
\big(i.e., $\big(K_{00}^{(1,10)}\big)^2, \big(K_{01}^{(1,10)}\big)^2, \big(K_{10}^{(1,10)}\big)^2, \big(K_{11}^{(1,10)}\big)^2$; $\big(K_{00}^{(2,10)}\big)^2, \big(K_{01}^{(2,10)}\big)^2, \big(K_{10}^{(2,10)}\big)^2, \big(K_{11}^{(2,10)}\big)^2$ : $\big(V_{n_2}(a_2), V_{n_2}(b_2)\big)^8$\big)
\end{footnotesize}
\big(here we consider not only the $2^2$ just found compact intervals $K_{00}^{(1,10)}, K_{01}^{(1,10)}, K_{10}^{(1,10)}, K_{11}^{(1,10)}$ from $U_1$, we also consider the $2^2$ newly chosen compact intervals $K_{00}^{(2,10)}, K_{01}^{(2,10)}, K_{10}^{(2,10)}, K_{11}^{(2,10)}$ from $U_2$\big), we obtain a positive integer $\big(k_{[(2-1)(10 \cdot 2^2+7 \cdot 2+6)/6]+1} =) \, k_{11} \, (> n_2 > k_{10}\big)$ in $\mathcal M$ and $2 \cdot 2^3$ pairwise disjoint compact intervals $K_{\al_0\al_1\al_2}^{(1,11)} \subset K_{\al_0\al_1}^{(1,10)}, \, \al_i = 0, 1, \, 0 \le i \le 2; \,\, K_{\al_0\al_1\al_2}^{(2,11)} \subset K_{\al_0\al_1}^{(2,10)}, \, \al_i = 0, 1, \, 0 \le i \le 2$, such that 
\begin{multline*}
$$
\text{each length is so smaller than} \,\,\, \frac 1{2^5} \,\,\, \text{that the open set} \\  
\qquad U_3 \setminus \big( \bigcup \big\{ K_{\al_0\al_1\al_2}^{(j,11)}: \, \al_i = 0, 1, \, 0 \le i \le 2, \, j = 1, 2 \big\} \big) \,\,\, \text{is nonempty and}
$$
\end{multline*}
\noindent for all $\al_i = 0, 1, \, 0 \le i \le 2$, 
$$
f^{k_{11}}\big(K_{\al_0\al_1\al_2}^{(1,11)} \cup K_{\al_0\al_1\al_2}^{(2,11)}\big) \subset W_{n_2}(\al_2), \,\, \text{where} \,\,\, W_{n_2}(0) = V_{n_2}(a_2) \,\,\, \text{and} \,\,\, W_{n_2}(1) = V_{n_2}(b_2).\vspace{.1in}
$$
\indent By applying Theorem 6(2) to the following $2(2 \cdot 2^3)$ intervals 
$$
K_{\al_0\al_1\al_2}^{(1,11)}, \, \al_i = 0, 1, \, 0 \le i \le 2, \,\, K_{\al_0\al_1\al_2}^{(2,11)}, \, \al_i = 0, 1, \, 0 \le i \le 2 \,\, : \,\, \big(V_{n_2}(a_2), V_{n_2}(b_2)\big)^8,
$$
$$
\big(\text{i.e.,} \,\,\, K_{000}^{(1,11)}, K_{001}^{(1,11)}, K_{010}^{(1,11)}, K_{011}^{(1,11)}, K_{100}^{(1,11)}, K_{101}^{(1,11)}, K_{110}^{(1,11)}, K_{111}^{(1,11)}; \,\, K_{000}^{(2,11)}, K_{001}^{(2,11)}, \cdots \,\,\, \text{and so on}\big)
$$
and by arguments similar to those in the second step of the first stage, we obtain a positive integer $k_{12} \, (> k_{11})$ in $\mathcal M$ and $2 \cdot 2^3$ pairwise disjoint compact intervals $K_{\al_0\al_1\al_2}^{(1,12)} \subset K_{\al_0\al_1\al_2}^{(1,11)}, \, \al_i = 0, 1, \, 0 \le i \le 2, \, K_{\al_0\al_1\al_2}^{(2,12)} \subset K_{\al_0\al_1\al_2}^{(2,11)}, \, \al_i = 0, 1, \, 0 \le i \le 2$, such that 
\begin{multline*}
$$
f^{k_{12}}\big(f(K_{\al_0\al_1\al_2}^{(1,12)}) \cup f(K_{\al_0\al_1\al_2}^{(2,12)})\big) = f^{k_{12}+1}\big(K_{\al_0\al_1\al_2}^{(1,12)} \cup K_{\al_0\al_1\al_2}^{(2,12)}\big) \subset W_{n_2}(\al_2), \\ \text{where} \,\,\, W_{n_2}(0) = V_{n_2}(a_2) \,\,\, \text{and} \,\,\, W_{n_2}(1) = V_{n_2}(b_2).\vspace{.15in}
$$
\end{multline*}
\indent Similarly, by applying Theorem 6(2) to the following $2(2 \cdot 2^3)$ intervals 
$$
K_{\al_0\al_1\al_2}^{(1,12)}, \, \al_i = 0, 1, \, 0 \le i \le 2, \,\, K_{\al_0\al_1\al_2}^{(2,12)}, \, \al_i = 0, 1, \, 0 \le i \le 2 \,\, : \,\, \big(V_{n_2}(a_2), V_{n_2}(b_2)\big)^8,
$$
we obtain a positive integer $k_{13} \, (> k_{12})$ in $\mathcal M$ and $2 \cdot 2^3$ pairwise disjoint compact intervals $K_{\al_0\al_1\al_2}^{(1,13)} \subset K_{\al_0\al_1\al_2}^{(1,12)}, \, \al_i = 0, 1, \, 0 \le i \le 2, \, K_{\al_0\al_1\al_2}^{(2,13)} \subset K_{\al_0\al_1\al_2}^{(2,12)}, \, \al_i = 0, 1, \, 0 \le i \le 2$, such that 
\begin{multline*}
$$
f^{k_{13}}\big(f^2(K_{\al_0\al_1\al_2}^{(1,13)}) \cup f^2(K_{\al_0\al_1\al_2}^{(2,13)})\big) = f^{k_{12}+2}\big(K_{\al_0\al_1\al_2}^{(1,13)} \cup K_{\al_0\al_1\al_2}^{(2,13)}\big) \subset W_{n_2}(\al_2), \\ \text{where} \,\,\, W_{n_2}(0) = V_{n_2}(a_2) \,\,\, \text{and} \,\,\, W_{n_2}(1) = V_{n_2}(b_2).\vspace{.1in}
$$
\end{multline*}
\indent By applying Theorem 6(2) to the following $2(2 \cdot 2^3)$ intervals 
$$
K_{\al_0\al_1\al_2}^{(1,13)}, \, \al_i = 0, 1, \, 0 \le i \le 2, \,\, K_{\al_0\al_1\al_2}^{(2,13)}, \, \al_i = 0, 1, \, 0 \le i \le 2 \,\, : \,\, \big(V_{n_2}(a_2)\big)^8, \,\,\, \big(V_{n_2}(b_2)\big)^8,
$$
we obtain a positive integer $k_{14} \, (> k_{13})$ in $\mathcal M$ and $2 \cdot 2^3$ pairwise disjoint compact intervals $K_{\al_0\al_1\al_2}^{(1,14)} \subset K_{\al_0\al_1\al_2}^{(1,13)}, \, \al_i = 0, 1, \, 0 \le i \le 2, \, K_{\al_0\al_1\al_2}^{(2,14)} \subset K_{\al_0\al_1\al_2}^{(2,13)}, \, \al_i = 0, 1, \, 0 \le i \le 2$, such that 
$$
f^{k_{14}}\big(\cup_{\al_i = 0, 1, \, 0 \le i \le 2} \, K_{\al_0\al_1\al_2}^{(1,14)}\big) \subset V_{n_2}(a_2) \qquad \text{and} \qquad f^{k_{14}}\big(\cup_{\al_i = 0, 1, \, 0 \le i \le 2} K_{\al_0\al_1\al_2}^{(2,14)}\big) \subset V_{n_2}(b_2).
$$
\indent By applying Theorem 6(2) to the following $2(2 \cdot 2^3)$ intervals 
$$
K_{\al_0\al_1\al_2}^{(1,14)}, \, \al_i = 0, 1, \, 0 \le i \le 2, \,\, K_{\al_0\al_1\al_2}^{(2,14)}, \, \al_i = 0, 1, \, 0 \le i \le 2 \,\, : \,\, \big(V_{n_2}(a_2)\big)^8, \,\,\, \big(V_{n_2}(b_2)\big)^8,
$$
we obtain a positive integer $k_{15} \, (> k_{14})$ in $\mathcal M$ and $2 \cdot 2^3$ pairwise disjoint compact intervals $K_{\al_0\al_1\al_2}^{(1,15)} \subset K_{\al_0\al_1\al_2}^{(1,14)}, \, \al_i = 0, 1, \, 0 \le i \le 2, \, K_{\al_0\al_1\al_2}^{(2,15)} \subset K_{\al_0\al_1\al_2}^{(2,14)}, \, \al_i = 0, 1, \, 0 \le i \le 2$, such that 
\begin{multline*}
$$
f^{k_{15}}\big(\cup_{\al_i = 0, 1, \, 0 \le i \le 2} \, f\big(K_{\al_0\al_1\al_2}^{(1,15)}\big)\big) = f^{k_{15}+1}\big(\cup_{\al_i = 0, 1, \, 0 \le i \le 2} \, K_{\al_0\al_1\al_2}^{(1,15)}\big) \subset V_{n_2}(a_2) \qquad \text{and} \\
\qquad f^{k_{15}}\big(\cup_{\al_i = 0, 1, \, 0 \le i \le 2} \, f\big(K_{\al_0\al_1\al_2}^{(2,15)}\big)\big) = f^{k_{15}+1}\big(\cup_{\al_i = 0, 1, \, 0 \le i \le 2} K_{\al_0\al_1\al_2}^{(2,15)}\big) \subset V_{n_2}(b_2).
$$
\end{multline*}
\indent By applying Theorem 6(2) to the following $2(2 \cdot 2^3)$ intervals 
$$
K_{\al_0\al_1\al_2}^{(1,15)}, \, \al_i = 0, 1, \, 0 \le i \le 2, \,\, K_{\al_0\al_1\al_2}^{(2,15)}, \, \al_i = 0, 1, \, 0 \le i \le 2 \,\, : \,\, \big(V_{n_2}(a_2)\big)^8, \,\,\, \big(V_{n_2}(b_2)\big)^8,
$$
we obtain a positive integer $k_{16} \, (> k_{15})$ in $\mathcal M$ and $2 \cdot 2^3$ pairwise disjoint compact intervals $K_{\al_0\al_1\al_2}^{(1,16)} \subset K_{\al_0\al_1\al_2}^{(1,15)}, \, \al_i = 0, 1, \, 0 \le i \le 2, \, K_{\al_0\al_1\al_2}^{(2,16)} \subset K_{\al_0\al_1\al_2}^{(2,15)}, \, \al_i = 0, 1, \, 0 \le i \le 2$, such that 
\begin{multline*}
$$
f^{k_{16}}\big(\cup_{\al_i = 0, 1, \, 0 \le i \le 2} \, f^2\big(K_{\al_0\al_1\al_2}^{(1,16)}\big)\big) = f^{k_{15}+2}\big(\cup_{\al_i = 0, 1, \, 0 \le i \le 2} \, K_{\al_0\al_1\al_2}^{(1,16)}\big) \subset V_{n_2}(a_2) \qquad \text{and} \\
\qquad f^{k_{16}}\big(\cup_{\al_i = 0, 1, \, 0 \le i \le 2} \, f^2\big(K_{\al_0\al_1\al_2}^{(2,16)}\big)\big) = f^{k_{16}+2}\big(\cup_{\al_i = 0, 1, \, 0 \le i \le 2} K_{\al_0\al_1\al_2}^{(2,16)}\bigg) \subset V_{n_2}(b_2).\vspace{.15in}
$$
\end{multline*}
\indent By applying Theorem 6(2) to the following $2(2 \cdot 2^3)$ intervals 
$$
K_{\al_0\al_1\al_2}^{(1,16)}, \, \al_i = 0, 1, \, 0 \le i \le 2, \,\, K_{\al_0\al_1\al_2}^{(2,16)}, \, \al_i = 0, 1, \, 0 \le i \le 2 \,\, : \,\, \big(V_{n_2}(v_1)\big)^{16},
$$
we obtain a positive integer $k_{17} \, (> k_{16})$ in $\mathcal M$ and $2 \cdot 2^3$ pairwise disjoint compact intervals $K_{\al_0\al_1\al_2}^{(1,17)} \subset K_{\al_0\al_1\al_2}^{(1,16)}, \, \al_i = 0, 1, \, 0 \le i \le 2, \, K_{\al_0\al_1\al_2}^{(2,17)} \subset K_{\al_0\al_1\al_2}^{(2,16)}, \, \al_i = 0, 1, \, 0 \le i \le 2$, such that 
$$
f^{k_{17}}\big(\cup \big\{ K_{\al_0\al_1\al_2}^{(j,17)}: \, \al_i = 0, 1, \, 0 \le i \le 2, \, j = 1, 2 \big\}\big) \, \subset V_{n_2}(v_1).
$$
\indent By applying Theorem 6(2) to the following $2(2 \cdot 2^3)$ intervals 
$$
K_{\al_0\al_1\al_2}^{(1,17)}, \, \al_i = 0, 1, \, 0 \le i \le 2, \,\, K_{\al_0\al_1\al_2}^{(2,17)}, \, \al_i = 0, 1, \, 0 \le i \le 2 \,\, : \,\, \big(V_{n_2}(v_2)\big)^{16},
$$
we obtain a positive integer $k_{18} \, (> k_{17})$ in $\mathcal M$ and $2 \cdot 2^3$ pairwise disjoint compact intervals $K_{\al_0\al_1\al_2}^{(1,18)} \subset K_{\al_0\al_1\al_2}^{(1,17)}, \, \al_i = 0, 1, \, 0 \le i \le 2, \, K_{\al_0\al_1\al_2}^{(2,18)} \subset K_{\al_0\al_1\al_2}^{(2,17)}, \, \al_i = 0, 1, \, 0 \le i \le 2$, such that 
$$
f^{k_{18}}\big(\cup \big\{ K_{\al_0\al_1\al_2}^{(j,18)}: \, \al_i = 0, 1, \, 0 \le i \le 2, \, j = 1, 2 \big\}\big) \, \subset V_{n_2}(v_2).
$$
\indent By applying Theorem 6(2) to the following $2(2 \cdot 2^3)$ intervals 
$$
K_{\al_0\al_1\al_2}^{(1,18)}, \, \al_i = 0, 1, \, 0 \le i \le 2, \,\, K_{\al_0\al_1\al_2}^{(2,18)}, \, \al_i = 0, 1, \, 0 \le i \le 2 \,\, : \,\, \big(U_2\big)^{16},
$$
we obtain a positive integer $k_{19} \, (> k_{18})$ and $2 \cdot 2^3$ pairwise disjoint compact intervals $K_{\al_0\al_1\al_2}^{(1,19)} \subset K_{\al_0\al_1\al_2}^{(1,18)}, \, \al_i = 0, 1, \, 0 \le i \le 2, \, K_{\al_0\al_1\al_2}^{(2,19)} \subset K_{\al_0\al_1\al_2}^{(2,18)}, \, \al_i = 0, 1, \, 0 \le i \le 2$, such that 
$$
f^{k_{19}}\big(\cup \big\{ K_{\al_0\al_1\al_2}^{(j,19)}: \, \al_i = 0, 1, \, 0 \le i \le 2, \, j = 1, 2 \big\}\big) \, \subset U_2 \quad \text{and} \quad 2! \,\,\, \text{divides} \,\,\, k_{19}. \vspace{.1in}
$$           
\indent By applying Theorem 6(2) to the following $2(2 \cdot 2^3)$ intervals 
$$
K_{\al_0\al_1\al_2}^{(j,19)}, \, \al_i = 0, 1, \, 0 \le i \le 2, \, j = 1, 2 \,\, : \,\, \big(V_{n_2}(y_1)\big)^{16},
$$
we obtain a positive integer $N^*$ such that 
$$
f^n\big(K_{\al_0\al_1\al_2}^{(j,19)}\big) \cap V_{n_2}(y_1) \ne \emptyset \,\,\, \text{for all} \,\,\, \al_i = 0, 1, \, 0 \le i \le 2, \, j = 1, 2, \, \text{and all} \,\,\, n \ge N^*.
$$
\indent Since $y_1$ is an $\omega$-limit point of $x_1$, the open neighborhood $V_{n_2}(y_1)$ of $y_1$ contains $f^i(x_1)$ for infinitely many $i$'s.  Let $k_{20}$ be one such $i$ which is $> N^*+k_{19}$.  Then $f^{k_{20}}(x_1) \in V_{n_2}(y_1)$ and $k_{20} > N*$. Let $K_{\al_0\al_1\al_2}^{(j,20)} \subset K_{\al_0\al_1\al_2}^{(j,19)}$, $\al_0 = 0, 1$, $\al_1 = 0, 1$, $j = 1, 2$, be $2 \cdot 2^3$ pairwise disjoint compact intervals such that 
$$
f^{k_{20}}\big(\cup \big\{ K_{\al_0\al_1\al_2}^{(j,20)}: \, \al_i = 0, 1, \, 0 \le i \le 2, \, j = 1, 2 \big\}\big) \, \subset V_{n_2}(y_1) \,\,\, \big(\text{and} \,\,\, f^{k_{20}}(x_1) \in V_{n_2}(y_1)\big).
$$
\indent By applying Theorem 6(2) to the following $2(2 \cdot 2^3)$ intervals 
$$
K_{\al_0\al_1\al_2}^{(j,20)}, \, \al_i = 0, 1, \, 0 \le i \le 2, \, j = 1, 2 \,\, : \,\, \big(V_{n_2}(y_1)\big)^{16},
$$
we obtain a positive integer $N^{**}$ such that 
$$
f^n\big(K_{\al_0\al_1\al_2}^{(j,20)}\big) \cap V_{n_2}(y_1) \ne \emptyset \,\,\, \text{for all} \,\,\, \al_i = 0, 1, \, 0 \le i \le 2, \, j = 1, 2, \, \text{and all} \,\,\, n \ge N^{**}. 
$$
\indent Since $y_1$ is an $\omega$-limit point of $x_1$, the open neighborhood $V_{n_2}(y_1)$ of $y_1$ contains $f^i(x_1)$ for infinitely many $i$'s.  Let $k_{21}$ be one such $i$ which is $> N^{**}+k_{20}$.  Then $f^{k_{21}}(x_1) \in V_{n_2}(y_1)$ and $k_{20}+1 > N^{**}$. Let $K_{\al_0\al_1\al_2}^{(j,21)} \subset K_{\al_0\al_1\al_2}^{(j,20)}$, $\al_0 = 0, 1$, $\al_1 = 0, 1$, $j = 1, 2$, be $2 \cdot 2^3$ pairwise disjoint compact intervals such that 
\begin{multline*}
$$
f^{k_{21}}\big(\cup \big\{ f\big(K_{\al_0\al_1\al_2}^{(j,21)}\big): \, \al_i = 0, 1, \, 0 \le i \le 2, \, j = 1, 2 \big\}\big) \\
= f^{k_{21}+1}\big(\cup \big\{ K_{\al_0\al_1\al_2}^{(j,21)}: \, \al_i = 0, 1, \, 0 \le i \le 2, \, j = 1, 2 \big\}\big) \subset V_{n_2}(y_1) \\ 
\big(\text{and} \,\,\, f^{k_{21}}(x_1) \in V_{n_2}(y_1)\big).\qquad\qquad\qquad\qquad\qquad\qquad\qquad\qquad
$$
\end{multline*}
\indent Similarly, we obtain a positive integer $k_{22} \, (> k_{21})$ and $2 \cdot 2^3$ pairwise disjoint compact intervals $K_{\al_0\al_1\al_2}^{(j,22)} \subset K_{\al_0\al_1\al_2}^{(j,21)}$, $\al_i = 0, 1$, $0 \le i \le 2$, $j = 1, 2$, such that 
$$
f^{k_{22}}\big(\cup \big\{ f^2\big(K_{\al_0\al_1\al_2}^{(j,22)}\big): \, \al_i = 0, 1, \, 0 \le i \le 2, \, j = 1, 2 \big\}\big) \subset V_{n_2}(y_1) \,\,\, \text{and} \,\,\, f^{k_{22}}(x_1) \in V_{n_2}(y_1).%\vspace{.1in}
$$
\indent By arguments similar to those in the previous three steps with $V_{n_2}(y_1)$ replaced by $\widehat V_{n_2}(y_1)$, we obtain successively 
\begin{multline*}
$$
\text{an integer} \,\,\, k_{23} \, (> k_{22}) \,\,\, \text{and} \,\,\, 2 \cdot 2^3 \,\,\, \text{pairwise disjoint compact intervals} \\      K_{\al_0\al_1\al_2}^{(j,23)} \subset K_{\al_0\al_1\al_2}^{(j,22)}, \, \al_i = 0, 1, \, 0 \le i \le 2, \, j = 1, 2, \, \text{and} \\
\text{an integer} \,\,\, k_{24} \, (> k_{23}) \,\,\, \text{and} \,\,\, 2 \cdot 2^3 \,\,\, \text{pairwise disjoint compact intervals}\qquad\qquad\qquad\qquad\qquad\quad\,\,\, \\ 
K_{\al_0\al_1\al_2}^{(j,24)} \subset K_{\al_0\al_1\al_2}^{(j,23)}, \, \al_i = 0, 1, \, 0 \le i \le 2, \, j = 1, 2, \, \text{and}  \\
\text{an integer} \,\,\, k_{25} \, (> k_{24}) \,\,\, \text{and} \,\,\, 2 \cdot 2^3 \,\,\, \text{pairwise disjoint compact intervals}\qquad\qquad\qquad\qquad\qquad\quad\,\,\, \\ 
K_{\al_0\al_1\al_2}^{(j,25)} \subset K_{\al_0\al_1\al_2}^{(j,24)}, \, \al_i = 0, 1, \, 0 \le i \le 2, \, j = 1, 2, \qquad\qquad\qquad\qquad\qquad\quad
$$
\end{multline*}
such that
$$
f^{k_{23}}\big(\cup \big\{ K_{\al_0\al_1\al_2}^{(j,23)}: \, \al_i = 0, 1, \, 0 \le i \le 2, \, j = 1, 2 \big\}\big) \, \subset \widehat V_{n_2}(y_1) \,\,\, \text{and} \,\,\, f^{k_{23}}(x_1) \in V_{n_2}(y_1),
$$
$$
f^{k_{24}}\big(\cup \big\{ f\big(K_{\al_0\al_1\al_2}^{(j,24)}\big): \, \al_i = 0, 1, \, 0 \le i \le 2, \, j = 1, 2 \big\}\big) \subset \widehat V_{n_2}(y_1) \,\,\, \text{and} \,\,\, f^{k_{14}}(x_1) \in V_{n_2}(y_1),\, \text{and}
$$
$$
f^{k_{25}}\big(\cup \big\{ f^2\big(K_{\al_0\al_1\al_2}^{(j,25)}\big): \, \al_i = 0, 1, \, 0 \le i \le 2, \, j = 1, 2 \big\}\big) \subset \widehat V_{n_2}(y_1) \,\,\, \text{and} \,\,\, f^{k_{25}}(x_1) \in V_{n_2}(y_1).\vspace{.1in}
$$
\indent By arguing as those in the previous six steps with $x_1$ replaced by $x_2$ and $y_1$ by $y_2$, we obtain successively
\begin{multline*}
$$
\text{an integer} \,\,\, k_{26} \, (> k_{25}) \,\,\, \text{and} \,\,\, 2 \cdot 2^3 \,\,\, \text{pairwise disjoint compact intervals} \\      K_{\al_0\al_1\al_2}^{(j,26)} \subset K_{\al_0\al_1\al_2}^{(j,25)}, \, \al_i = 0, 1, \, 0 \le i \le 2, \, j = 1, 2, \, \text{and} \\
\text{an integer} \,\,\, k_{27} \, (> k_{26}) \,\,\, \text{and} \,\,\, 2 \cdot 2^3 \,\,\, \text{pairwise disjoint compact intervals}\qquad\qquad\qquad\qquad\qquad\quad\,\,\, \\ 
K_{\al_0\al_1\al_2}^{(j,27)} \subset K_{\al_0\al_1\al_2}^{(j,26)}, \, \al_i = 0, 1, \, 0 \le i \le 2, \, j = 1, 2, \, \text{and}  \\
\text{an integer} \,\,\, k_{28} \, (> k_{27}) \,\,\, \text{and} \,\,\, 2 \cdot 2^3 \,\,\, \text{pairwise disjoint compact intervals}\qquad\qquad\qquad\qquad\qquad\quad\,\,\, \\      
K_{\al_0\al_1\al_2}^{(j,28)} \subset K_{\al_0\al_1\al_2}^{(j,27)}, \, \al_i = 0, 1, \, 0 \le i \le 2, \, j = 1, 2, \, \text{and} \\
\text{an integer} \,\,\, k_{29} \, (> k_{28}) \,\,\, \text{and} \,\,\, 2 \cdot 2^3 \,\,\, \text{pairwise disjoint compact intervals}\qquad\qquad\qquad\qquad\qquad\quad\,\,\, \\ 
K_{\al_0\al_1\al_2}^{(j,29)} \subset K_{\al_0\al_1\al_2}^{(j,28)}, \, \al_i = 0, 1, \, 0 \le i \le 2, \, j = 1, 2, \, \text{and}  \\
\text{an integer} \,\,\, k_{30} \, (> k_{29}) \,\,\, \text{and} \,\,\, 2 \cdot 2^3 \,\,\, \text{pairwise disjoint compact intervals}\qquad\qquad\qquad\qquad\qquad\quad\,\,\, \\ 
K_{\al_0\al_1\al_2}^{(j,30)} \subset K_{\al_0\al_1\al_2}^{(j,29)}, \, \al_i = 0, 1, \, 0 \le i \le 2, \, j = 1, 2, \, \text{and}  \\
\text{an integer} \,\,\, k_{31} \, (> k_{30}) \,\,\, \text{and} \,\,\, 2 \cdot 2^3 \,\,\, \text{pairwise disjoint compact intervals}\qquad\qquad\qquad\qquad\qquad\quad\,\,\, \\ 
K_{\al_0\al_1\al_2}^{(j,31)} \subset K_{\al_0\al_1\al_2}^{(j,30)}, \, \al_i = 0, 1, \, 0 \le i \le 2, \, j = 1, 2, \qquad\qquad\qquad\qquad\qquad\quad
$$
\end{multline*}
such that
$$
f^{k_{26}}\big(\cup \big\{ K_{\al_0\al_1\al_2}^{(j,26)}: \, \al_i = 0, 1, \, 0 \le i \le 2, \, j = 1, 2 \big\}\big) \, \subset V_{n_2}(y_2) \,\,\, \text{and} \,\,\, f^{k_{26}}(x_2) \in V_{n_2}(y_2),  
$$
$$
f^{k_{27}}\big(\cup \big\{ f\big(K_{\al_0\al_1\al_2}^{(j,27)}\big): \, \al_i = 0, 1, \, 0 \le i \le 2, \, j = 1, 2 \big\}\big) \subset V_{n_2}(y_2) \,\,\, \text{and} \,\,\, f^{k_{27}}(x_2) \in V_{n_2}(y_2), 
$$
$$
f^{k_{28}}\big(\cup \big\{ f^2\big(K_{\al_0\al_1\al_2}^{(j,28)}\big): \, \al_i = 0, 1, \, 0 \le i \le 2, \, j = 1, 2 \big\}\big) \subset V_{n_2}(y_2) \,\,\, \text{and} \,\,\, f^{k_{28}}(x_2) \in V_{n_2}(y_2),
$$
$$
f^{k_{29}}\big(\cup \big\{ K_{\al_0\al_1\al_2}^{(j,29)}: \, \al_i = 0, 1, \, 0 \le i \le 2, \, j = 1, 2 \big\}\big) \, \subset \widehat V_{n_2}(y_2) \,\,\, \text{and} \,\,\, f^{k_{29}}(x_2) \in V_{n_2}(y_2),  
$$
$$
f^{k_{30}}\big(\cup \big\{ f\big(K_{\al_0\al_1\al_2}^{(j,30)}\big): \, \al_i = 0, 1, \, 0 \le i \le 2, \, j = 1, 2 \big\}\big) \subset \widehat V_{n_2}(y_2) \,\,\, \text{and} \,\,\, f^{k_{30}}(x_2) \in V_{n_2}(y_2), \, \text{and}
$$
$$
f^{k_{31}}\big(\cup \big\{ f^2\big(K_{\al_0\al_1\al_2}^{(j,31)}\big): \, \al_i = 0, 1, \, 0 \le i \le 2, \, j = 1, 2 \big\}\big) \subset \widehat V_{n_2}(y_2) \,\,\, \text{and} \,\,\, f^{k_{31}}(x_2) \in V_{n_2}(y_2).
$$
\indent On the other hand, by applying Theorem 6(2) to the following $2(2 \cdot 2^3)$ intervals
$$
K_{\al_0\al_1\al_2}^{(j,31)}, \, \al_i = 0, 1, \, 0 \le i \le 2, \, j = 1, 2 \,\, : \,\, \big(V_{n_2}(y_1)\big)^{16},
$$
we obtain a positive integer $\widehat N$ such that 
$$
f^n\big(K_{\al_0\al_1\al_2}^{(j,31)}\big) \cap V_{n_2}(y_1) \ne \emptyset \,\,\, \text{for all} \,\,\, \al_i = 0, 1, \, 0 \le i \le 2, \, j = 1, 2, \, \text{and all} \,\,\, n \ge \widehat N.
$$
\indent Now, since $y_1$ is an $\omega$-limit point of $x_1$, the open neighborhood $V_{n_2}(y_1)$ of $y_1$ contains $f^i(x_1)$ for infinitely many $i$'s.  Let $k_{32}+1$ be one such $i$ which is $> \widehat N+k_{31}+1$.  Then $f^{k_{32}}\big(f(x_1)\big) = f^{k_{32}+1}(x_1) \in V_{n_2}(y_1)$ and $k_{32} > \widehat N$.  Let $K_{\al_0\al_1\al_2}^{(j,32)} \subset K_{\al_0\al_1\al_2}^{(j,31)}$, $\al_i = 0, 1$, $0 \le i \le 2$, $j = 1, 2$, be $2 \cdot 2^3$ pairwise disjoint compact intervals such that 
\begin{multline*}
$$
f^{k_{32}}\big(\cup \big\{ K_{\al_0\al_1\al_2}^{(j,32)}: \, \al_i = 0, 1, \, 0 \le i \le 2, \, j = 1, 2 \big\}\big) \subset V_{n_2}(y_1) \\ \big(\text{and} \,\,\, f^{k_{32}}\big(f(x_1)\big) = f^{k_{32}+1}(x_1) \in V_{n_2}(y_1)\big).
$$
\end{multline*}
\indent By applying Theorem 6(2) to the following $2(2 \cdot 2^3)$ intervals 
$$
K_{\al_0\al_1\al_2}^{(j,32)}, \, \al_i = 0, 1, \, 0 \le i \le 2, \, j = 1, 2 \,\, : \,\, \big(V_{n_2}(y_1)\big)^{16},
$$
we obtain a positive integer $\widehat N^*$ such that 
$$
f^n\big(K_{\al_0\al_1\al_2}^{(j,32)}\big) \cap V_{n_2}(y_1) \ne \emptyset \,\,\, \text{for all} \,\,\, \al_i = 0, 1, \, 0 \le i \le 2, \, j = 1, 2, \, \text{and all} \,\,\, n \ge \widehat N^*.
$$
\indent Now, since $y_1$ is an $\omega$-limit point of $x_1$, the open neighborhood $V_{n_2}(y_1)$ of $y_1$ contains $f^i(x_1)$ for infinitely many $i$'s.  Let $k_{33}+2$ be one such $i$ which is $> \widehat N^*+k_{32}+2$.  Then $f^{k_{33}}\big(f^2(x_1)\big) = f^{k_{33}+2}(x_1) \in V_{n_2}(y_1)$ and $k_{33} > \widehat N^*$.  Let $K_{\al_0\al_1\al_2}^{(j,33)} \subset K_{\al_0\al_1\al_2}^{(j,32)}$, $\al_i = 0, 1$, $0 \le i \le 2$, $j = 1, 2$, be $2 \cdot 2^3$ pairwise disjoint compact intervals such that 
\begin{multline*}
$$
f^{k_{33}}\big(\cup \big\{ K_{\al_0\al_1\al_2}^{(j,33)}: \, \al_i = 0, 1, \, 0 \le i \le 2, \, j = 1, 2 \big\}\big) \subset V_{n_2}(y_1) \\ \big(\text{and} \,\,\, f^{k_{33}}\big(f^2(x_1)\big) = f^{k_{33}+2}(x_1) \in V_{n_2}(y_1)\big).
$$
\end{multline*}

\indent By arguments similar to the previous two steps with {\it apptopriate} $V_{n_2}(y_1)$ replaced by $\widehat V_{n_2}(y_1)$, we obain successively
\begin{multline*}
$$
\text{an integer} \,\,\, k_{34} \, (> k_{33}) \,\,\, \text{and} \,\,\, 2 \cdot 2^3 \,\,\, \text{pairwise disjoint compact intervals}\qquad\qquad\qquad\qquad\qquad\quad\,\,\, \\ 
K_{\al_0\al_1\al_2}^{(j,34)} \subset K_{\al_0\al_1\al_2}^{(j,33)}, \, \al_i = 0, 1, \, 0 \le i \le 2, \, j = 1, 2, \, \text{and}  \\
\text{an integer} \,\,\, k_{35} \, (> k_{34}) \,\,\, \text{and} \,\,\, 2 \cdot 2^3 \,\,\, \text{pairwise disjoint compact intervals}\qquad\qquad\qquad\qquad\qquad\quad\,\,\, \\ 
K_{\al_0\al_1\al_2}^{(j,35)} \subset K_{\al_0\al_1\al_2}^{(j,34)}, \, \al_i = 0, 1, \, 0 \le i \le 2, \, j = 1, 2, \qquad\qquad\qquad\qquad\qquad\quad
$$
\end{multline*}
such that
\begin{multline*}
$$
f^{k_{34}}\big(\cup \big\{ K_{\al_0\al_1\al_2}^{(j,34)}: \, \al_i = 0, 1, \, 0 \le i \le 2, \, j = 1, 2 \big\}\big) \subset \widehat  V_{n_2}(y_1) \\ \text{and} \,\,\, f^{k_{34}}\big(f(x_1)\big) = f^{k_{34}+1}(x_1) \in V_{n_2}(y_1), \, \text{and}
$$
\end{multline*}
\begin{multline*}
$$
f^{k_{35}}\big(\cup \big\{ K_{\al_0\al_1\al_2}^{(j,35)}: \, \al_i = 0, 1, \, 0 \le i \le 2, \, j = 1, 2 \big\}\big) \subset \widehat  V_{n_2}(y_1) \\ \text{and} \,\,\, f^{k_{35}}\big(f^2(x_1)\big) = f^{k_{35}+2}(x_1) \in V_{n_2}(y_1).
$$
\end{multline*}
\indent By arguments similar to the previous four steps with $x_1$ replaced by $x_2$ and $y_1$ by $y_2$, we obtain successively
\begin{multline*}
$$
\text{an integer} \,\,\, k_{36} \, (> k_{35}) \,\,\, \text{and} \,\,\, 2 \cdot 2^3 \,\,\, \text{pairwise disjoint compact intervals} \\      K_{\al_0\al_1\al_2}^{(j,36)} \subset K_{\al_0\al_1\al_2}^{(j,35)}, \, \al_i = 0, 1, \, 0 \le i \le 2, \, j = 1, 2, \, \text{and} \\
\text{an integer} \,\,\, k_{37} \, (> k_{36}) \,\,\, \text{and} \,\,\, 2 \cdot 2^3 \,\,\, \text{pairwise disjoint compact intervals}\qquad\qquad\qquad\qquad\qquad\quad\,\,\, \\ 
K_{\al_0\al_1\al_2}^{(j,37)} \subset K_{\al_0\al_1\al_2}^{(j,36)}, \, \al_i = 0, 1, \, 0 \le i \le 2, \, j = 1, 2, \, \text{and}  \\
\text{an integer} \,\,\, k_{38} \, (> k_{37}) \,\,\, \text{and} \,\,\, 2 \cdot 2^3 \,\,\, \text{pairwise disjoint compact intervals}\qquad\qquad\qquad\qquad\qquad\quad\,\,\, \\      
K_{\al_0\al_1\al_2}^{(j,38)} \subset K_{\al_0\al_1\al_2}^{(j,37)}, \, \al_i = 0, 1, \, 0 \le i \le 2, \, j = 1, 2, \, \text{and} \\
\text{an integer} \,\,\, k_{39} \, (> k_{38}) \,\,\, \text{and} \,\,\, 2 \cdot 2^3 \,\,\, \text{pairwise disjoint compact intervals}\qquad\qquad\qquad\qquad\qquad\quad\,\,\, \\ 
K_{\al_0\al_1\al_2}^{(j,39)} \subset K_{\al_0\al_1\al_2}^{(j,38)}, \, \al_i = 0, 1, \, 0 \le i \le 2, \, j = 1, 2,  \qquad\qquad\qquad\qquad\qquad\quad
$$
\end{multline*}
such that
\begin{multline*}
$$
f^{k_{36}}\big(\cup \big\{ K_{\al_0\al_1\al_2}^{(j,36)}: \, \al_i = 0, 1, \, 0 \le i \le 2, \, j = 1, 2 \big\}\big) \subset V_{n_2}(y_2) \\ \text{and} \,\,\, f^{k_{36}}\big(f(x_2)\big) = f^{k_{36}+1}(x_2) \in V_{n_2}(y_2), \, \text{and}
$$
\end{multline*}
\begin{multline*}
$$
f^{k_{37}}\big(\cup \big\{ K_{\al_0\al_1\al_2}^{(j,37)}: \, \al_i = 0, 1, \, 0 \le i \le 2, \, j = 1, 2 \big\}\big) \subset V_{n_2}(y_2) \\ \text{and} \,\,\, f^{k_{37}}\big(f^2(x_2)\big) = f^{k_{37}+2}(x_2) \in V_{n_2}(y_2).
$$
\end{multline*}
\begin{multline*}
$$
f^{k_{38}}\big(\cup \big\{ K_{\al_0\al_1\al_2}^{(j,38)}: \, \al_i = 0, 1, \, 0 \le i \le 2, \, j = 1, 2 \big\}\big) \subset \widehat  V_{n_2}(y_2) \\ \text{and} \,\,\, f^{k_{38}}\big(f(x_2)\big) = f^{k_{38}+1}(x_2) \in V_{n_2}(y_2), \, \text{and}
$$
\end{multline*}
\begin{multline*}
$$
f^{k_{39}}\big(\cup \big\{ K_{\al_0\al_1\al_2}^{(j,39)}: \, \al_i = 0, 1, \, 0 \le i \le 2, \, j = 1, 2 \big\}\big) \subset \widehat  V_{n_2}(y_2) \\ \text{and} \,\,\, f^{k_{39}}\big(f^2(x_2)\big) = f^{k_{39}+2}(x_2) \in V_{n_2}(y_2).
$$
\end{multline*}

Since, for each $\al_i = 0, 1$, $0 \le i \le 2$ and each $j = 1, 2$, we have
$K_{\al_0\al_1}^{(j,10)} \supset K_{\al_0\al_1\al_2}^{(j,11)} \supset K_{\al_0\al_1\al_2}^{(j,12)} \supset K_{\al_0\al_1\al_2}^{(j,13)} \supset \cdots \supset K_{\al_0\al_1\al_2}^{(j,37)} \supset K_{\al_0\al_1\al_2}^{(j,38)} \supset K_{\al_0\al_1\al_2}^{(j,39)}$,
we can summarize the arguments in the (above) second stage as follows:

{\footnotesize \[
\left(\begin{array}{l}
\textrm{\qquad Let $n_2$ be any positive integer $> k_{10}$.}\\  
\textrm{Take these $2^2$ pairwise disjoint compact intervals $K_{\al_0\al_1}^{(1,10)}$, $\al_0=0, 1$, $\al_1=0, 1$ in $U_1$ obtained in the}\\ 
\textrm{first stage and take any $2^2$ pairwise disjoint compact intervals $K_{\al_0\al_1}^{(2,10)}$, $\al=0, 1$, $\al_1=0, 1$ in}\\
\textrm{$U_2 \setminus \big(K_{00}^{(1,10)} \cup K_{01}^{(1,10)} \cup K_{10}^{(1,10)} \cup K_{11}^{(1,10)}\big)$.  By applying Theorem 6(2) to the following $2 \cdot 2^4$ intervals}\\
\textrm{\scriptsize{$\big(K_{00}^{(1,10)}\big)^2, \, \big(K_{01}^{(1,10)}\big)^2, \, \big(K_{10}^{(1,10)}\big)^2, \, \big(K_{11}^{(1,10)}\big)^2; \, \big(K_{00}^{(2,10)}\big)^2, \, \big(K_{01}^{(2,10)}\big)^2, \, \big(K_{10}^{(2,10)}\big)^2, \, \big(K_{11}^{(2,10)}\big)^2 : \, \big(V_{n_2}(a_2), \, V_{n_2}(b_2)\big)^8$,}} \\ 
\textrm{we obtain a positive integer $k_{11}$ $(> n_2)$ in $\mathcal M$ and $2 \cdot 2^3$ pairwise disjoint compact intervals}\\
\textrm{$K_{\al_0\al_1\al_2}^{(1,11)} \subset K_{\al_0\al_1}^{(1,10)}$, $\al_i = 0,1$, $0 \le i \le 2$ \,\, and \,\, $K_{\al_0\al_1\al_2}^{(2,11)} \subset K_{\al_0\al_1}^{(2,10)}$, $\al_i = 0,1$, $0 \le i \le 2$ such that}\\
\textrm{\qquad each length is so smaller than $1/2^5$ that the open set}\\
\textrm{\qquad\qquad\qquad $U_3 \setminus \big( \bigcup \big\{ K_{\al_0\al_1\al_2}^{(j,11)}: \, \al_i = 0, 1, \, 0 \le i \le 2, \, j = 1, 2 \big\} \big)$ is nonempty and,}\\
\textrm{for all $\al_0 = 0, 1$, $\al_1 = 0, 1$, $\al_2 = 0,1$, and $j = 1, 2$, we have} \\
\textrm{\qquad\qquad\quad $f^{k_{11}}\big(K_{\al_0\al_1\al_2}^{(j,11)}\big) \subset W_{n_2}(\al_2)$, where $W_{n_2}(0) = V_{n_2}(a_2)$ and $W_{n_2}(1) = V_{n_2}(b_2)$.}\\
\textrm{Start with $K_{\al_0\al_1\al_2}^{(j,11)}$, $\al_i=0, 1$, $0 \le i \le 2$, $j = 1, 2$, by applying Theorem 6(2) \scriptsize{successively and appropriately,}}\\
\textrm{we obtain $(2 \cdot 2+1)^2 +2^2= 29$ positive integers $(k_{10} < n_2 <)$ $k_{11} < k_{12} < \cdots < k_{38} < k_{39}$ and}\\
\textrm{\quad\quad $2 \cdot 2^3$ pairwise disjoint compact intervals $K_{\al_0\al_1\al_2}^{(j,39)}$, $\al_i= 0, 1$, $0 \le i \le 2$, $j =1,2$, such that}\\
\textrm{$\{ k_{11}, k_{12}, k_{13}, \cdots, k_{18}  \} \subset \mathcal M$, $2!$ divides $k_{19}$ and} \\
\textrm{$f^{k_{11}}\big(K_{\al_0\al_1\al_2}^{(1,39)} \cup K_{\al_0\al_1\al_2}^{(2,39)}\big) \subset W_{n_2}(\al_2)$, where $W_{n_2}(0) = V_{n_2}(a_2)$ and $W_{n_2}(1) = V_{n_2}(b_2)$,}\\
\textrm{$f^{k_{12}}\big(f(K_{\al_0\al_1\al_2}^{(1,39)}) \cup f(K_{\al_0\al_1\al_2}^{(2,39)})\big) \subset W_{n_2}(\al_2)$, where $W_{n_2}(0) = V_{n_2}(a_2)$ and $W_{n_2}(1) = V_{n_2}(b_2)$,}\\
\textrm{$f^{k_{13}}\big(f^2(K_{\al_0\al_1\al_2}^{(1,39)}) \cup f^2(K_{\al_0\al_1\al_2}^{(2,39)})\big) \subset W_{n_2}(\al_2)$, where $W_{n_2}(0) = V_{n_2}(a_2)$ and $W_{n_2}(1) = V_{n_2}(b_2)$,}\\
\textrm{$f^{k_{14}}\big(\bigcup_{\al_i = 0, 1, \, 0 \le i \le 2} \, K_{\al_0\al_1\al_2}^{(1,39)}\big) \subset V_{n_2}(a_2) \,\,\, \text{and} \,\,\, f^{k_{14}}\big(\bigcup_{\al_i = 0, 1, \, 0 \le i \le 2} K_{\al_0\al_1\al_2}^{(2,39)}\big) \subset V_{n_2}(b_2)$,}\\
\textrm{$f^{k_{15}}\big(\bigcup_{\al_i = 0, 1, \, 0 \le i \le 2} \, f(K_{\al_0\al_1\al_2}^{(1,39)})\big) \subset V_{n_2}(a_2) \,\,\, \text{and} \,\,\, f^{k_{15}}\big(\bigcup_{\al_i = 0, 1, \, 0 \le i \le 2} f(K_{\al_0\al_1\al_2}^{(2,39)})\big) \subset V_{n_2}(b_2)$,}\\
\textrm{$f^{k_{16}}\big(\bigcup_{\al_i = 0, 1, \, 0 \le i \le 2} \, f^2(K_{\al_0\al_1\al_2}^{(1,39)})\big) \subset V_{n_2}(a_2) \,\,\, \text{and} \,\,\, f^{k_{16}}\big(\bigcup_{\al_i = 0, 1, \, 0 \le i \le 2} f^2(K_{\al_0\al_1\al_2}^{(2,39)})\big) \subset V_{n_2}(b_2)$,}\\
\textrm{$f^{k_{17}}\big(\bigcup \big\{ K_{\al_0\al_1\al_2}^{(j,39)}: \, \al_i = 0, 1, \, 0 \le i \le 2, \, j = 1, 2 \big\}\big) \, \subset V_{n_2}(v_1)$, and}\\
\textrm{$f^{k_{18}}\big(\bigcup \big\{ K_{\al_0\al_1\al_2}^{(j,39)}: \, \al_i = 0, 1, \, 0 \le i \le 2, \, j = 1, 2 \big\}\big) \, \subset V_{n_2}(v_2)$,}\\
\textrm{$f^{k_{19}}\big(\bigcup \big\{ K_{\al_0\al_1\al_2}^{(j,39)}: \, \al_i = 0, 1, \, 0 \le i \le 2, \, j = 1, 2 \big\}\big) \, \subset U_2$,}\\
\textrm{$f^{k_{20}}\big(\bigcup \big\{ K_{\al_0\al_1\al_2}^{(j,39)}: \, \al_i = 0, 1, \, 0 \le i \le 2, \, j = 1, 2 \big\}\big) \, \subset V_{n_2}(y_1)$ \, \text{and} \, $f^{20}(x_1) \in V_{n_2}(y_1)$,}\\
\textrm{$f^{k_{21}}\big(\bigcup \big\{ f\big(K_{\al_0\al_1\al_2}^{(j,39)}\big): \, \al_i = 0, 1, \, 0 \le i \le 2, \, j = 1, 2 \big\}\big) \, \subset V_{n_2}(y_1)$ \, \text{and} \, $f^{21}(x_1) \in V_{n_2}(y_1)$,}\\
\textrm{$f^{k_{22}}\big(\bigcup \big\{ f^2\big(K_{\al_0\al_1\al_2}^{(j,39)}\big): \, \al_i = 0, 1, \, 0 \le i \le 2, \, j = 1, 2 \big\}\big) \, \subset V_{n_2}(y_1)$ \, \text{and} \, $f^{22}(x_1) \in V_{n_2}(y_1)$,}\\
\textrm{$f^{k_{23}}\big(\bigcup \big\{ K_{\al_0\al_1\al_2}^{(j,39)}: \, \al_i = 0, 1, \, 0 \le i \le 2, \, j = 1, 2 \big\}\big) \, \subset \widehat V_{n_2}(y_1)$ \, \text{and} \, $f^{23}(x_1) \in V_{n_2}(y_1)$,}\\
\textrm{$f^{k_{24}}\big(\bigcup \big\{ f\big(K_{\al_0\al_1\al_2}^{(j,39)}\big): \, \al_i = 0, 1, \, 0 \le i \le 2, \, j = 1, 2 \big\}\big) \, \subset \widehat V_{n_2}(y_1)$ \, \text{and} \, $f^{24}(x_1) \in V_{n_2}(y_1)$,}\\
\textrm{$f^{k_{25}}\big(\bigcup \big\{ f^2\big(K_{\al_0\al_1\al_2}^{(j,39)}\big): \, \al_i = 0, 1, \, 0 \le i \le 2, \, j = 1, 2 \big\}\big) \, \subset \widehat V_{n_2}(y_1)$ \, \text{and} \, $f^{25}(x_1) \in V_{n_2}(y_1)$,}\\
\textrm{$f^{k_{26}}\big(\bigcup \big\{ K_{\al_0\al_1\al_2}^{(j,39)}: \, \al_i = 0, 1, \, 0 \le i \le 2, \, j = 1, 2 \big\}\big) \, \subset V_{n_2}(y_2)$ \, \text{and} \, $f^{26}(x_2) \in V_{n_2}(y_2)$,}\\
\textrm{$f^{k_{27}}\big(\bigcup \big\{ f\big(K_{\al_0\al_1\al_2}^{(j,39)}\big): \, \al_i = 0, 1, \, 0 \le i \le 2, \, j = 1, 2 \big\}\big) \, \subset V_{n_2}(y_2)$ \, \text{and} \, $f^{27}(x_2) \in V_{n_2}(y_2)$,}\\
\textrm{$f^{k_{28}}\big(\bigcup \big\{ f^2\big(K_{\al_0\al_1\al_2}^{(j,39)}\big): \, \al_i = 0, 1, \, 0 \le i \le 2, \, j = 1, 2 \big\}\big) \, \subset V_{n_2}(y_2)$ \, \text{and} \, $f^{28}(x_2) \in V_{n_2}(y_2)$,}\\
\textrm{$f^{k_{29}}\big(\bigcup \big\{ K_{\al_0\al_1\al_2}^{(j,39)}: \, \al_i = 0, 1, \, 0 \le i \le 2, \, j = 1, 2 \big\}\big) \, \subset \widehat V_{n_2}(y_2)$ \, \text{and} \, $f^{29}(x_2) \in V_{n_2}(y_2)$,}\\
\textrm{$f^{k_{30}}\big(\bigcup \big\{ f\big(K_{\al_0\al_1\al_2}^{(j,39)}\big): \, \al_i = 0, 1, \, 0 \le i \le 2, \, j = 1, 2 \big\}\big) \, \subset \widehat V_{n_2}(y_2)$ \, \text{and} \, $f^{30}(x_2) \in V_{n_2}(y_2)$,}\\
\textrm{$f^{k_{31}}\big(\bigcup \big\{ f^2\big(K_{\al_0\al_1\al_2}^{(j,39)}\big): \, \al_i = 0, 1, \, 0 \le i \le 2, \, j = 1, 2 \big\}\big) \, \subset \widehat V_{n_2}(y_2)$ \, \text{and} \, $f^{31}(x_2) \in V_{n_2}(y_2)$.}\\
\textrm{$f^{k_{32}}\big(\bigcup \big\{ K_{\al_0\al_1\al_2}^{(j,39)}: \, \al_i = 0, 1, \, 0 \le i \le 2, \, j = 1, 2 \big\}\big) \, \subset V_{n_2}(y_1)$ \, \text{and} \, $f^{32}\big(f(x_1)\big) \in V_{n_2}(y_1)$,}\\
\textrm{$f^{k_{33}}\big(\bigcup \big\{ K_{\al_0\al_1\al_2}^{(j,39)}: \, \al_i = 0, 1, \, 0 \le i \le 2, \, j = 1, 2 \big\}\big) \, \subset V_{n_2}(y_1)$ \, \text{and} \, $f^{33}\big(f^2(x_1)\big) \in V_{n_2}(y_1)$,}\\
\textrm{$f^{k_{34}}\big(\bigcup \big\{ K_{\al_0\al_1\al_2}^{(j,39)}: \, \al_i = 0, 1, \, 0 \le i \le 2, \, j = 1, 2 \big\}\big) \, \subset \widehat V_{n_2}(y_1)$ \, \text{and} \, $f^{34}\big(f(x_1)\big) \in V_{n_2}(y_1)$,}\\
\textrm{$f^{k_{35}}\big(\bigcup \big\{ K_{\al_0\al_1\al_2}^{(j,39)}: \, \al_i = 0, 1, \, 0 \le i \le 2, \, j = 1, 2 \big\}\big) \, \subset \widehat V_{n_2}(y_1)$ \, \text{and} \, $f^{35}\big(f^2(x_1)\big) \in V_{n_2}(y_1)$,}\\
\textrm{$f^{k_{36}}\big(\bigcup \big\{ K_{\al_0\al_1\al_2}^{(j,39)}: \, \al_i = 0, 1, \, 0 \le i \le 2, \, j = 1, 2 \big\}\big) \, \subset V_{n_2}(y_2)$ \, \text{and} \, $f^{36}\big(f(x_2)\big) \in V_{n_2}(y_2)$,}\\
\textrm{$f^{k_{37}}\big(\bigcup \big\{ K_{\al_0\al_1\al_2}^{(j,39)}: \, \al_i = 0, 1, \, 0 \le i \le 2, \, j = 1, 2 \big\}\big) \, \subset V_{n_2}(y_2)$ \, \text{and} \, $f^{37}\big(f^2(x_2)\big) \in V_{n_2}(y_2)$,}\\
\textrm{$f^{k_{38}}\big(\bigcup \big\{ K_{\al_0\al_1\al_2}^{(j,39)}: \, \al_i = 0, 1, \, 0 \le i \le 2, \, j = 1, 2 \big\}\big) \, \subset \widehat V_{n_2}(y_2)$ \, \text{and} \, $f^{38}\big(f(x_2)\big) \in V_{n_2}(y_2)$,}\\
\textrm{$f^{k_{39}}\big(\bigcup \big\{ K_{\al_0\al_1\al_2}^{(j,39)}: \, \al_i = 0, 1, \, 0 \le i \le 2, \, j = 1, 2 \big\}\big) \, \subset \widehat V_{n_2}(y_2)$ \, \text{and} \, $f^{39}\big(f^2(x_2)\big) \in V_{n_2}(y_2)$.}\\
\end{array}\right)
\]}

\indent By proceeding in this manner indefinitely, at the $\ell^{th}$ stage for each integer $\ell \ge 1$, we obtain $(2\ell+1)^2+\ell^2$ positive integers (let $k_0 = 0$)
$$
\big(k_{(\ell-1)(10\ell^2+7\ell+6)/6} < n_\ell <\big) \, k_{[(\ell-1)(10\ell^2+7\ell+6)/6]+1} < k_{[(\ell-1)(10\ell^2+7\ell+6)/6]+2} < \cdots < k_{\ell(10\ell^2+27\ell+23)/6}
$$ 
with 
$$
\big\{ k_{[(\ell-1)(10\ell^2+7\ell+6)/6]+i}: 1 \le i \le \ell(\ell+2) \big\} \subset \mathcal M \,\,\, \text{and} \,\,\, {\ell}! \,\,\, \text{dividing} \,\,\, k_{[(\ell-1)(10\ell^2+7\ell+6)/6]+(\ell+1)^2} 
$$ 
and $\ell \cdot 2^{\ell+1}$ pairwise disjoint compact intervals 
$$
K_{\al_0\al_1\al_2\cdots \al_\ell}^{(j, \,\ell(10\ell^2+27\ell+23)/6)} \,\, \big(\subset K_{\al_0\al_1\al_2\cdots \al_{\ell-1}}^{(j,\, (\ell-1)(10\ell^2+7\ell+6)/6)}\big), \,\, \al_i = 0, 1, \,\, 0 \le i \le \ell, \,\, 1 \le j \le \ell,
$$ 
such that each length is so smaller than $\frac 1{2^{2\ell+1}}$ that the open set 
$$
U_{\ell+1} \setminus \ \bigg( \bigcup \big\{ K_{\al_0\al_1\al_2\cdots \al_\ell}^{(j,\,\ell(10\ell^2+27\ell+23)/6)}: \, \al_i = 0, 1, \, 0 \le i \le \ell, \, 1 \le j \le \ell \big\} \bigg) \,\,\, \text{is nonempty} \,\,\, \text{and}
$$
the following hold:

\noindent
for each $0 \le r \le \ell-1$, 
\begin{multline*}
$$
\qquad\qquad \text{if} \,\,\, r = 0, \, \text{then}, \, \text{for all} \,\,\, 0 \le s \le \ell, \\
f^{k_{[(\ell-1)(10\ell^2+7\ell+6)/6)]+r(\ell+1)+1+s}}\bigg(\bigcup_{j=1}^\ell \big\{ f^s\big(K_{\al_0\al_1\al_2\cdots\al_\ell}^{(j,\,\ell(10\ell^2+27\ell+23)/6)}\big): \, \al_i = 0, 1, \, 0 \le i \le \ell \big\}\bigg) \subset W_{n_\ell}(\al_\ell), \\
\text{where} \,\,\, W_{n_\ell}(0) = V_{n_\ell}(a_\ell) \,\,\, \text{and} \,\,\, W_{n_\ell}(1) = V_{n_\ell}(b_\ell) \,\,\, \text{and},
$$
\end{multline*}
\begin{multline*}
$$
\qquad\qquad \text{if} \,\,\, 1 \le r \le \ell-1, \, \text{then} \,\,\, \big(\text{void if} \,\,\, \ell = 1\big), \, \text{for all} \,\,\, 0 \le s \le \ell, \\
f^{k_{[(\ell-1)(10\ell^2+7\ell+6)/6]+r(\ell+1)+1+s}}\bigg(\bigcup_{j_1=1}^r \big\{ f^s\big(K_{\al_0\al_1\al_2\cdots\al_\ell}^{(j_1,\,\ell(10\ell^2+27\ell+23)/6)}\big): \, \al_i = 0, 1, \, 0 \le i \le \ell \big\}\bigg) \subset V_{n_\ell}(a_\ell) \,\,\, \text{and} \qquad \\ 
f^{k_{[(\ell-1)(10\ell^2+7\ell+6)/6]+r(\ell+1)+1+s}}\bigg(\bigcup_{j_2=r+1}^\ell \big\{ f^s\big(K_{\al_0\al_1\al_2\cdots\al_\ell}^{(j_2,\,\ell(10\ell^2+27\ell+23)/6)}\big): \, \al_i = 0, 1, \, 0 \le i \le \ell \big\}\bigg) \subset V_{n_\ell}(b_\ell) \,\,\, \text{and},
$$
\end{multline*}
for each $1 \le m \le \ell$, 
$$
\qquad f^{k_{[(\ell-1)(10\ell^2+7\ell+6)/6]+\ell(\ell+1)+m}}\bigg(\bigcup_{j=1}^\ell \big\{ K_{\al_0\al_1\al_2\cdots\al_\ell}^{(j, \,\ell(10\ell^2+27\ell+23)/6)}: \, \al_i = 0, 1, \, 0 \le i \le \ell \big\}\bigg) 
\subset V_{n_\ell}(v_m),\, \text{and},
$$
$$
f^{k_{[(\ell-1)(10\ell^2+7\ell+6)/6]+(\ell+1)^2}}\bigg(\bigcup_{j=1}^\ell \big\{ K_{\al_0\al_1\al_2\cdots\al_\ell}^{(j, \,\ell(10\ell^2+27\ell+23)/6)}: \, \al_i = 0, 1, \, 0 \le i \le \ell \big\}\bigg) \, \subset U_\ell \,\,\, \text{and}, \qquad\qquad\quad\,\,
$$
for each $1 \le m \le \ell$,  
\begin{multline*}
$$
f^{k_{[(\ell-1)(10\ell^2+7\ell+6)/6]+(\ell+1)^2+(2m-2)(\ell+1)+1+s}}\bigg(\bigcup_{j=1}^\ell \big\{ f^s\big(K_{\al_0\al_1\al_2\cdots\al_\ell}^{(j, \, \ell(10\ell^2+27\ell+23)/6)}\big): \, \al_i = 0, 1, \, 0 \le i \le \ell \big\}\bigg) \\ 
\subset V_{n_\ell}(y_m) \,\,\, \text{and} \,\,\, f^{k_{[(\ell-1)(10\ell^2+7\ell+6)/6]+(\ell+1)^2+(2m-2)(\ell+1)+1+s}}(x_m) \in V_{n_\ell}(y_m) \,\,\, \text{for each} \,\,\, 0 \le s \le \ell, \, \text{and} \\
\,\,\,\,\, f^{k_{[(\ell-1)(10\ell^2+7\ell+6)/6]+(\ell+1)^2+(2m-1)(\ell+1)+1+s}}\bigg(\bigcup_{j=1}^\ell \big\{f^s\big(K_{\al_0\al_1\al_2\cdots\al_\ell}^{(j, \, \ell(10\ell^2+27\ell+23)/6)}\big): \, \al_i = 0, 1, \, 0 \le i \le \ell \big\}\bigg) \qquad\quad\,\, \\ 
\quad\,\,\, \subset \widehat V_{n_\ell}(y_m) \,\,\, \text{and} \,\,\, f^{k_{[(\ell-1)(10\ell^2+7\ell+6)/6]+(\ell+1)^2+(2m-1)(\ell+1)+1+s}}(x_m) \in V_{n_\ell}(y_m) \,\,\, \text{for each} \,\,\, 0 \le s \le \ell, \text{and} \qquad\,\,\,
$$
\end{multline*}
for each $1 \le m \le \ell$, 
\begin{multline*}
$$
f^{k_{[(\ell-1)(10\ell^2+7\ell+6)/6]+(\ell+1)^2+2\ell(\ell+1)+(2m-2)\ell+t}}\bigg(\bigcup_{j=1}^\ell \big\{ K_{\al_0\al_1\al_2\cdots\al_\ell}^{(j, \,\ell(10\ell^2+27\ell+23)/6)}: \, \al_i = 0, 1, \, 0 \le i \le \ell \big\}\bigg) \\ 
\,\, \subset V_{n_\ell}(y_m) \,\,\, \text{and} \,\,\, f^{k_{[(\ell-1)(10\ell^2+7\ell+6)/6]+(\ell+1)^2+2\ell(\ell+1)+(2m-2)\ell+t}}\big(f^{t}(x_m)\big) \in V_{n_\ell}(y_m) \,\,\, \text{for each} \,\,\, 1 \le t \le \ell, \, \text{and} \\
\,\,\,\,\, f^{k_{[(\ell-1)(10\ell^2+7\ell+6)/6]+(\ell+1)^2+2\ell(\ell+1)+(2m-1)\ell+t}}\bigg(\bigcup_{j=1}^\ell \big\{ K_{\al_0\al_1\al_2\cdots\al_\ell}^{(j, \,\ell(10\ell^2+27\ell+23)/6)}: \, \al_i = 0, 1, \, 0 \le i \le \ell \big\}\bigg) \qquad\qquad\,\,\,\,\, \\ 
\,\,\, \subset \widehat V_{n_\ell}(y_m) \,\,\, \text{and} \,\,\, f^{k_{[(\ell-1)(10\ell^2+7\ell+6)/6]+(\ell+1)^2+2\ell(\ell+1)+(2m-1)\ell+t}}\big(f^{t}(x_m)\big) \in V_{n_\ell}(y_m) \,\,\, \text{for each} \,\,\, 1 \le t \le \ell.\qquad\qquad\,
$$
\end{multline*}

Let
$$
\Sigma_2 = \{ \al : \al = \al_0\al_1\al_2 \cdots, \al_i = 0 \,\, \text{or} \,\, 1, \, i \ge 0 \}.
$$
\indent For any $\al = \al_0\al_1\al_2 \cdots \in \Sigma_2$, let 
{\large 
$$
\omega_\al = \breve \al_0 \breve \al_1 \breve \al_2 \cdots = \al_0 \,\, \al_0\al_1 \,\, \al_0\al_1\al_2 \,\, \al_0\al_1\al_2\al_3 \cdots \in \Sigma_2
$$} 
and, for each $j \ge 1$, let 
$$
K_\al^{(j)} = \bigcap_{\ell \, \ge \, j} \, K_{\al_0\al_1\al_2\cdots \al_\ell}^{(j, \, \ell(10\ell^2+27\ell+23)/6)} \, \big(\subset U_j\big).
$$ 
\big(The \, $j$ \, in the superscript of $K_\al^{(j)}$ indicates that the set $K_\al^{(j)}$ lies in $U_j$\big).  Then, since the sequence $< K_{\al_0\al_1\al_2\cdots \al_\ell}^{(j, \, \ell(10\ell^2+27\ell+23)/6)}>_{\ell \ge j}$ of compact intervals is nested and their lengths decrease to zero, each $K_\alpha^{(j)}$ consists of exactly one point, say $K_\al^{(j)} = \{ x_\al^{(j)} \}$.  Furthermore, since the total sum of the lengths of all these $\ell \cdot 2^{\ell+1}$ compact intervals $K_{\al_0\al_1\al_2\cdots \al_\ell}^{(j, \,\ell(10\ell^2+27\ell+23)/6)}$, $\al_i = 0, 1$, $0 \le i \le \ell$, $1 \le j \le \ell$, is $< \ell \cdot 2^{\ell+1} \cdot \frac 1{2^{2\ell+1}} = \frac \ell{2^\ell}$ and tends to zero as $\ell$ tends to $\infty$, the set $\{ x_\alpha^{(j)} : \alpha \in \Sigma_2 \}$ is a Cantor set with Lebesgue measure zero.

For each $j \ge 1$, let 
$$
\mathcal S^{(j)} = \{ x_{\omega_\al}^{(j)} : \al \in \Sigma_2 \} = \{ x_{\breve \al_0 \breve \al_1 \breve \al_2 \cdots}^{(j)}: \al \in \Sigma_2 \} \, \big(\subset U_j\big).
$$
\big(The \, $j$ \, in the superscript of $\mathcal S^{(j)}$ indicates that the set $\mathcal S^{(j)}$ lies in $U_j$\big).  Then it is easy to see that $\mathcal S^{(j)}$ is also a Cantor set with Lebesgue measure zero.  

For each $j \ge 1$ and all $\ell \ge j$, it follows from the definitions of \,\, $k_{[(\ell-1)(10\ell^2+7\ell+6)/6]+(\ell+1)^2}$ \,\, and the corresponding compact intervals (at the $\ell^{th}$ stage) $K_{\breve \al_0\breve \al_1\breve \al_2\cdots\breve \al_\ell}^{(j, \, (\ell-1)(10\ell^2+7\ell+6)/6)+(\ell+1)^2)} (\supset K_{\breve \al_0\breve \al_1\breve \al_2\cdots\breve \al_\ell}^{(j, \, \ell(10\ell^2+27\ell+23)/6)})$, $\al_n = 0, 1, \, 0 \le n \le \ell, \, 1 \le j \le \ell$, that
\begin{multline*}
$$
{\ell}! \,\,\, \text{divides} \,\,\, k_{[(\ell-1)(10\ell^2+7\ell+6)/6]+(\ell+1)^2} \,\,\, \text{and} \\ 
f^{k_{[(\ell-1)(10\ell^2+7\ell+6)/6]+(\ell+1)^2}}\big(\mathcal S^{(j)}\big) \subset f^{k_{[(\ell-1)(10\ell^2+7\ell+6)/6]+(\ell+1)^2}}\big(\mathcal S^{(1)} \cup \mathcal S^{(2)} \cup \mathcal S^{(3)} \cup \cdots \cup \mathcal S^{(\ell)} \big) \\
\subset f^{k_{[(\ell-1)(10\ell^2+7\ell+6)/6]+(\ell+1)^2}}\big(\cup_{j=1}^\ell \, \big\{ K_{\breve \al_0\breve \al_1\breve \al_2\cdots\breve \al_\ell}^{(j, \, (\ell-1)(10\ell^2+27\ell+23)/6+(\ell+1)^2)}: \, \al_n = 0, 1, \, 0 \le n \le \ell \big\}\big) \\
\subset f^{k_{[(\ell-1)(10\ell^2+7\ell+6)/6]+(\ell+1)^2}}\big(\cup_{j=1}^\ell \, \big\{ K_{\breve \al_0\breve \al_1\breve \al_2\cdots\breve \al_\ell}^{(j, \, (\ell-1)(10\ell^2+7\ell+6)/6+(\ell+1)^2)}: \, \al_n = 0, 1, \, 0 \le n \le \ell \big\}\big) \subset U_\ell.
$$  
\end{multline*}
\indent Since the collection $\{ U_1, U_2, U_3, \cdots \}$ is a countable open base for $I$ and since, for each $j \ge 1$, $f^{k_{[(\ell-1)(10\ell^2+7\ell+6)/6]+(\ell+1)^2}}\big(\mathcal S^{(j)}\big) \subset U_\ell$ for all $\ell \ge j$ and ${\ell}!$ divides $k_{[(\ell-1)(10\ell^2+7\ell+6)/6]+(\ell+1)^2}$, we obtain that, for each integer $n \ge 1$, $n$ divides $k_{[(\ell-1)(10\ell^2+7\ell+6)/6]+(\ell+1)^2}$ foe all $\ell \ge n$.  Therefore, for each $j \ge 1$ and each point $x$ in $\mathcal S^{(j)}$, $x$ is a transitive point of $f^n$ for all $n \ge 1$, i.e., each $x$ in $\mathcal S^{(j)}$ is a totally transitive point of $f$.  

Since, for each $j \ge 1$, we have $\mathcal S^{(j)} \subset U_j$, and, since $\{ U_1, U_2, U_3, \cdots \}$ is a countable open base for $I$, each nonempty open set in $I$ contains countably infinitely many $U_j$'s and so, contains countably infinitely many $\mathcal S^{(j)}$'s.  This proves Part (1).  

For each integer $r \ge 1$, since 
$$
\mathcal S^{(1)} \cup \mathcal S^{(2)} \cup \cdots \cup \mathcal S^{(r)} \subset \cup_{j=1}^\ell \,\, \big\{ K_{\breve \al_0\breve \al_1\breve \al_2\cdots\breve \al_\ell}^{(j, \, \ell(10\ell^2+27\ell+23)/6)}: \al_i = 0, 1, \, 0 \le i \le \ell \big\} \,\,\, \text{for each integer} \,\,\, \ell \ge r,
$$
and since, for each integer $m \ge 1$ and all integers $\ell \ge \max \{ r, m \}$,   
\begin{multline*}
$$
k_{[(\ell-1)(10\ell^2+7\ell+6)/6]+\ell(\ell+1)+m} \in \mathcal M \,\,\, \text{and} \\
f^{k_{[(\ell-1)(10\ell^2+7\ell+6)/6]+\ell(\ell+1)+m}}\big(\cup_{j=1}^\ell \,\, \big\{K_{\breve \al_0\breve \al_1\breve \al_2\cdots\breve \al_\ell}^{(j, \, \ell(10\ell^2+27\ell+23)/6)}:  \al_i = 0, 1, \, 0 \le i \le \ell \big\}\big) \subset \, V_{n_\ell}(v_m),
$$
\end{multline*}
and since the interval $V_{n_\ell}(v_m)$ shrinks to the point $v_m$ as $n_\ell$ tends to $\infty$, we obtain that, for each integer $r \ge 1$, 
$$
\text{the union} \,\,\, \mathcal S^{(1)} \cup \mathcal S^{(2)} \cup \cdots \cup \mathcal S^{(r)} \,\,\, \text{is} \,\,\, \mathcal M\text{-synchronouly proximal to} \,\,\, v_m \,\,\, \text{for each} \,\,\, m \ge 1. 
$$
\indent On the other hand, for each $r \ge 1$, each $m \ge 1$, and each $\ell \ge \max \{ r, m \}$, it follows from the definitions of 
$k_{[(\ell-1)(10\ell^2+7\ell+6)/6]+(\ell+1)^2+(2m-2)(\ell+1)+1}$ and the corresponding compact intervals 
$$
K_{\breve \al_0\breve \al_1\breve \al_2\cdots\breve \al_\ell}^{(j, \, (\ell-1)(10\ell^2+7\ell+6)/6)+(\ell+1)^2+(2m-2)(\ell+1)+1)} \, \big(\supset K_{\breve \al_0\breve \al_1\breve \al_2\cdots\breve \al_\ell}^{(j, \, \ell(10\ell^2+27\ell+23)/6)}\big), \, \al_n = 0, 1, \, 0 \le n \le \ell, \, 1 \le j \le \ell, 
$$
that
\begin{multline*}
$$
f^{k_{[(\ell-1)(10\ell^2+7\ell+6)/6]+(\ell+1)^2+(2m-2)(\ell+1)+1}}\big(\mathcal S^{(1)} \cup \mathcal S^{(2)} \cup \mathcal S^{(3)} \cup \cdots \cup \mathcal S^{(r)}\big) \\
\subset f^{k_{[(\ell-1)(10\ell^2+7\ell+6)/6]+(\ell+1)^2+(2m-2)(\ell+1)+1}}\big(\cup_{j=1}^\ell \, \big\{ K_{\breve \al_0\breve \al_1\breve \al_2\cdots\breve \al_\ell}^{(j, \, \ell(10\ell^2+27\ell+23)/6)}: \, \al_n = 0, 1, \, 0 \le n \le \ell \big\}\big) \qquad\, \\ 
\,\,\,\,\,\, \subset f^{k_{[(\ell-1)(10\ell^2+7\ell+6)/6]+(\ell+1)^2+(2m-2)(\ell+1)+1}}\big(\cup_{j=1}^\ell \, \big\{ K_{\breve \al_0\breve \al_1\breve \al_2\cdots\breve \al_\ell}^{(j, \, (\ell-1)(10\ell^2+7\ell+6)/6+(\ell+1)^2+(2m-2)(\ell+1)+1)}: \, \al_n = 0, 1, \\
\qquad\, 0 \le n \le \ell \big\}\big) \subset V_{n_\ell}(y_m) \,\,\, \text{and} \,\,\, f^{k_{[(\ell-1)(10\ell^2+7\ell+6)/6]+(\ell+1)^2+(2m-2)(\ell+1)+1}}(x_m) \in V_{n_\ell}(y_m) \,\,\, \text{and} \qquad\quad \\ 
$$
\end{multline*}
\indent Since the length of the interval $V_{n_\ell}(y_m)$ shrinks to zero as $n_\ell$ tends to $\infty$, we obtain that, for each $r \ge 1$, the union $\mathcal S^{(1)} \cup \mathcal S^{(2)} \cup \cdots \cup \mathcal S^{(r)}$ is {\it dynamically} synchronously proximal to the point $x_m$ for each $m \ge 1$.  

Furthermore, for each $r \ge 1$, each $m \ge 1$, each $t \ge 1$ and all $\ell \ge \max \{ r, m, t \}$, it follows from the definitions of 
$k_{[(\ell-1)(10\ell^2+7\ell+6)/6]+(\ell+1)^2+2\ell(\ell+1)+(2m-2)\ell+t}$ and the corresponding compact intervals 
$$
K_{\breve \al_0\breve \al_1\breve \al_2\cdots\breve \al_\ell}^{(j, \, (\ell-1)(10\ell^2+7\ell+6)/6)+(\ell+1)^2+2\ell(\ell+1)+(2m-2)\ell+t)} \big(\supset K_{\breve \al_0\breve \al_1\breve \al_2\cdots\breve \al_\ell}^{(j, \, \ell(10\ell^2+27\ell+23)/6)}\big), \al_n = 0, 1, 0 \le n \le \ell, 1 \le j \le \ell, 
$$
that
\begin{multline*}
$$
\,\,\, f^{k_{[(\ell-1)(10\ell^2+7\ell+6)/6]+(\ell+1)^2+2\ell(\ell+1)+(2m-2)\ell+t}}\big(\mathcal S^{(1)} \cup \mathcal S^{(2)} \cup \mathcal S^{(3)} \cup \cdots \cup \mathcal S^{(r)}\big) \qquad\qquad\qquad\qquad \\
\subset f^{k_{[(\ell-1)(10\ell^2+7\ell+6)/6]+(\ell+1)^2+2\ell(\ell+1)+(2m-2)\ell+t}}\big(\cup_{j=1}^\ell \, \big\{ K_{\breve \al_0\breve \al_1\breve \al_2\cdots\breve \al_\ell}^{(j, \, \ell(10\ell^2+27\ell+23)/3)}: \, \al_n = 0, 1, \, 0 \le n \le \ell \big\}\big) \\ 
\qquad\,\,\,\, \subset V_{n_\ell}(y_m) \,\,\, \text{and} \,\,\, f^{k_{[(\ell-1)(10\ell^2+7\ell+6)/6]+(\ell+1)^2+2\ell(\ell+1)+(2m-2)\ell+t}}\big(f^t(x_m)\big) \in V_{n_\ell}(y_m). \qquad\qquad\qquad\quad\quad
$$
\end{multline*}
\indent Since the length of the interval $V_{n_\ell}(y_m)$ shrinks to zero as $n_\ell$ tends to $\infty$, we obtain that, for each $r \ge 1$, the union $\mathcal S^{(1)} \cup \mathcal S^{(2)} \cup \cdots \cup \mathcal S^{(r)}$ is {\it dynamically} synchronously proximal to the point $f^t(x_m)$ for each $m \ge 1$ and each $t \ge 1$.  This, combined with the above, implies that the union $\mathcal S^{(1)} \cup \mathcal S^{(2)} \cup \cdots \cup \mathcal S^{(r)}$ is {\it dynamically} synchronously proximal to the point $f^i(x_m)$ for each $m \ge 1$ and each $i \ge 0$.  Therefore, Part (2) holds.   

Let $\mathbb S = \bigcup_{j=1}^\infty \, \mathcal S^{(j)}$.  Then, since $\mathcal S^{(j)} \subset U_j$ for each $j \ge 1$ and since $\{U_1, U_2, U_3, \cdots \}$ is a countable open base for $I$, the set $\mathbb S$ is dense in $I$.  For any two distinct points $c$ and $d$ in $\mathbb S = \bigcup_{j=1}^\infty \, \mathcal S^{(j)}$, since there exists a positive integer $\ell$ such that $\{ c, d \} \subset \mathcal S^{(1)} \cup \mathcal S^{(2)} \cup \cdots \cup \mathcal S^{(\ell)}$, the set $\{ c, d \}$, as just shown above, is $\mathcal M$-synchronouly proximal to, say, $v_1$ and so, we obtain that 
$$
\liminf_{\substack{n \to \infty \\ n \in \mathcal M}} \big|f^n(c) - f^n(d)\big| = 0.
$$
\indent On the other hand, to determine $\limsup_{\substack{n \to \infty \\ n \in \mathcal M}} \big|f^n(c) - f^n(d)\big|$, we have two cases to consider:

Case 1. Suppose both $c$ and $d$ belong to $\mathcal S^{(j)} = \{ x_{\omega_\al}^{(j)}: \al \in \Sigma_2 \}$ for some positive integer $j$. Then, $c = x_{\omega_{\be}}^{(j)}$ and $d = x_{\omega_{\gamma}}^{(j)}$ for some $\be \ne \gamma$ in $\Sigma_2$.  Write $\omega_{\be} = \breve \be_0 \breve \be_1 \breve \be_2 \breve \be_3 \cdots$ and $\omega_{\gamma} = \breve \gamma_0 \breve \gamma_1 \breve \gamma_2 \breve \gamma_3 \cdots$.  Then $\breve \be_\ell \ne \breve \gamma_\ell$ for infinitely many $\ell$'s.  For any $i \ge 0$ and any such $\ell \ge \max \{ i,j \}$, since 
$$
c = x_{\omega_{\be}}^{(j)} \in K_{\breve \be_0 \breve \be_1 \breve \be_2 \breve \be_3 \cdots\breve \be_\ell}^{(j, \, \ell(10\ell^2+27\ell+23)/6)} \,\,\, \text{and} \,\,\, d = x_{\omega_{\gamma}}^{(j)} \in K_{\breve \gamma_0 \breve \gamma_1 \breve \gamma_2 \breve \gamma_3 \cdots\breve \gamma_\ell}^{(j, \, \ell(10\ell^2+27\ell+23)/6)},
$$
it follows from the defintion of $k_{[(\ell-1)(10\ell^2+7\ell+6)/6]+1+i} \in \mathcal M$ (at the $\ell^{th}$ stage) that 
\begin{multline*}
$$
f^{k_{[(\ell-1)(10\ell^2+7\ell+6)/6)]+1+i}}\big(f^i(c)\big) \in W_{n_\ell}(\breve \be_\ell) \,\,\, \text{and} \,\,\, f^{k_{[(\ell-1)(10\ell^2+7\ell+6)/6]+1+i}}\big(f^i(d)\big) \in W_{n_\ell}(\breve \gamma_\ell), \\ \text{where} \,\,\, W_{n_\ell}(0) = V_{n_\ell}(a_\ell) \,\,\, \text{and} \,\,\,  W_{n_\ell}(1) = V_{n_\ell}(b_\ell).
$$
\end{multline*}
Since $\breve \be_\ell \ne \breve \gamma_\ell$, we have $\{ \breve \be_\ell, \breve \gamma_\ell \} = \{ 0, 1 \}$.  Since the intervals $V_{n_\ell}(a_\ell)$ and $V_{n_\ell}(b_\ell)$ are at least $|a_\ell - b_\ell| - \frac 3{n_\ell}$ apart, so are the points $f^{k_{[(\ell-1)(10\ell^2+7\ell+6)/6]+1+i}}\big(f^i(c)\big)$ and $f^{k_{[(\ell-1)(10\ell^2+7\ell+6)/6]+1+i}}\big(f^i(d)\big)$.  

Case 2. Suppose $c \in \mathcal S^{(j_1)}$ and $d \in \mathcal S^{(j_2)}$ with $1 \le j_1 < j_2$.  Then, for any $i \ge 0$ and any $\ell \ge \max \{ i, j_2 \}$, it follows from the definition of $k_{[(\ell-1)(10\ell^2+7\ell+6)/6]+j_1(\ell+1)+1+i} \in \mathcal M$ (at the $\ell^{th}$ stage) that 
$$
f^{k_{[(\ell-1)(10\ell^2+7\ell+6)/6]+j_1(\ell+1)+1+i}}\big(f^i(c)\big) \in V_{n_\ell}(a_\ell) \,\,\, \text{and} \,\,\, f^{k_{[(\ell-1)(10\ell^2+7\ell+6)/6]+j_1(\ell+1)+1+i}}\big(f^i(d)\big) \in V_{n_\ell}(b_\ell).
$$
So, the points $f^{k_{[(\ell-1)(10\ell^2+7\ell+6)/6]+j_1(\ell+1)+1+i}}\big(f^i(c)\big)$ and $f^{k_{[(\ell-1)(10\ell^2+7\ell+6)/6]+j_1(\ell+1)+1+i}}\big(f^i(d)\big)$ are at least $|a_\ell - b_\ell| - \frac 3{n_\ell}$ apart.  

By combining the results in Case 1 and Case 2 above with $i = 0$ and by the definition of the sequences $<V_{n_\ell}(a_n)>$ and $<V_{n_\ell}(b_n)>$, we obtain that 
$$
\limsup_{\substack{n \to \infty \\ n \in \mathcal M}} \big|f^n(c) - f^n(d)\big| = \beta = \begin{cases}
                       \infty, & \text{if $I$ is unbounded,} \\
                       |I|, & \text{if $I$ is bounded.} \\
                        \end{cases}
$$
This shows that 
$$
\mathbb S \,\,\, \text{is a dense} \,\,\, \beta\text{-scrambled set (with respect to} \,\,\, \mathcal M) \,\,\, \text{of} \,\,\, f \,\,\, \text{in} \,\,\, I.
$$
Therefore, Part (3) holds.

By combining the results in Case 1 and Case 2 above with $i \ge 0$ and by the definition of the sequences $<V_{n_\ell}(a_n)>$ and $<V_{n_\ell}(b_n)>$, we obtain that 
$$
\limsup_{\substack{n \to \infty \\ n \in \mathcal M}} \big|f^n\big(f^i(c)\big) - f^n\big(f^i(d)\big)\big| = \beta = \begin{cases}
                       \infty, & \text{if $I$ is unbounded,} \\
                       |I|, & \text{if $I$ is bounded.} \\
                        \end{cases}
$$
This shows that Part (5) holds.

Now, for any integers $r \ge 1$, $m \ge 1$, $s \ge 0$, and all $\ell \ge \max \{ r, m, s \}$, it follows from the definitions of 
$$
k_{[(\ell-1)(10\ell^2+7\ell+6)/6]+(\ell+1)^2+(2m-2)(\ell+1)+1+s} \,\,\, \text{and} \,\,\, k_{[(\ell-1)(10\ell^2+7\ell+6)/6]+(\ell+1)^2+(2m-1)(\ell+1)+1+s}
$$ 
that 
\begin{multline*}
$$
f^{k_{[(\ell-1)(10\ell^2+7\ell+6)/6]+(\ell+1)^2+(2m-2)(\ell+1)+1+s}}\big(\cup_{j=1}^\ell \big\{ f^s\big(K_{\breve \al_0\breve \al_1\breve \al_2\cdots\breve \al_\ell}^{(j, \ell(10\ell^2+27\ell+23)/6)}\big): \, \al_i = 0, 1, \, 0 \le i \le \ell \big\}\big) \\ 
\subset V_{n_\ell}(y_m) \,\,\, \text{and} \,\,\, f^{k_{[(\ell-1)(10\ell^2+7\ell+6)/6]+(\ell+1)^2+(2m-2)(\ell+1)+1+s}}(x_m) \in V_{n_\ell}(y_m) \,\,\, \text{and} \\
f^{k_{[(\ell-1)(10\ell^2+7\ell+6)/6]+(\ell+1)^2+(2m-1)(\ell+1)+1+s}}\big(\cup_{j=1}^\ell \big\{ f^s\big(K_{\breve \al_0\breve \al_1\breve \al_2\cdots\breve \al_\ell}^{(j, \ell(10\ell^2+27\ell+23)/6)}\big): \, \al_i = 0, 1, \, 0 \le i \le \ell \big\}\big) \qquad\,\,\,\, \\ 
\subset \widehat V_{n_\ell}(y_m) \,\,\, \text{and} \,\,\, f^{k_{[(\ell-1)(10\ell^2+7\ell+6)/6]+(\ell+1)^2+(2m-1)(\ell+1)+1+s}}(x_m) \in V_{n_\ell}(y_m). \qquad\qquad\qquad\,\,\,\,\,
$$
\end{multline*}
Therefore, we obtain that, for any point $x$ in $\mathcal S^{(r)}$, 
$$
f^{k_{[(\ell-1)(10\ell^2+7\ell+6)/6]+(\ell+1)^2+(2m-2)(\ell+1)+1+s}}\big(\{ f^s(x), x_m \}\big) \subset V_{n_{\ell}}(y_m) \,\,\, \text{and}, 
$$
\begin{multline*}
$$
f^{k_{[(\ell-1)(10\ell^2+7\ell+6)/6]+(\ell+1)^2+(2m-1)(\ell+1)+1+s}}\big(f^s(x)\big) \in \widehat V_{n_\ell}(y_m) \,\,\, \text{and} \\ f^{k_{[(\ell-1)(10\ell^2+7\ell+6)/6]+(\ell+1)^2+(2m-1)(\ell+1)+1+s}}(x_m) \in V_{n_\ell}(y_m).
$$
\end{multline*}
\indent Since the length of the interval $V_{n_\ell}(y_m)$ shrinks to zero as $n_\ell$ tends to $\infty$, we obtain that 
$$
\liminf_{n \to \infty} \big|f^n\big(f^s(x)\big) - f^n(x_m)\big| = 0.
$$
\indent Since, by the choices of $\widehat V_{n_\ell}(y_m)$'s and $V_{n_\ell}(y_m)$'s, we have $\lim_{\ell \to \infty} dist\big(\widehat V_{n_\ell}(y_m), V_{n_\ell}(y_m)\big) \ge \frac \beta2$ and so, we obtain that 
$$
\limsup_{n \to \infty} \big|f^n\big(f^s(x)\big) - f^n(x_m)\big| \ge \frac \beta2.\vspace{.1in}
$$
\indent Furthermore, for any integers $r \ge 1$, $m \ge 1$, $t \ge 1$ and all $\ell \ge \max \{ r, m, t \}$, it follows from the definitions of  
$$
k_{[(\ell-1)(10\ell^2+7\ell+6)/6]+(\ell+1)^2+2\ell(\ell+1)+(2m-2)\ell+t} \,\,\, \text{and} \,\,\, k_{[(\ell-1)(10\ell^2+7\ell+6)/6]+(\ell+1)^2+2\ell(\ell+1)+(2m-1)\ell+t}
$$ 
that 
\begin{multline*}
$$
f^{k_{[(\ell-1)(10\ell^2+7\ell+6)/6]+(\ell+1)^2+2\ell(\ell+1)+(2m-2)\ell+t}}\bigg(\bigcup_{j=1}^\ell \big\{ K_{\breve \al_0\breve \al_1\breve \al_2\cdots\breve \al_\ell}^{(j, \ell(10\ell^2+27\ell+23)/6)}: \, \al_i = 0, 1, \, 0 \le i \le \ell \big\}\bigg) \\ 
\subset V_{n_\ell}(y_m) \,\,\, \text{and} \,\,\, f^{k_{[(\ell-1)(10\ell^2+7\ell+6)/6]+(\ell+1)^2+2\ell(\ell+1)+(2m-2)\ell+t}}\big(f^t(x_m)\big) \in V_{n_\ell}(y_m) \,\,\, \text{and} \\
 f^{k_{[(\ell-1)(10\ell^2+7\ell+6)/6]+(\ell+1)^2+2\ell(\ell+1)+(2m-1)\ell+t}}\big(\cup_{j=1}^\ell \big\{ K_{\breve \al_0\breve \al_1\breve \al_2\cdots\breve \al_\ell}^{(j, \ell(10\ell^2+27\ell+23)/6)}\big): \, \al_i = 0, 1, \, 0 \le i \le \ell \big\}\big) \qquad\quad\,\,\, \\ 
\subset \widehat V_{n_\ell}(y_m) \,\,\, \text{and} \,\,\, f^{k_{[(\ell-1)(10\ell^2+7\ell+6)/6]+(\ell+1)^2+2\ell(\ell+1)+(2m-1)\ell+t}}\big(f^t(x_m)\big) \in V_{n_\ell}(y_m). \qquad\qquad\quad\,\,\,\,
$$
\end{multline*}
Therefore, we obtain that, for any point $x$ in $\mathcal S^{(r)}$,  
$$
f^{k_{[(\ell-1)(10\ell^2+7\ell+6)/6]+(\ell+1)^2+2\ell(\ell+1)+(2m-2)\ell+t}}\big(\{ x, f^t(x_m) \}\big) \subset V_{n_{\ell}}(y_m) \,\,\, \text{and},
$$
\begin{multline*}
$$
f^{k_{[(\ell-1)(10\ell^2+7\ell+6)/6]+(\ell+1)^2+2\ell(\ell+1)+(2m-1)\ell+t}}(x) \in \widehat V_{n_\ell}(y_m) \,\,\, \text{and} \\ f^{k_{[(\ell-1)(10\ell^2+7\ell+6)/6]+(\ell+1)^2+2\ell(\ell+1)+(2m-1)\ell+t}}(f^t\big(x_m)\big) \in V_{n_\ell}(y_m).
$$
\end{multline*}
\indent Since the length of the interval $V_{n_\ell}(y_m)$ shrinks to zero as $n_\ell$ tends to $\infty$, we obtain that 
$$
\liminf_{n \to \infty} \big|f^n(x) - f^n\big(f^t(x_m)\big)\big| = 0.
$$
\indent Since, by the choices of $\widehat V_{n_\ell}(y_m)$'s, we have $\lim_{\ell \to \infty} dist\big(\widehat V_{n_\ell}(y_m), V_{n_\ell}(y_m)\big) \ge \frac \beta2$ and so, we obtain that 
$$
\limsup_{n \to \infty} \big|f^n(x) - f^n\big(f^t(x_m)\big)\big| \ge \frac \beta2.
$$
\indent This, combined with the above, implies Part (4).

Now we show that $\delta = \inf_{n \ge 1} \big\{ \sup \{ \big|f^n(x) - x\big|: x \in I \} \big\} > 0$.  

Let $K$ and $K'$ be two disjoint compact intervals in $I$.  Then, since $f$ is mixing, there exists a positive integer $N$ such that 
$$
f^n(K) \cap K' \ne \emptyset \,\,\, \text{for all} \,\,\, n \ge N.
$$
Consequently, for each $n \ge N$, the number $\delta_n = \sup \big\{ \big|f^n(x) - x\big|: x \in I \big\} \ge dist(K, K') > 0$.  

On the other hand, let $u$ be a transitive point of $f$.  Then $f^n(u) \ne u$ for all $n \ge 1$.  In particular, for each $1 \le n \le N-1$, the number $\delta_n = \sup \{ \big|f^n(x) - x\big|: x \in I \} \ge \big|f^n(u) - u\big| > 0$.  Therefore, we obtain that 
$$
\delta = \inf_{n \ge 1} \big\{ \sup_{x \in I} \big\{ \big|f^n(x) - x\big| \big\} \big\} \ge \min \big\{ \big|f(u) - u\big|, \big|f^2(u) - u\big|, \cdots, \big|f^{N-1}(u) - u\big|, dist(K, K') \big\} > 0.
$$

Finally, let $z$ be a fixed point of $f$.  Then, since $\delta > 0$, $f$ is not uniformly rigid and it follows from {\bf\cite{fo}} that $f$ has an invariant $\va$-scrambled set for some positive number $\va$.  Here, we show that $f$ has s dense invariant $\delta$-scrambled set (with respect to $\mathcal M$) of totally transitive points of $f$ in $I$.

For any positive integers $n$ and $k$, let $c_{n,k}$ be a point in $I$ such that $\big|f^n(c_{n,k}) - c_{n,k} \big| > \delta_n - 1/k$.  Since $f$ is a continuous mixing map on $I$, it follows from the above that $f$ has totally transitive points.  Therefore, for each $m \ge 1$, $f^m$ is transitive on $I$ and so, $f^m(I)$ is dense in $I$.  Consequently, for each $n \ge 1$, each $m \ge 1$ and each $k \ge 1$, there exists a point $c_{n,m,k}$ in $I$ such that the point $f^m(c_{n,m,k})$ is so close to $c_{n,k}$ that 
$$
\big|f^n\big(f^m(c_{n,m,k})\big) - f^m(c_{n,m,k})\big| > \delta_n - 1/k.
$$
For each $n \ge 1$ and each $k \ge 1$, let $c_{n,0,k} = c_{n,k}$.

Let $v_1 = z$. We arrange these countably infinitely many points $c_{n,m,k}, \, n \ge 1$, $m \ge 0$, $k \ge 1$, in a sequence and call this sequence 
$$
v_2, \, v_3, \, v_4, \, \cdots \,\,\, \text{(note that we have defined} \,\,\, v_1 = z)
$$
and let $\mathbb S' = \bigcup_{t \ge 1} S'^{(t)}$ be the resulting $\beta$-scrambled set (with respect to $\mathcal M$) of totally transitive points of $f$ obtained in the theorem.   

Let $c$ and $d$ be any points in $\mathbb S'$ and let $j > i \ge 0$ be any integers.  

Since, for each $\ell \ge 1$, the union $\mathcal S'^{(1)} \cup \mathcal S'^{(2)} \cup \mathcal S'^{(3)} \cup \cdots \cup \mathcal S'^{(\ell)}
$ is $\mathcal M$-synchronously proximal to $v_t$ for each $t \ge 1$, so is the set $\{ c, d \}$.  

Since $v_1 = z$ is a fixed point of $f$ and since the set $\{ c, d \}$ is $\mathcal M$-synchronously proximal to $v_1 = z$, we easily obtain that 
\begin{equation*}
\liminf_{\substack{n \to \infty \\ n \in \mathcal M}} \big|f^n\big(f^i(c)\big) - f^n\big(f^j(d)\big)\big| = 0. \tag{$\dagger$}
\end{equation*}
\indent On the other hand, since the set $\{ c, d \}$ is $\mathcal M$-synchronously proximal to $v_t$ for each $t \ge 2$, $\{ c, d \}$ is $\mathcal M$-synchronously proximal to $c_{n,m,k}$ for each $n \ge 1$, each $m \ge 0$, and each $k \ge 1$.  Consequently, for each $k \ge 1$, there exists a strictly increasing sequence $r_{1,k}, r_{2,k}, r_{3,k}, \cdots$ of positive integers, depending on $i$ and $j$, in $\mathcal M$ such that $\lim_{\ell \to \infty} f^{r_{\ell,k}}(c) = f^{r_{\ell,k}}(d) = c_{j-i, i,k}$.  In particular, for each $k \ge 1$, we have    
\begin{multline*}
$$
\limsup_{\substack{n \to \infty \\ n \in \mathcal M}} \big|f^n\big(f^j(c)\big) - f^n\big(f^i(d)\big)\big| \ge \lim_{\ell \to \infty} \big|f^{r_{\ell,k}}\big(f^j(c)\big) - f^{r_{\ell,k}}\big(f^i(d)\big)\big| \\
= \lim_{\ell \to \infty} \big|f^j\big(f^{r_{\ell,k}}(c)\big) - f^i\big(f^{r_{\ell,k}}(d)\big)\big|  
= \big|f^j\big(c_{j-i,i,k}\big) - f^i\big(c_{j-i,i,k}\big)\big| = \big|f^{j-i}\big(f^i(c_{j-i,i,k})\big) - f^i(c_{j-i,i,k})\big| \\
> \delta_{j-i} \, - \, \frac 1k \ge \delta \, - \, \frac 1k
$$
\end{multline*}
which implies that 
$$
\limsup_{\substack{n \to \infty \\ n \in \mathcal M}} \big|f^n\big(f^j(c)\big) - f^n\big(f^i(d)\big)\big| \ge \delta.
$$
\indent Furthermore, if $c \ne d$, then, it follows from Part (5) that
$$
\limsup_{\substack{n \to \infty \\ n \in \mathcal M}} \big|f^n\big(f^i(c)\big) - f^n\big(f^i(d)\big)\big| = \beta \ge \delta.
$$
\indent This, combined with the above case when $j > i \ge 0$ and the above $(\dagger)$ with $\liminf$, implies that the set 
$\widehat {\mathbb S'} = \bigcup_{s=0}^\infty \, f^s(\mathbb S')$
is a dense {\it invariant} $\delta$-scrambled set (with respect to $\mathcal M$) of totally transitive points of $f$ in $I$.  (Since each point $u$ of $\mathbb S'$ is a totally transitive point of $f$, for each $n \ge 1$, the set $\{ (f^n)^k(u): k \ge 0 \}$ is dense in $I$ and so is the set $f^s\big(\{ (f^n)^k(u): k \ge 0 \}\big) = \{ (f^n)^k\big(f^s(u)\big): k \ge 0 \}$ for each integer $s \ge 0$).
\hfill\sq

For continuous bitransitive maps $f$ on an interval, Theorem 8 guarantees the existence of dense extremely scrambled sets of totally transitive points of $f$.  We may ask if there exist dense {\it invariant} extremely scrambled sets of totally transitive points of $f$.  In the following, we provide various examples to demonstrate all possibilities (some examples are mimicking those in {\bf\cite{bo, c}}):  

\noindent
{\bf Example 1.} Let $I = (-\infty, \infty)$ and let $1/3 > \va_0 > \va_1 > \va_2 > \cdots$ and $1/3 > \va_{-1} > \va_{-2} > \va_{-3} > \cdots$ be two strictly decreasing sequences of positive numbers such that $\lim_{n \to \infty} \va_n = 0 = \lim_{n \to \infty} \va_{-n}$.  Let $f(x) : (-\infty, \infty) \rightarrow (-\infty, \infty)$ be the continuous map defined by putting, for each integer $i$, 
$$
\begin{cases}
(1) \,\, f(i) = i, \, f(i + 1/3) = i+1 + \va_i, \, f(i + 2/3) = i - \va_i, \, f(i+1) = i+1 \,\,\, \text{and} \cr
(2) \,\, \text{$f$ is linear on each of the intervals: $[i, i+1/3], [i+1/3, i+2/3]$ and $[i+2/3, i+1]$}.\cr
\end{cases}
$$
\indent Let $K$ be any compact interval in $(-\infty, \infty)$.  If $K$ contains the point $i+1/3$ or the point $i+2/3$ for some integer $i$, then the lengths ratio $|f(K)|/|K|$ is $> 1.5$.  Otherwise, the lengths ratio $|f(K)|/|K|$ is $> 3$.  Consequently, for some positive integer $j$, the length $|f^j(K)|$ is $> 1$ and so, $f^j(K)$ contains the interval $[k. k+1]$ for some integer $k$.  But then, since $f\big([k, k+1]\big) \supsetneq [k, k+1]$, we see that, for any positive integer $n$, there exists a positive integer $m$ such that $f^\ell([k, k+1]) \supset [-n, n]$ for all $\ell \ge m$.  Therefore, $f$ is a continuous bitransitive map on $I = (-\infty, \infty)$.  However, for every point $x$, we have $|f(x) - x| < 4/3$.  In particular, if $\{ x_0, f(x_0) \}$ is a $\delta$-scrambled set of $f$, then $\delta \le \limsup_{n \to \infty} |f^n\big(f(x_0)\big)-f^n(x_0)| < 4/3$.  This implies that $f(x)$ cannot have {\it invariant} $\infty$-scrambled sets on the unbounded interval $I = (-\infty, \infty)$ although, by Theorem 8, it has dense $\infty$-scrambled sets of totally transitive points in $I = (-\infty, \infty)$.  

\noindent
{\bf Example 2.} Let $I = (-\infty, \infty)$ and let $\hat f(x) : (-\infty, \infty) \to (-\infty, \infty)$ be the continuous symmetric map defined by putting 
$$
\begin{cases}
(1) \,\, \hat f(0) = 0, \cr
(2) \,\, \hat f(m) = -(m+2) \,\,\, \text{for all positive odd integer $m$}, \cr
(3) \,\, \hat f(n) = 1 \,\,\, \text{for all positive even integer $n$}, \cr
(4) \,\, \hat f(x) \,\,\, \text{is linear on the interval $[k, k+1]$ for each integer $k \ge 0$}, \cr
(5) \,\, \hat f(x) = -\hat f(-x) \,\,\, \text{for all $x \le 0$}. \cr
\end{cases}
$$
\indent Let $K$ be any compact interval in $(-\infty, \infty)$.  If $K$ contains the point $i$ for some integer $i$, then the lengths ratio $|\hat f(K)|/|K|$ is $> 1.5$.  Otherwise, the lengths ratio $|\hat f(K)|/|K|$ is $\ge 3$.  Consequently, for some positive integer $j$, the length of $|\hat f^j(K)|$ is $> 1$ and so, ${\hat f}^j(K)$ contains the interval $[k. k+1]$ for some integer $k$.  Since 
$$
\begin{cases}
(a) \,\, \hat f^2\big([-1, 0]\big) = \hat f\big([0, 3]\big) \supset \hat f\big([1, 2]\big) = [-3, 1] \supset [-1, 1], \\
(b) \,\, \hat f^2\big([0, 1]\big) = \hat f\big([-3, 0]\big) \supset \hat f\big([-2, -1]\big) \supset [-1, 3] \supset [-1, 1], \\
(c) \,\, \hat f\big([k, k+1]\big) \supset [-1, 1] \,\,\, \text{for each $k \ne -1, 0$, and} \\
(d) \,\, \hat f\big([-m, m]\big) = [-(m+2), m+2] \,\,\, \text{for each positive odd integer $m$},\\
\end{cases}
$$
we see that, for any integer $n \ge 1$, $\hat f^\ell\big([k, k+1]\big) \supset [-n, n]$ for all sufficiently large integer $\ell$.  So, $\hat f$ is a continuous bitransitive map on $I = (-\infty, \infty)$ such that the point $x = 0$ is the {\it unique} fixed point of $\hat f(x)$ and, for each integer $n \ge 1$, there is a non-recurrent point $v_n$ in the interval $[n, n+1]$ such that $\hat f(v_n) = 0$ and so, the $\delta$ defined in Theorem 8 equals $\infty$.  Therefore, $\hat f(x)$ has dense {\it invariant} $\infty$-scrambled sets of totally transitive points but has no {\it invariant} extremely scrambled sets \big(since, for each $x_0 \ge 0$, $\big\{ x_0, \hat f(x_0) \big\} \subset [1, \infty)$ while, for each $x_0 \le 0$, $\big\{ x_0, \hat f(x_0) \big\} \subset [-1, \infty)$\big) in $I = (-\infty, \infty)$.  

\noindent
{\bf Example 3.} Let $I = (0, \infty)$ and let $g(x) : (0, \infty) \to (0, \infty)$ be the continuous map defined by putting, for all positive integers $m$ and all positive integers $n$, 
$$
\begin{cases}
(1) \,\, g(1/2^m) = 1/2^m, \\
(2) \,\, g\big(1/2^m + 1/(3 \cdot 2^m)\big) = 1/2^{m-1} + 1/(3 \cdot 2^m), \\
(3) \,\, g\big(1/2^m + 2/(3 \cdot 2^m)\big) = 1/2^m - 1/(3 \cdot 2^m), \\
(4) \,\, \text{$g(x)$ is linear on $[1/2^{k+1}, 1/2^k]$ for each integer $k \ge 0$}, \\
(5) \,\, g(n) = n, \\
(6) \,\, g(n+\frac 13) = n+1+\frac 13, \\
(7) \,\, g(n+\frac 23) = n - \frac 13, \, \text{and} \\
(8) \,\, \text{$g(x)$ is linear on $[n, n +1]$}. \\ 
\end{cases}
$$
Then $g(x)$ is a continuous bitransitive map from $(0, \infty)$ onto itself such that $|g(x) - x| < 2$ for all $x > 0$.  Therefore, $g(x)$ cannot have {\it invariant} $\infty$-scrambled sets in $I = (0, \infty)$ although, by Theorem 8, it has dense $\infty$-scrambled sets of totally transitive points.  Furthermore, if we let $g(0) = 0$, then $g(x)$ is a continuous bitransitive map from $[0, \infty)$ onto itself which does not have {\it invariant} $\infty$-scrambled sets on $I = [0, \infty)$.  By an easy conjugation and modification, this also shows that there exist continuous bitransitive maps on {\it any} bounded interval (compact or not, if compact then both endpoints are fixed points) which does not have {\it invariant} extremely scrambled sets in $I$.  

\noindent
{\bf Example 4.} Let $I = (0, \infty)$ and let $\hat g(x) : (0, \infty) \rightarrow (0, \infty)$ be the continuous map defined by putting, for all positive integers $m$, 
$$
\begin{cases}
(1) \,\, \hat g(1/2^m) = 1/2^m, \\
(2) \,\, \hat g\big(1/2^m + 1/(3 \cdot 2^m)\big) = 1/2^{m-1} + 1/(3 \cdot 2^m), \\
(3) \,\, \hat g\big(1/2^m + 2/(3 \cdot 2^m)\big) = 1/2^m - 1/(3 \cdot 2^m), \\
(4) \,\, \text{$\hat g(x)$ is linear on $[1/2^{i+1}, 1/2^i]$ for each integer $i \ge 0$}, \\
(5) \,\, \hat g(1) = 1, \\
(6) \,\, \text{$\hat g(n) = 5/6$ for all odd integers $n \ge 3$}, \\
(7) \,\, \text{$\hat g(n) = n+2$ for all even integers $n \ge 2$, \, and} \\ 
(8) \,\, \text{$\hat g(x)$ is linear on $[j, j +1]$ for each positive integer $j$}.\\  
\end{cases}
$$
Then $\hat g(x)$ is a continuous bitransitive map from $(0, \infty)$ onto itself.  Since the point $x = 1$ is a fixed point of $\hat g(x)$ and since, there exists a point $v_n$ in the interval $[n, n+1]$, for each even $n \ge 2$, such that $\hat g(v_n) = 1$, the $\delta$ in Theorem 8 equals $\infty$.  So, $\hat g(x)$ has dense {\it invariant} $\infty$-scrambled sets of totally transitive points but has no {\it invariant} extremely scrambled sets in $I = (0, \infty)$ \big(since, for each $x_0 \ge 1$, $\big\{ x_0, \hat f(x_0) \big\} \subset [5/6, \infty)$ while, for each $0 \le x_0 \le 1$, $\big\{ x_0, \hat f(x_0) \big\} \subset [0, 2]$\big).  By an easy conjugation, this also shows that there exist continuous bitransitive maps on bounded open interval or bounded half-open interval which does not have {\it invariant} extremely scrambled sets.  

\noindent
{\bf Example 5.} Let $I = (-\infty, \infty)$ and let $\breve f(x) : I \to I$ be the continuous symmetric map defined by putting 
$$
\begin{cases}
(1) \,\, \breve f(0) = 0, \\
(2) \,\, \text{$\breve f(m) = 2m+1$ for all odd integers $m \ge 1$}, \\
(3) \,\, \text{$\breve f(n) = -2n$ for all even integers $n \ge 2$}, \\
(4) \,\, \text{$\breve f(x)$ is linear on $[k, k+1]$ for each integer $k \ge 0$, and} \\
(5) \,\, \text{$\breve f(x) = -\tilde f(-x)$ for all $x \le 0$}.\\
\end{cases}
$$
By arguing as those in Example 1, we obtain that $\breve f$ is a continuous bitransitive map on $I = (-\infty, \infty)$.  Let $<z_{-n}>_{n \ge 1}$ be a strictly decreasing unbounded sequence of fixed points of $\breve f(x)$ in the interval $(-\infty, -2]$ and let $<v_n>_{n \ge 1}$ be a strictly increasing unbounded sequence of {\it non-recurrent} points in $[2, \infty)$ such that $\tilde f(v_n) = z_{-n}$ for all positive integers $n$.  Then, the set $\mathbb S$ in Theorem 8 can be chosen so that $\widehat {\mathbb S}$ is also a dense uncountable {\it invariant extremely scrambled} set of totally transitive points of $\tilde f$ in $I = (-\infty, \infty)$.  

\noindent
{\bf Example 6.} Let $I = (0, \infty)$ and let $h(x) : (0, \infty) \rightarrow (0, \infty)$ be the continuous map defined by putting 
$$
\begin{cases}
(1) \,\, \text{$h(x) = 3x$ for $0 < x \le 1$}, \\
(2) \,\, \text{$h(m) = 3^{-m/2}$ for each even integer $m \ge 2$}, \\
(3) \,\, \text{$h(n) = n+2$ for each odd integer $n \ge 3$, and} \\
(4) \,\, \text{$h(x)$ is linear on $[k, k+1]$ for each integer $k \ge 1$}. \\
\end{cases}
$$
Then, by arguing as those in Example 1, we obtain that $h(x)$ is a continuous bitransitive map on $I = (0, \infty)$.  $h(x)$ is also a continuous bitransitive map from the closed unbounded interval $[0, \infty)$ onto itself if we let $h(0) = 0$ and so, the point $z = 0$ is the {\it unique} point in $[0, \infty)$ such that $h(z) = 0$. In fact, for the map $h(x)$, whether defined on $(0, \infty)$ or on $[0, \infty)$, there exist a strictly increasing {\it unbounded} sequence $<w_n>_{n \ge 1}$ of points such that $\lim_{n \to \infty} h(w_n) = 0$.  Therefore, the extremely scrambled set $\mathbb S$ in Theorem 8 can be chosen to satisfy that the larger set $\widehat {\mathbb S}$ is {\it also} a dense {\it invariant} extremal scrambled set of totally transitive points of $h$.  By an easy conjugation, this also shows that there exist continuous bitransitive maps on bounded open interval or bounded half-open interval which have dense {\it invariant} extremely scrambled sets of totally transitive points.  For continuous bitransitive map on compact intervals, the tent map $T(x) = 2x$ for $0 \le x \le 1/2$; and $T(x) = 2 - 2x$ for $1/2 \le x \le 1$ is such an example that it has dense {\it invariant} extremely scrambled sets of totally transitive points in $[0, 1]$.    

In the following, we apply Theorems 4 $\&$ 8 to show that every continuous transitive map on $I$ has a dense uncountable {\it invariant} $\delta$-scrambled set of transitive points for some positive number $\delta$ (cf. {\bf\cite{bal, liye, yuan}}).

\noindent
{\bf Corollary 9.}
{\it Assume that $f : I \longrightarrow I$ is a continuous transitive map.  Then $f$ has a dense uncountable {\it invariant} $\delta$-scrambled set of transitive points in $I$, where
$$
\delta = \begin{cases}
               \inf_{n \ge 1} \big\{\sup_{x \in I} |f^n(x)-x| \big\}, \,\, \text{if $f^2$ is also transitive on $I$}, \cr
               \max \big\{\inf_{n \ge 1} \{\sup_{x \in I \cap (-\infty, \breve z]} |f^{2n}(x)-x| \}, \, \inf_{n \ge 1} \{\sup_{x \in I                 \cap [\breve z, \infty)} |f^{2n}(x)-x| \}\big\}, \\
               \hspace*{1.2in} \text{if $f^2$ is not transitive on $I$}, \text{and $\breve z$ is the unique fixed point of $f$ in $I$}. \cr
       \end{cases}
$$}
\noindent
\!\!{\it Proof.}
If $f^2$ is transitive on $I$, then the existence of dense uncountable {\it invariant} $\delta$-scrambled sets of transitive points of $f$ in $I$ follows from Theorem 8.  So, in the following, we assume that $f$ is transitive but $f^2$ is not transitive on $I$.  Then by Theorem 4, there exists a unique fixed point $\breve z$ of $f$ (which lies in the interior of $I$) such that \\
[10pt]\centerline{$f\bigl(I \cap (-\infty, \breve z]\bigr) \subset I \cap [\breve z, \infty) \,\,\, \text{and} \,\,\, f\bigl(I \cap [\breve z, \infty)\bigr) \subset I \cap (-\infty, \breve z]$} \\
[10pt]and, on each of the intervals $I \cap (-\infty, \breve z]$ and $I \cap [\breve z, \infty)$, $f^2$ is transitive and has at least two fixed points (the fixed point $\breve z$ and a period-2 point of $f$) and so $f^2$ is {\it bitransitive}.  Consequently, it follows from Theorem 8 that both $\inf_{n \ge 1} \big\{ \sup_{x \in I \cap (-\infty, \breve z]} |f^{2n}(x)-x| \big\}$ and $\inf_{n \ge 1} \big\{ \sup_{x \in I \cap [\breve z, \infty)} |f^{2n}(x)-x| \big\}$ are positive.  Without loss of generality, we may assume that \\
[10pt]\centerline{$(0 <) \inf_{n \ge 1} \big\{ \sup_{x \in I \cap (-\infty, \breve z]} |f^{2n}(x)-x| \big\} \le \inf_{n \ge 1} \big\{ \sup_{x \in I \cap [\breve z, \infty)} |f^{2n}(x)-x| \big\}$} \\
[10pt]and consider the map $f^2$ on $I \cap [\breve z, \infty)$.  So, $\delta = \inf_{n \ge 1} \big\{ \sup_{x \in I \cap [\breve z, \infty)} |f^{2n}(x)-x| \big\}$ in this case.  Then, by choosing the fixed point $z$ in the proof of Theorem 8 to be the {\it unique} fixed point $\breve z$ of $f$ in the interior of $I$ and since $f^2$ is bitransitive on $I \cap [\breve z, \infty)$, Theorem 8 guarantees that $f^2$ has a dense uncountable $(\sup I- \breve z)$-scrambled set $\mathbb S$ of totally transitive points of $f^2$ in $I \cap [\breve z, \infty)$ such that (a) the set $\widehat {\mathbb S} = \bigcup_{i \ge 0} f^{2i}(\mathbb S)$ is an invariant $\delta$-scrambled set of $f^2$ \big(and so $f^2(\widehat {\mathbb S}) \subset \widehat {\mathbb S}\big)$ in $I \cap [\breve z, \infty)$ and (b) for any two distinct points $u_1$ and $u_2$ in $\mathbb S$, the set $\{ u_1, u_2 \}$ is synchronously proximal (with respect to $f^2$) to the unique fixed point $\breve z$ of $f$ and so, for any two distinct points $v_1$ and $v_2$ in $\widehat {\mathbb S}$, the set $\{ v_1, v_2 \}$ is synchronously proximal (with respect to $f^2$) to the unique fixed point $\breve z$ of $f$.  

In the following, we show that the set $\bigcup_{i \ge 0} f^i(\widehat {\mathbb S}) = \widehat {\mathbb S} \cup f(\widehat {\mathbb S})$ is a dense uncountable {\it invariant} $\delta$-scrambled set of transitive points of $f$ in $I$. 

Let $u$ and $u'$ be any two distinct points in $\widehat {\mathbb S} \cup f(\widehat {\mathbb S})$.

If both $u$ and $u'$ lie in the set $\widehat {\mathbb S}$, then since $\widehat {\mathbb S}$ is a $\delta$-scrambled set of $f^2$ in $I \cap [\breve z, \infty)$, we obtain that $\{ u, u' \}$ is a $\delta$-scrambled set of $f$.

If both $u$ and $u'$ lie in the set $f(\widehat {\mathbb S})$, then there exist two distinct points $\hat u$ and $\hat u'$ in $\widehat {\mathbb S} \subset I \cap [\breve z, \infty)$ such that $f(\hat u) = u$ and $f(\hat u') = u'$.  Since $\{ \hat u, \hat u' \} (\subset \widehat {\mathbb S})$ is a $\delta$-scrambled set of $f^2$, the set $\big\{ f(\hat u), f(\hat u') \big\} \big(= \{ u, u' \}\big)$ is a $\delta$-scrambled set of $f$.  

Now suppose $u \in \widehat {\mathbb S}$ and $u' \in f(\widehat {\mathbb S})$.  Let $\hat u'$ be a point in $\widehat {\mathbb S}$ such that $f(\hat u') = u'$.  So, both $u$ and $\hat u'$ lie in $\widehat {\mathbb S}$.  It follows from (b) above that $\big\{ u, \hat u' \big\}$ is synchronously proximal (with respect to $f^2$) to the fixed point $\breve z$ of $f$.  Consequently, we have 
$$
0 \le \liminf_{n \to \infty} \big|f^n(u) - f^n(u')\big| \le \liminf_{n \to \infty} \big|f^{2n}(u) - f^{2n}(\hat u')\big| = \liminf_{n \to \infty} \big|f^{2n}(u) - f\big(f^{2n}(\hat u')\big)\big| = 0.
$$
On the other hand, it is easy to see that, for each $i \ge 0$, $f^i(u)$ and $f^i(u')$ lie on opposite sides of $\breve z$.  Since the orbit $O_{f^2}(u)$ of $u$ with respect to $f^2$ is dense in $I \cap [\breve z, \infty)$ and since the orbit $O_{f^2}(u')$ of $u'$ with respect to $f^2$ is contained entirely in $I \cap (-\infty, \breve z]$, we easily obtain that 
$$
\limsup_{n \to \infty} \big|f^n(u) - f^n(u')\big| \ge \limsup_{n \to \infty} \big|f^{2n}(u) - f^{2n}(u')\big| \ge \sup I - \breve z \ge \delta.
$$
This implies that the set $\{ u, u' \}$ is a $\delta$-scrambled set of $f$.

Therefore, the set $\mathbb S \cup f(\mathbb S)$ is a dense uncountable {\it invariant} $\delta$-scrambled set of transitive points of $f$ in $I$.
\hfill\sq

\noindent
{\bf Remark 4.}
We can generalize Theorem 8 for continuous {\it bitransitive} maps on an interval to continuous {\it weakly mixing} maps and continuous {\it mixing} maps on infinite separable locally compact metric spaces without isolated points.  The proofs of these three results are almost the same except that we use different lemmas to prove them.  We refer to {\bf\cite{du5}} for details.  

\noindent
{\bf Remark 5.}
In {\bf\cite{du3, du4, du5, v2, v1}} (see also {\bf\cite{wang}}), we call a map $f$ {\it chaotic} if there exists a positive number $\va$ such that for every point $x$ and $$\text{every nonempty open set} \,\,\, V \,\,\, \text{(not necessarily an open neighborhood of} \,\,\, x),$$ there is a point $y$ in $V$ such that $$\liminf_{n \to \infty} \big|f^n(x) - f^n(y)\big| = 0 \,\,\, \text{and} \,\,\, \limsup_{n \to \infty} \big|f^n(x) - f^n(y)\big| \ge \va.$$ Therefore, continuous bitransitive maps on a bounded interval are {\it chaotic} in this sense while $$\text{not all continuous transitive maps on an interval are {\it chaotic}}$$ is demonstrated by the transitive map 
$\hat f : [0, 1] \to [0, 1]$ defined by 
$$
\hat f(x) = \begin{cases}
            1/2 + 2x & \text{for $0 \le x \le 1/4$}, \\
            -2x + 3/2 & \text{for $1/4 \le x \le 1/2$, and} \\
            -x + 1 & \text{for $1/2 \le x \le 1$}. \\
            \end{cases}
$$ together with the period-2 point $x = 1/6$ and the open interval $(1/2, 1)$.
The following is a {\it chaotic} map which is {\it not} a transitive map and yet has similar properties as those stated in Theorem 8.

\noindent
{\bf Example 7.}
{\it Assume that $g(x)$ is a continuous map from $[0, 3]$ onto itself defined by  
$$
g(x) =  \begin{cases}
               2x, & \text{if \quad $0 \le x \le 1$}, \qquad\qquad\qquad \qquad\qquad\qquad \cr
               4-2x, & \text{if \quad $1 \le x \le 2$}, \qquad\qquad\qquad \qquad\qquad\qquad \cr
               2x-4, & \text{if \quad $2 \le x \le 3$}. \qquad\qquad\qquad \qquad\qquad\qquad \cr
       \end{cases}
$$
Let $X = \{ x_i : i = 1, 2, \cdots \}$ be a countably infinite set of points in $[0, 3]$.  For each positive integer $i$, let $y_i$ be an $\omega$-limit point of $x_i$ and let $\rho$ be a real number such that 
$$
1 \le \rho \le \inf \big\{ \max \big\{ |y_i - 0|, |2-y_i| \big\} : i \ge 1 \big\}.
$$  
Then, there exists a dense uncountable 2-scrambled set $\mathcal Y$ of $g$ in $[0, 3]$ which is a countably infinite union of synchronously proximal Cantor sets $C_i, i \ge 1$ such that each nonempty open interval in $[0, 3]$ contains counatbly infinite many such Cantor sets $C_{k_i}, i \ge 1$ and the larger set $\widehat {\mathcal Y} = \bigcup_{i \ge 0} \, f^i(\mathcal Y)$ is a dense uncountable {\rm invariant} 2-scrambled set of $g$ in $[0, 3]$ with the property that, for any points $x$ in $X$ and $y$ in $\widehat {\mathcal Y}$, we have $\liminf_{n \to \infty} \big|g^n(x)-g^n(y)\big| = 0$ and $\limsup_{n \to \infty} \big|g^n(x)-g^n(y)\big| \ge \rho$.}

\noindent
\begin{proof}
For each point $x$ in the interval $[1, 3]$, let $x^* = 4 - x$.  Since $(x+x^*)/2 = 2$, the two points $x$ and $x^*$ are symmetric with respect to the point 2 and so, $g(x) = g(x^*)$.  In particular, for each $x \in [1, 3]$, $x$ and $x^*$ have the same $\omega$-limit sets.  

For any nonempty set $A$ of points in the interval $[0, 3]$, let $$A' = \bigl(A \cap [0, 2]\bigr) \cup \big\{ x^* : x \in A \,\,\, \text{and} \,\,\, x \in [2, 3] \big\} \subset [0, 2].$$  Let $X = \{ x_i : i = 1, 2, \cdots \}$ be any given countably infinite subset of the interval $[0, 3]$ and, for each positive integer $i$, let $y_i$ be an $\omega$-limit point of $x_i$.  Then $X'$ is a countably infinite set in the interval $[0, 2]$ and, for each positive integer $i$, since $g(x^*) = g(x)$, $y_i$ is also an $\omega$-limit point of $x_i^*$.  

Since $g(x)$ is a continuous bitransitive map (in fact, $g(x)$ is the tent map) from the interval $[0, 2]$ onto itself and the point $v = 2$ is a non-recurrent point of $g$ with $$\delta = \inf \big\{ |v - g^i(v)| : i = 1, 2, \cdots \big\} = 2,$$ we can apply Theorem 8 to the map $g(x)$ on the interval $[0, 2]$ with the same given $\rho$ to obtain a dense uncountable 2-scrambled set $Y$ in the interval $[0, 2]$ which is a union of countably infinitely many synchronously proximal Cantor sets $C_i, i \ge 1$ in the interval $[0, 2]$ such that each nonempty subinterval of the interval $[0, 2]$ contains {\it countably infinitely many} such Cantor sets $C_i$'s and the larger set $\widehat Y = \bigcup_{n=0}^\infty \, g^n(Y)$ is an {\it invariant} 2-scrambled set of $g$ in the interval $[0, 2]$ with the property that, for any points $x$ in $X'$ and $y$ in $\widehat Y$, we have 
$$
\liminf_{n \to \infty} \big|g^n(x) - g^n(y)\big| = 0 \quad \text{and} \quad \limsup_{n \to \infty} \big|g^n(x) - g^n(y)\big| \ge \rho.
$$
\indent Now let $\mathcal C = \{ C_i : i = 1, 2, \cdots \}$ and let $U_1, U_2, \cdots$ be an enumeration of all open intervals in the interval $[1, 2]$ with rational endpoints.  Then each $U_i$ contains countably infinitely many members of $\mathcal C$.  It is clear that the (symmetric with respect to the point $x = 2$) open intervals $U_1^*, U_2^*, \cdots$ are an enumeration of all open intervals in the interval $[2, 3]$ with rational endpoints.

Let $C_{i_1}$ be a Cantor set in $\mathcal C$ which lies in $U_1$.  Then $U_2 \setminus C_{i_1}$ is a nonempty open set.  Let $C_{i_2}$ be a Cantor set in $\mathcal C \setminus \{ C_{i_1} \}$ which lies in $U_2$.  Then $U_3 \setminus \big(C_{i_1} \cup C_{i_2}\big)$ is a nonempty open set.  Let $C_{i_3}$ be a Cantor set in $U_3 \setminus \big(C_{i_1} \cup C_{i_2}\big)$.  Inductively, after choosing $C_{i_k}, 1 \le k \le n-1$, we let $C_{i_n}$ be a Cantor set in $\mathcal C \setminus \{ C_{i_k} : 1 \le k \le n-1 \}$ which lies in $U_n$.  It is clear that the union $\bigcup_{k \ge 1} C_{i_k}$is dense in the interval $[1, 2]$ and, since eack $C_{i_k}^*$ is contained in $U_k^*$, the union $\bigcup_{k \ge 1} C_{i_k}^*$ is dense in the interval $[2, 3]$.  Furthermore, since each $U_k$ contains countably infinitely many $U_{i_j}$ $(\supset C_{i_j})$, each $U_k^*$ contains countable infinitely many $C_{i_j}^*$ and, since $g(C_{i_k}^*) = g(C_{i_k})$, each $C_{i_k}^*$ is a synchronously proximal Cantor set in the interval $[2, 3]$.  Let \\
[5pt]\centerline{$\mathcal Y = \bigcup_{k \ge 1} C_{i_k}^* \,\,\, \bigcup \,\,\, \big( \bigcup_{n \ge 1} C_n \setminus \bigcup_{k \ge 1} C_{i_k} \big)$.} \\
[5pt]Then, $g(\mathcal Y) = g(Y)$ and, since $Y$ is a dense uncountable 2-scrambled set of $g$ in the interval $[0, 2]$, we obtain that the set $\mathcal Y$ is a {\it dense} uncountable 2-scrambled set of $g$ in the interval $[0, 3]$ which is a union of countably infinitely many synchronously proximal Cantor sets and the larger set $\widehat {\mathcal Y} = \bigcup_{i \ge 0} \, f^i(\mathcal Y)$ is a dense uncountable {\it invariant} 2-scrambled set of $g$ in $[0, 3]$ with the property that, for any points $x$ in $X$ and $y$ in $\widehat {\mathcal Y}$, we have 
\[
\liminf_{n \to \infty} \big|g^n(x) - g^n(y)\big| = 0 \quad \text{and} \quad \limsup_{n \to \infty} |g^n(x) - g^n(y)| \ge \rho.\qedhere \]
\end{proof}
\indent In the above Theorem 8, we obtain dense scrambled sets which are unions of countably infinitely many Cantor sets with zero Lebesgue measure.  By the following result of Oxtoby and Ulam {[\bf{25}}, Theorem 9] \big(see also {\bf\cite{go}}\big), there exists a homeomorphism $h$ of $(-\infty, \infty)$ onto itself which takes these scrambled sets to sets of full Lebesgue measure.  Therefore, these maps are topologically conjugate to maps with scrambled sets of full Lebesgue measure (cf. {\bf\cite{bab}}).  We omit the details.  

\noindent
{\bf Theorem 10} {\rm (Oxtoby \& Ulam 1941}{\bf\cite{ox}}).
{\it Let $B$ be any subset of the $n$-dimensional euclidean space $\mathcal R^n$.  In order that there exists an automorphism $h$ of $\mathcal R^n$ (i.e., a homeomorphism from $\mathcal R^n$ onto itself) such that $h(B)$ has Lebesgue measure zero it is necessary and sufficient that the complement of $B$ contains a sequence of perfect sets whose union is dense in $\mathcal R^n$.}

In {\bf\cite{bab}}, Babilonov\'a shows that each continuous bitransitive map on the unit interval $I = [0, 1]$ is topologically conjugate to a continuous bitransitive map $g: I \longrightarrow I$ which is extremely chaotic almost everywhere.  Since every homeomorphism from an interval to itself must be strictly monotonic, by combining Theorems 8 $\&$ 10, we extent the result of Babilonov\'a from any compact interval to any interval in the real line.

\noindent
{\bf Corollary 11.}
{\it Let $I$ denote any interval in the real line.  Let $f$ be a continuous bitransitive map from $I$ onto itself. Then $f$ is topologically conjugate to a map in $C(I, I)$ which is extremely chaotic almost everywhere.}

By applying Theorem 10 to Examples 5 $\&$ 6 above, we obtain the following result pertaining to the existence of {\it invariant} extremely scrambled sets of transitive points of full Lebesgue measure.

\noindent
{\bf Corollary 12.}
{\it Let $I$ denote any interval in the real line.  Then there exists a continuous bitransitive map from $I$ onto itself which has an invariant extremely scrambled set of totally transitive points of full Lebesgue measure.} 

\noindent
{\bf Acknowledgement}\\
This work was partially supported by the Ministry of Science and Technology of Taiwan.

%\bibliographystyle{plain}
%\bibliography{duref}

\end{document}